\documentclass{article}

\usepackage{PRIMEarxiv}

\usepackage[utf8]{inputenc} % allow utf-8 input
\usepackage[T1]{fontenc}    % use 8-bit T1 fonts
\usepackage{hyperref}       % hyperlinks
\usepackage{url}            % simple URL typesetting
\usepackage{booktabs}       % professional-quality tables
\usepackage{amsfonts}       % blackboard math symbols
\usepackage{nicefrac}       % compact symbols for 1/2, etc.
\usepackage{microtype}      % microtypography
\usepackage{lipsum}
\usepackage{fancyhdr}       % header
\usepackage{graphicx}       % graphics
\graphicspath{{media/}}     % organize your images and other figures under media/ folder

\usepackage{amsmath}
\usepackage{subcaption}
\usepackage{amssymb}
\usepackage[style=authoryear, maxcitenames=1]{biblatex}
\usepackage{graphicx}
\usepackage{subcaption}
\usepackage{booktabs} % For formal tables
\usepackage{array}    % For table wrapping
\usepackage{algorithm}
\usepackage{algorithmic}
\usepackage{longtable}
\usepackage{todonotes}
\usepackage{lscape}
\usepackage{hyperref}
\usepackage{cleveref}
\usepackage{tikz}
\usetikzlibrary{decorations.pathreplacing, calligraphy}
\usepackage{titlesec}
\usepackage{titlesec}

\titleformat{\paragraph}[block]
  {\normalfont\normalsize\bfseries}{\theparagraph}{1em}{}
\titlespacing*{\paragraph}{0pt}{1\baselineskip}{0.5\baselineskip}
\title{Decentralized Modeling of Vehicular Maneuvers and Interactions at Urban Junctions
}

\author{
  Saeed Rahmani \\
  Department of Transport and Planning \\
  Delft University of Technology \\
  Delft, NL\\
  \texttt{s.rahmani@tudelft.nl} \\
   \And
  Simeon C. Calvert \\
  Department of Transport and Planning \\
  Delft University of Technology \\
  Delft, NL\\
  \texttt{s.c.calvert@tudelft.nl} \\
     \And
  Bart van Arem \\
  Department of Transport and Planning \\
  Delft University of Technology \\
  Delft, NL\\
  \texttt{b.vanarem@tudelft.nl}\\
}
\begin{document}
\maketitle

\begin{abstract}
Modeling and evaluation of automated vehicles (AVs) in mixed-autonomy traffic is essential prior to their safe and efficient deployment. This is especially important at urban junctions where complex multi-agent interactions occur. Current approaches for modeling vehicular maneuvers and interactions at urban junctions have limitations in formulating non-cooperative interactions and vehicle dynamics within a unified mathematical framework. Previous studies either assume predefined paths or rely on cooperation and central controllability, limiting their realism and applicability in mixed-autonomy traffic. This paper addresses these limitations by proposing a modeling framework for trajectory planning and decentralized vehicular control at urban junctions. The framework employs a bi-level structure where the upper level generates kinematically feasible reference trajectories using an efficient graph search algorithm with a custom heuristic function, while the lower level employs a predictive controller for trajectory tracking and optimization. Unlike existing approaches, our framework does not require central controllability or knowledge sharing among vehicles. The vehicle kinematics are explicitly incorporated at both levels, and acceleration and steering angle are used as control variables. This intuitive formulation facilitates analysis of traffic efficiency, environmental impacts, and motion comfort. The framework's decentralized structure accommodates operational and stochastic elements, such as vehicles' detection range, perception uncertainties, and reaction delay, making the model suitable for safety analysis. Numerical and simulation experiments across diverse scenarios demonstrate the framework's capability in modeling accurate and realistic vehicular maneuvers and interactions at various urban junctions, including unsignalized intersections and roundabouts.
\end{abstract}

% keywords can be removed

%% main text
\section{Introduction}
\label{sec:introduction}
The advent of partially and fully automated vehicles (AVs) promises to substantially transform transportation systems by enhancing safety, efficiency, and accessibility (\cite{dannemiller2023autonomous, tafidis2022safety}). However, the widespread adoption of AVs fundamentally depends on their ability to safely and efficiently navigate complex urban environments (\cite{shetty2021safety}). In particular, urban junctions, such as unsignalized intersections and roundabouts, stand out due to their intricate multi-directional interactions and ambiguous lane-following behavior, which represent a particularly challenging environment for automated vehicles (\cite{zhang2023platoon, chen2023nearly}). This complexity is further amplified in mixed-autonomy settings, where AVs must coexist and interact effectively with human-driven vehicles, whose behaviors are inherently heterogeneous and harder to predict. 

Evaluating whether AVs can safely and efficiently navigate such complex environments requires analytical tools and models that can accurately emulate individual vehicle driving behaviors, as well as emergent behaviors arising from multi-agent interactions. Accordingly, an effective modeling approach should satisfy several key requirements, including the accurate representation of vehicle dynamics to guarantee the kinematic feasibility of generated maneuvers; operation under decentralized control paradigms that represent the autonomous decision-making characteristics of mixed-automated-human traffic; and computational tractability to support extensive verification and optimization across diverse scenarios. Beyond these foundational requirements, a realistic modeling framework should accommodate various characteristics and algorithmic components of AVs, such as explicit representation of sensor range limitations and trajectory prediction algorithms that fundamentally shape AVs' behavior and decision-making in real-world deployments.

Current methods fall short in meeting these requirements when it comes to modeling urban junctions, such as unsignalized intersections and roundabouts. Traditional lane-based traffic flow models oversimplify the complex dynamics at intersections by reducing interactions to one-dimensional, leader-follower, and relying on predefined paths (\cite{alemdar2021interdisciplinary, arafat2021benefits, rahmani2023bi}). Moreover, many of these frameworks are inherently collision-free by formulation, which prevents their use for safety and accident analysis (\cite{zhang2024car}). Other approaches like Cellular Automata and Social Force models, while capable of representing non-lane-based flow, struggle to accurately capture continuous vehicle kinematics and non-holonomic motion constraints (\cite{singh2020cellular,zhao2023microscopic}). While more recent optimal control techniques are promising, their existing implementations frequently depend on the assumption of full cooperation or centralized control of all vehicles (\cite{zhao2023microscopic, yu2021model}). This inherent assumption makes them not applicable for modeling the decentralized decision-making that characterizes mixed-autonomy traffic. In addition, jointly optimizing trajectories across multiple interacting agents makes the resulting problem prohibitively expensive for large-scale, real-time simulations. More recently, machine learning methods have been explored for modeling multi-agent scenarios (\cite{tang2019towards, di2021survey}). However, such methods often lack interpretability, are difficult to configure for testing targeted behavioral changes, and require large datasets to learn realistic, generalizable behaviors, which are currently unavailable (\cite{di2021survey}).

Motivated by these limitations, we develop a bi-level modeling framework that addresses the key challenges of modeling vehicular maneuvers and decentralized interactions at urban junctions in mixed-autonomy traffic. The framework operates in a fully decentralized manner without imposing assumptions of cooperation or knowledge sharing among vehicles. This architecture more accurately captures the autonomous decision-making processes inherent in mixed-automated-human traffic and enables safety analysis through emergent phenomena. The framework explicitly incorporates vehicle dynamics through a kinematic vehicle model and receding horizon control scheme, with acceleration and steering angle serving as control variables. These intuitive and interpretable control inputs make the framework well-suited for diverse applications, including traffic efficiency, environmental impact, and motion comfort analysis. The bi-level architecture of the framework is designed to achieve computational efficiency through the separation of planning and control horizons. This design is computationally efficient and supports real-time simulation with complexity that scales only linearly with the number of vehicles. This property addresses the limitations of computationally expensive optimization-based methods that rely on offline computation and scale exponentially. Last but not least, the bi-level and modular structure of the proposed formulation accommodates key characteristics and behavioral factors such as detection range, planning horizon, prediction uncertainties, and delayed responses from drivers or the vehicles. This makes it particularly suitable for safety impact assessment and accident analysis. To facilitate reproducibility and broader adoption, we make the framework publicly available as open-source software\footnote{\url{https://github.com/SaeedRahmani/AV-Simulation-at-Intersections}}. Taken together, the theoretical contributions and practical applicability of this framework represent a significant advancement in vehicular maneuver and interaction modeling in the era of mixed autonomy.

The remainder of this paper is structured as follows: Section \ref{sec:literature} provides a review of existing approaches for modeling vehicular maneuvers and interactions at intersections. Section \ref{sec:methodology} presents the mathematical formulation of the proposed framework. Section \ref{sec:experiments} describes the numerical and simulation experiments designed to validate the framework's theoretical properties and computational performance. Finally, Section \ref{sec:conclusion} discusses theoretical implications and outlines directions for future research.

\section{Related Works}\label{sec:literature}
Current approaches to modeling vehicular movements at intersections can be broadly categorized into four categories (\cite{zhao2020two}): lane-based traffic flow models, cellular automata (CA), social force models, and optimal control techniques.

\subsection{Lane-based Traffic Flow Models}
Lane-based traffic flow models, initially developed to understand and predict vehicular traffic dynamics on roadways, primarily encompass car-following and gap acceptance models. The main body of research in this domain has been focused on calibrating car-following models at intersections, with a special focus on signalized intersections.
\textcite{mathew2010calibration} and \textcite{alemdar2021interdisciplinary} calibrated car-following models in VISSIM for signalized intersections. However, these models were limited to longitudinal interactions, neglecting crucial lateral dynamics such as turning radius and axle configurations. Moreover, they often assume predefined fixed paths for vehicles, which oversimplifies the complex interactions and maneuvers at intersections.
For unsignalized intersections, \textcite{pollatschek2002decision} developed an analytical decision model for gap acceptance behavior on minor roads. However, this approach is confined to modeling driver decisions about gap acceptance and rejection, without addressing the actual vehicular movements and 2D interactions. Similarly, \textcite{arafat2021benefits} calibrated a gap acceptance model in VISSIM to study left-turning assistant systems, but vehicular movement was still represented using a one-dimensional car-following model.
While these studies provide valuable insights into high-level intersection operations, they rely heavily on predefined paths and primarily consider unidirectional movements. The interaction mechanisms in these models are typically limited to bi-vehicle (follower-leader) interactions, which fail to capture the multi-vehicle, multi-directional interactions common at intersections. Furthermore, these traffic flow models are inherently collision-free, limiting their applicability in safety analysis.

\subsection{Cellular Automata Models}
Cellular automata (CA) models offer an alternative approach to modeling vehicle movements at intersections. These models utilize a discrete lattice of cells, where each cell can be in a certain state, and vehicles occupy one or multiple cells at a time.
\textcite{foulaadvand2007vehicular} adapted the Nagel-Schreckenberg CA model to characterize conflicting vehicle traffic flow at road junctions. However, their model was limited to one-dimensional vehicle movements. \cite{wang2021multi} and \cite{chai2015fuzzy} employed CA models to evaluate intersection traffic flow performance but focused primarily on modeling car-following and lane-changing behavior before vehicles enter the intersection. Recent advancements have extended CA models to optimize intersection flow (\cite{zhu2018modeling, cruz2019automated}) and study traffic accidents (\cite{marzoug2022modeling}) at intersections. However, these studies have been largely confined to signalized intersections. For a comprehensive overview of CA models in intersection studies, readers are directed to the survey by \textcite{singh2020cellular}. 

While CA models can represent non-lane-based traffic flows and capture more detailed vehicle interactions compared to traditional traffic flow models, their application has been primarily limited to signalized intersections. Moreover, the fixed-size, grid-based nature of CA models constrains their ability to accurately represent vehicle kinematics and spatial variations in movement. Increasing the resolution of cells to capture more accurate vehicular movements significantly increases computational costs, limiting their practicality for real-time simulations (\cite{singh2020cellular}).

\subsection{Social Force Models}
Social force models describe interactions among agents as a series of forces influencing their behavior. While primarily applied to modeling pedestrian and cyclist interactions (\cite{chen2018social, shrivas2024modified, golchoubian2023pedestrian}), these models have also been adapted for vehicular movements at intersections.
\textcite{ma2017two} developed a social force model for vehicles at intersections, but their assumptions were based on lane-based environments, limiting the model's applicability to true 2D environments. \textcite{yang2018model} used social force theory to model vehicle movements around work zones within intersections, but did not consider interactions with crossing vehicles.

While social force models offer a more comprehensive consideration of interactions, they have several shortcomings in modeling vehicular traffic. The resulting trajectories are primarily determined by the combination of several "forces" rather than explicitly modeling driver decision-making or vehicle dynamics (\cite{zhao2023microscopic}). This reduces their interpretability and applicability to safety and efficiency impact assessments and limits their capability to explicitly consider the non-holonomic constraints of vehicles. Additionally, the high complexity of these models challenges their scalability to dense traffic situations (\cite{ma2017two}), and their numerous parameters make calibration demanding.
\subsection{Machine Learning Methods}
Machine learning approaches have gained significant attention for modeling vehicular maneuvers and interactions at intersections, offering data-driven alternatives to traditional rule-based methods. Recent studies have explored deep reinforcement learning for autonomous intersection navigation (\cite{li2022continuous,wu2024recent}), with multi-agent reinforcement learning frameworks showing promise for decentralized coordination without central control (\cite{zhao2024survey}). Neural network architectures, particularly LSTM encoder-decoder models and attention-based transformers, have demonstrated effectiveness in trajectory prediction at complex urban intersections (\cite{chen2023milestones, hou2019interactive}). Generative adversarial networks have also been applied to traffic simulation and behavior modeling (\cite{chen2022combining}). However, these machine learning methods face several critical limitations that restrict their applicability in the scope of intersection modeling: they often lack interpretability and transparency in decision-making processes, require extensive training datasets that may not be readily available for diverse intersection scenarios, struggle with generalization to unseen traffic conditions or geographic regions, and pose significant challenges for safety analysis and verification due to their "black-box" nature (\cite{di2021survey, haydari2020deep}). These limitations motivate the development of more interpretable and mathematically grounded approaches, such as the optimal control framework proposed in this study.

\subsection{Optimal Control Techniques}
Optimal control theory has emerged as a popular technique for modeling vehicle movements at intersections (\cite{wang2023optimal, yu2019managing}). This approach offers great flexibility in capturing multi-dimensional interactions and movements while considering vehicle dynamics. Many optimal control-based models for intersections focus on coordinating multiple vehicles using a central controller or assuming connectivity and cooperation among vehicles (\cite{wu2022intersection, yu2021model, pan2022convex, chen2021mixed, malikopoulos2018decentralized}). However, the assumption of full cooperation, connectivity, or controllability among vehicles limits their applicability in mixed-traffic environments where not all vehicles are connected or autonomous. In the robotics community, optimal control algorithms have been developed for single autonomous vehicles navigating among human-driven vehicles (\cite{yu2021model, stano2023model}). These trajectory optimization and path tracking methods are designed for real-world operation, processing vast amounts of sensor data and capturing high levels of detail, such as vehicle mechanical constraints, tire forces, etc. While necessary for real-world applications, this level of detail can be computationally expensive and may not be suitable for analytical studies and real-time microsimulation studies, especially when dealing with multiple vehicles (\cite{samak2021control, yu2021model}).

Accordingly, efforts have been made to reduce the computational load while preserving key characteristics of vehicle dynamics and interactions. \textcite{zhan2017spatially, liu2017speed} simplified the trajectory optimization problem to a speed optimization problem by assuming pre-planned paths and predefined passing orders for vehicles.  \textcite{shi2023trajectory} extended this work by removing the assumption of passing priority among the vehicles. However, they still assumed a pre-defined trajectory for the vehicles, which simplifies the problem into a one-dimensional speed optimization problem and ignores the non-holonomic relationships between the vehicle's speed and other control inputs such as steering angle. \textcite{zhao2020two} proposed a two-dimensional model for describing single-vehicle maneuvers within an intersection, later extending it to incorporate interactions among vehicles (\cite{zhao2023microscopic}). While successful in planning 2D trajectories, this model lacks some important features and is developed based on some relatively strong assumptions. Firstly, the proposed model by \textcite{zhao2023microscopic} assumes full cooperation among vehicles and that they share their utility functions. This limits its applicability to mixed traffic environments. Also, this assumption implies that there will be no collision among vehicles, which makes the proposed model not suitable for safety analysis. 

The literature review reveals that existing approaches for modeling vehicular maneuvers at urban junctions have significant limitations. Lane-based models oversimplify intersection dynamics through one-dimensional representations and predefined paths. Cellular automata and social force models lack precision in incorporating vehicle kinematics and non-holonomic constraints. Machine learning approaches suffer from interpretability issues and extensive data requirements that hinder safety-critical deployment. Current optimal control techniques predominantly assume full cooperation or centralized coordination, making them unsuitable for mixed-autonomy traffic. Critically, no existing framework successfully combines kinematically feasible trajectory planning with decentralized collision avoidance while maintaining computational tractability for real-time applications, motivating the bi-level optimal control framework proposed in this study.

\section{Methodology}\label{sec:methodology}
\begin{figure}[b]
    \centering
    \includegraphics[width=0.4\linewidth]{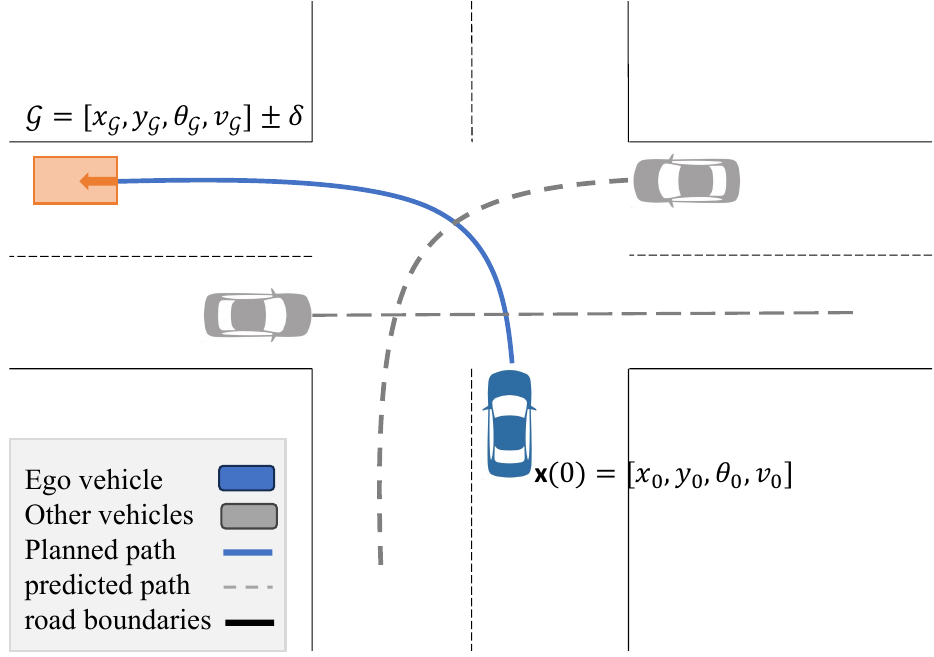}
    \caption{Problem visualization for an imaginary scenario including the ego vehicle and two other vehicles}
    \label{fig:problem_scenario}
\end{figure}

We consider an urban junction scenario involving multiple vehicles operating under decentralized control without inter-vehicle communication or centralized coordination (Figure~\ref{fig:problem_scenario}). The ego vehicle seeks to navigate from its initial state $x(0)$ to a goal region $\mathcal{G}$ while respecting kinematic constraints and avoiding collisions. We formulate the problem as finite-horizon optimal control on \(\mathrm{SE}(2) \times \mathbb{R}\), where each vehicle $i$ is characterized by its state vector $\mathbf{x}_i(t) = [x_i(t), y_i(t), \theta_i(t), v_i(t)]^{\top}$. This state vector represents the planar position $(x_i, y_i)$, heading angle $\theta_i$, and velocity $v_i$ at time $t$. The control inputs are longitudinal acceleration  \(a\) and steering angle \(\delta\). 

While this formulation is mathematically complete, solving the resulting problem is computationally expensive in a multi-vehicle setup once interactions and collision constraints are included. To obtain a computationally tractable solution, we propose a bi-level modeling approach that decomposes the problem into two tractable subproblems: coordinated spatial path planning and time-parameterized motion control. As illustrated in Fig.~\ref{fig:framework}, the proposed modeling framework comprises a global planner that generates a kinematically feasible reference trajectory considering the intersection geometry, traffic rules, and user-defined criteria, and a local motion controller that tracks the reference trajectory while avoiding dynamic obstacles by applying a real-time collision avoidance strategy. The two components interact continuously through a feedback mechanism, enabling adaptive replanning when trajectory deviations occur or new obstacles are detected. The mathematical formulation of this bi-level framework is as follows:
\begin{figure}[t]
    \centering    \includegraphics[width=1.0\textwidth]{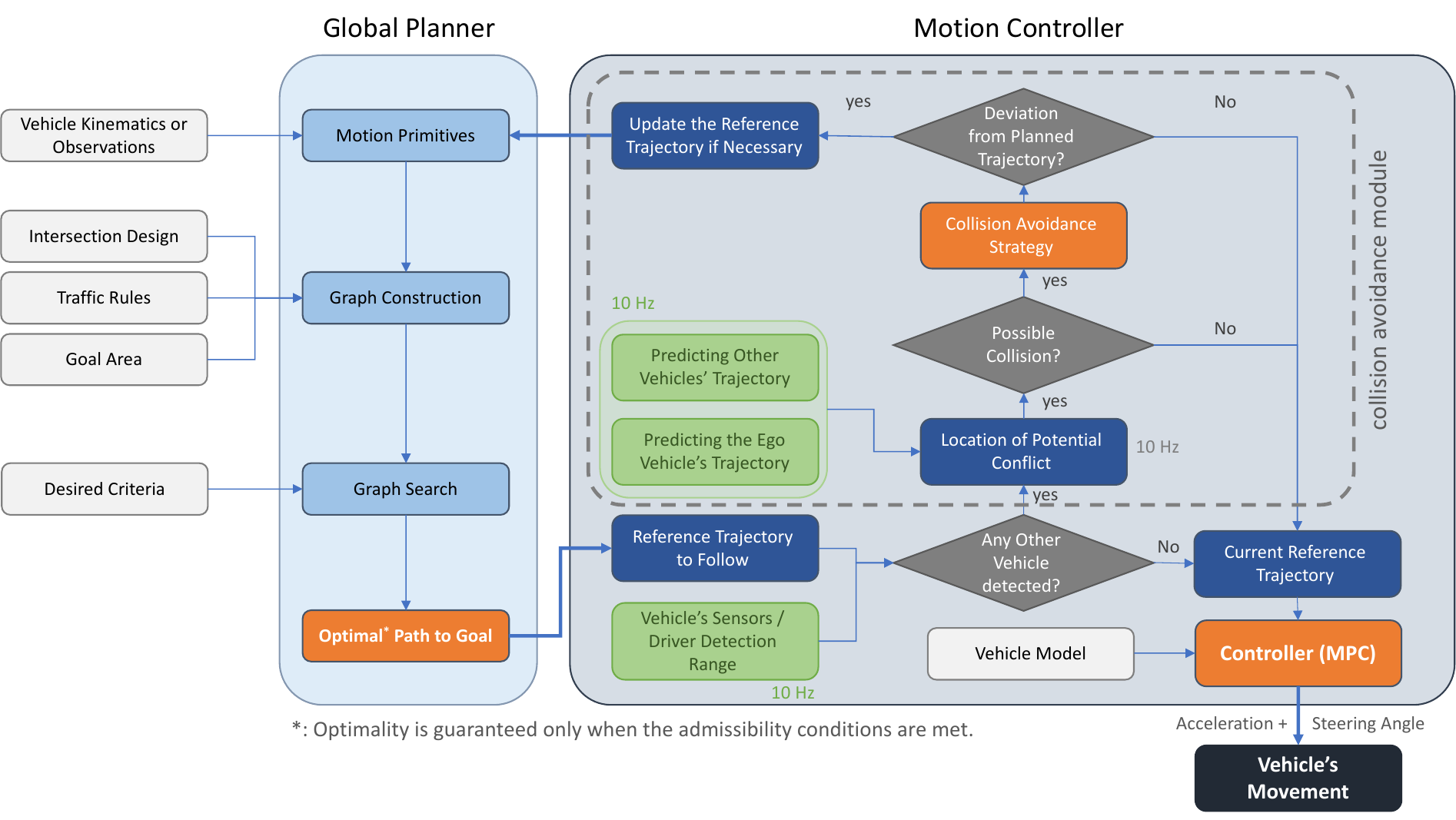}
    \caption{The proposed framework for modeling the vehicular maneuvers and interactions at urban junctions. The global planner generates optimal reference trajectories through graph search, while the local motion controller executes real-time collision avoidance and trajectory tracking with continuous feedback between levels.}
    \label{fig:framework}
\end{figure}

At the higher level, the planner synthesizes a reference trajectory $\tau^*$, including geometric waypoints and associated reference speeds, by searching a lattice of pre-computed motion primitives:
\begin{equation}
    \tau^* = \arg\min_{\{\mathbf{q}_i\}} \sum_{i=1}^{N_p} J_{\text{global}}(\mathbf{q}_i, v_i^{\mathrm{ref}})
    \quad \text{s.t.}\quad \mathbf{q}_i \in \mathcal{Q}_{\text{feasible}}, \; \tau^* \in \mathcal{T}_{\text{safe}} \cap \mathcal{T}_{\text{legal}}
\end{equation}
where 
\begin{equation}
    \tau^* = \bigl\{(\mathbf{q}_i, v_i^{\mathrm{ref}})\mid \mathbf{q}_i=(x_i, y_i, \theta_i)\in\mathrm{SE}(2),\; v_i^{\mathrm{ref}}\in\mathbb{R}\bigr\}
\end{equation}
In this formulation, the planner optimises only the geometric component $\mathbf{q}_i$ that represents the waypoints connected by motion primitives. $v_i^{\mathrm{ref}}$ are reference speeds assigned to each waypoint; $x_i$ and $y_i$ denote the global coordinates of the $i$-th waypoint; $\theta_i$ is its yaw angle; and $N_p$ is the planning horizon. $\mathcal{Q}_{\text{feasible}}$ represents the space of kinematically admissible waypoint configurations respecting vehicle motion limits, $\mathcal{T}_{\text{safe}}$ denotes the collision-free trajectory space excluding static obstacles, and $\mathcal{T}_{\text{legal}}$ encompasses the trajectory space compliant with traffic regulations.

At the subsequent operational level, a motion controller utilizes this spatially optimised reference trajectory to determine optimal control inputs while accounting for dynamic obstacles and vehicular interactions. This optimization occurs within the complete $\mathrm{SE}(2)\times\mathbb{R}$ state space at each temporal instance $t_k$:
\begin{equation}
    \min_{\mathbf{u}(\cdot)} \int_{t_k}^{t_k+N_c} \ell\bigl(\mathbf{x}(t),\mathbf{u}(t),\mathbf{x}_{\text{ref}}(t)\bigr)\,\mathrm{d}t
    \quad\text{s.t.}\quad
    \mathbf{\dot x} = f(\mathbf{x},\mathbf{u}),\;
    \|\mathbf{u}(t)\|\le \mathbf{u}_{\max},
\end{equation}
where
\begin{equation}
    \mathbf{x}(t)=\bigl[x(t),\,y(t),\,\theta(t),\,v(t)\bigr]^{\!\top}\in\mathrm{SE}(2)\times\mathbb{R},
    \qquad
    \mathbf{u}(t)=\bigl[a(t),\,\delta(t)\bigr]^{\!\top}.
\end{equation}
In these equations, $\mathbf{x}_{\text{ref}}(t)$ is the reference state at time $t$ obtained from $\tau^*$. The details about how the time parametrized $\mathbf{x}_{\text{ref}}(t)$ is obtained from $\tau$, is presented in Section~\ref{sec:motion controller}. $\ell$ is a cost function defined over the vehicle's state $\mathbf{x}$, reference state $\mathbf{x}_{ref}$, and the control input $\mathbf{u}$; $f$ encodes the vehicle dynamics. $a(t)$ is the longitudinal acceleration, and $\delta(t)$ is the steering angle at time $t$. To ensure tractability, the control horizon $N_c$ is chosen significantly shorter than the planning horizon $N_p$. To compensate for this reduced horizon and maintain safety, a dedicated collision-avoidance module predicts the states of the surrounding agents over an extended detection horizon and enforces safety constraints within the controller (the dotted rectangle in Figure~\ref{fig:framework}). Should avoidance maneuvers induce trajectory deviations beyond defined tolerances, the reference trajectory is recomputed through closed-loop feedback before being fed to the controller. 

The framework architecture aligns conceptually with Michon’s hierarchical model of the driving task (\cite{michon1985critical}) and modular control architectures for automated vehicles (\cite{schwarting2018planning}) by separating strategic and tactical reasoning from operational execution while maintaining feedback across levels. This structural correspondence and the decentralized control architecture make the framework particularly suitable for modeling mixed-automated-human interactions and dynamics. 
% The bi-level design enhances both scalability and interpretability. Each vehicle independently plans and optimizes its motion, resulting in computational complexity that grows linearly with the number of vehicles.
Additionally, the hierarchical separation of strategic reasoning from operational execution makes the decision-making process transparent and traceable across levels. This modular structure provides flexibility for different system configurations through tunable parameters such as detection ranges, prediction horizons, and cost-function weights at each level. In the next sections, the two main modules of the proposed framework are discussed in detail, starting with the global planner.
% \begin{figure}[t]
%     \centering
%     \includegraphics[width=0.6\textwidth]{images/michon_framework.pdf}
%     \caption{Michon's Framework for driving tasks (\cite{michon1985critical})}
%     \label{fig:michon_framework}
% \end{figure}

\subsection{Global Planner}\label{global planner}
\begin{figure}[b]
    \centering    \includegraphics[width=1.0\textwidth]{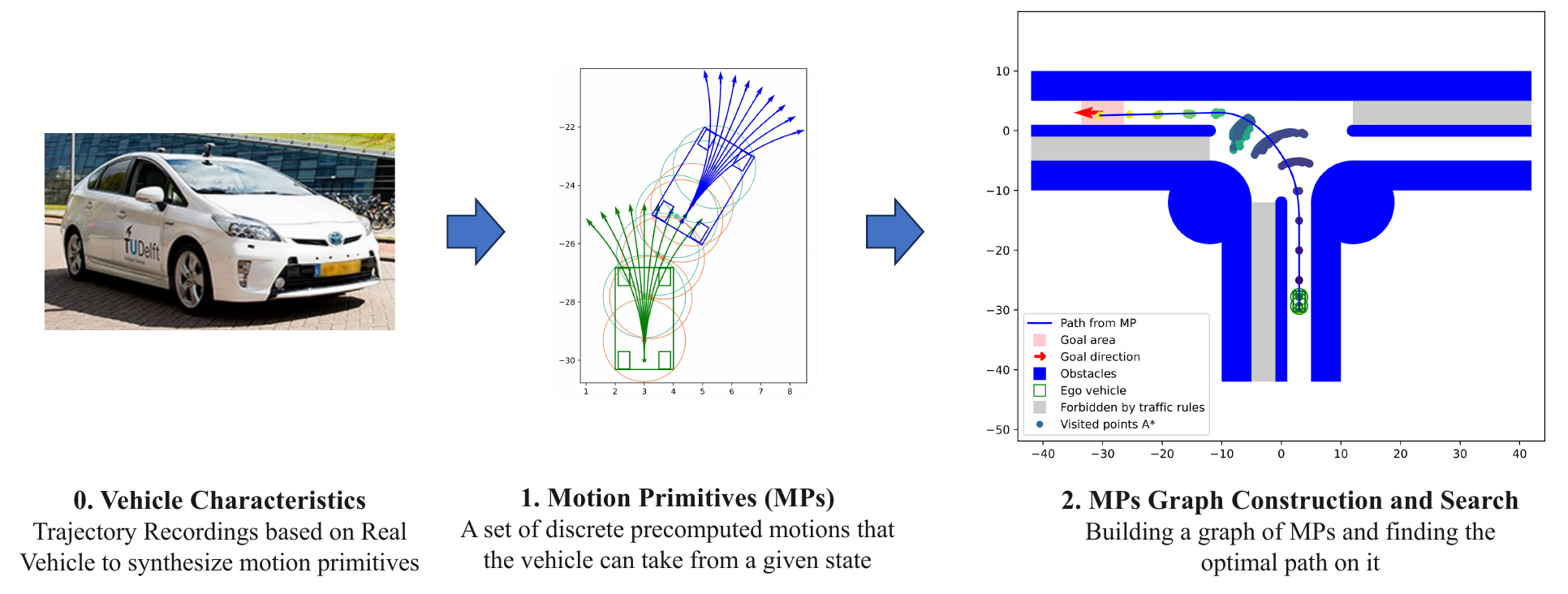}
    \caption{The components of global planner: 1. motion primitives, 2. graph construction and search.}
    \label{fig:Global_planner}
\end{figure}

The global planner architecture integrates two fundamental components: a motion primitives generator that synthesizes kinematically constrained vehicle maneuvers, and a unified graph construction and search algorithm that embeds these primitives into a traversable state space and computes the minimum-cost trajectory subject to user-defined optimization criteria. These components are depicted in Figure \ref{fig:Global_planner} and further elaborated in the following sections.

\subsubsection{Motion Primitives}\label{motion primitives}
Motion primitives constitute the fundamental discretisation units for our kinematically feasible trajectory-generation. We define these primitives as trajectory segments that explicitly satisfy the vehicle’s non-holonomic constraints. As a proof of concept, in this study, we formalise the motion-primitive generation through a state-space representation based on a high-fidelity nonlinear model derived from TU Delft’s autonomous vehicle platform (\cite{IVG_Facilities, Spahn_gym_envs_urdf}). The system dynamics are governed by
\begin{equation}
    \dot x(t)=v(t)\cos\bigl(\theta(t)\bigr), \quad
    \dot y(t)=v(t)\sin\bigl(\theta(t)\bigr), \quad
    \dot\theta(t)=\frac{v(t)}{L}\tan\bigl(\delta(t)\bigr),
    \label{eq:state_equations}
\end{equation}
where $L$ is the wheelbase and $\delta(t)$ the steering angle. For discretisation, the admissible steering domain $[-\delta_{\max}, +\delta_{\max}]$ is partitioned into $n\in\mathbb{N}$ discrete values $\{\delta_i\}_{i=1}^{n}$. Each motion primitive is then derived through numerical integration of the system dynamics over the interval $\Delta t$ with constant $\delta_i$:
\begin{equation}
    \mathbf{q}_i(t+\Delta t)=\mathbf{q}_i(t)+\int_{t}^{t+\Delta t}f_\mathbf{q}\bigl(\mathbf{q}_i(t'),\delta_i\bigr)\,dt',
    \label{eq:motion_primitives_combined}
\end{equation}
where $\mathbf{q}_i(t)=[\,x_i(t),y_i(t),\theta_i(t)\,]^{\top}$,
$f_\mathbf{q}=[\,\dot x_i,\dot y_i,\dot\theta_i\,]^{\top}$ follows from Equation~\ref{eq:state_equations}, and $t'$ is the dummy integration variable within $[t,t+\Delta t]$.

\begin{figure}[t]
    \centering 
  \begin{subfigure}[b]{0.265\textwidth}
    \includegraphics[width=\textwidth]{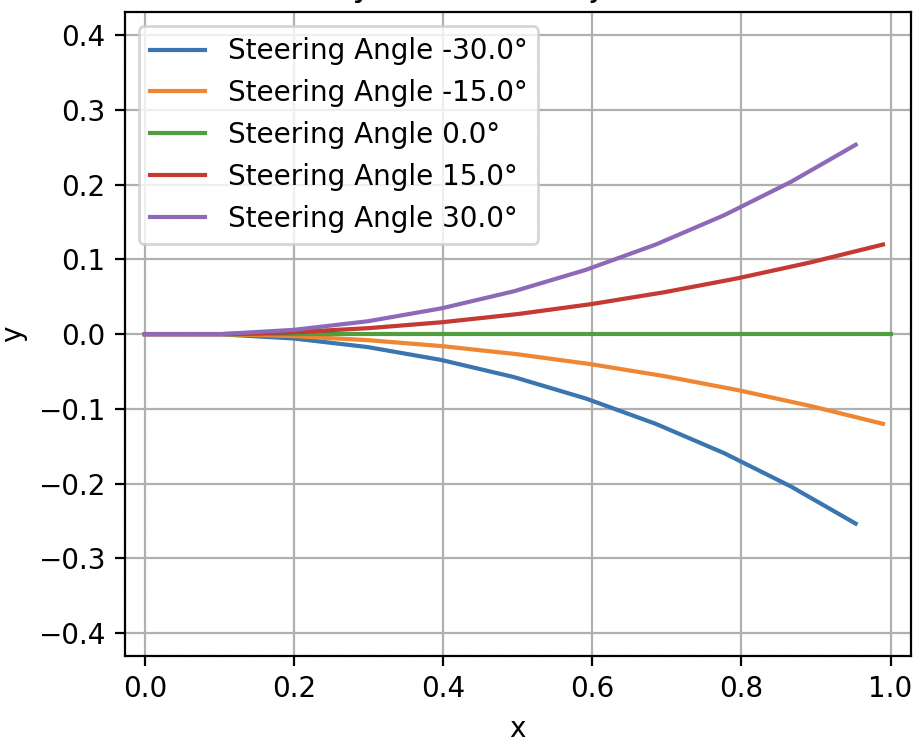}
    \caption{}
    \label{fig:motion_primitives_a}
  \end{subfigure}
  \begin{subfigure}[b]{0.3\textwidth}
    \includegraphics[width=\textwidth]{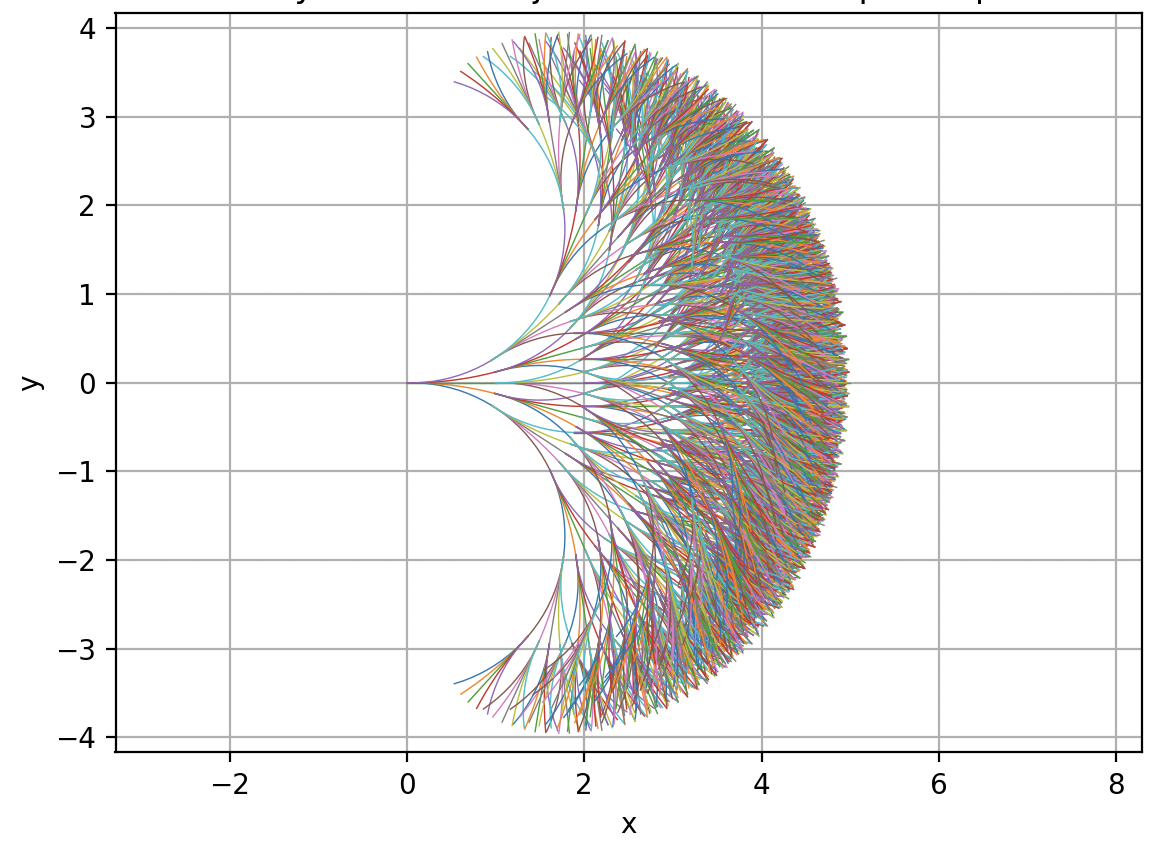}
    \caption{}
    \label{fig:motion_primitives_b}
  \end{subfigure}
  \begin{subfigure}[b]{0.275\textwidth}
    \includegraphics[width=\textwidth]{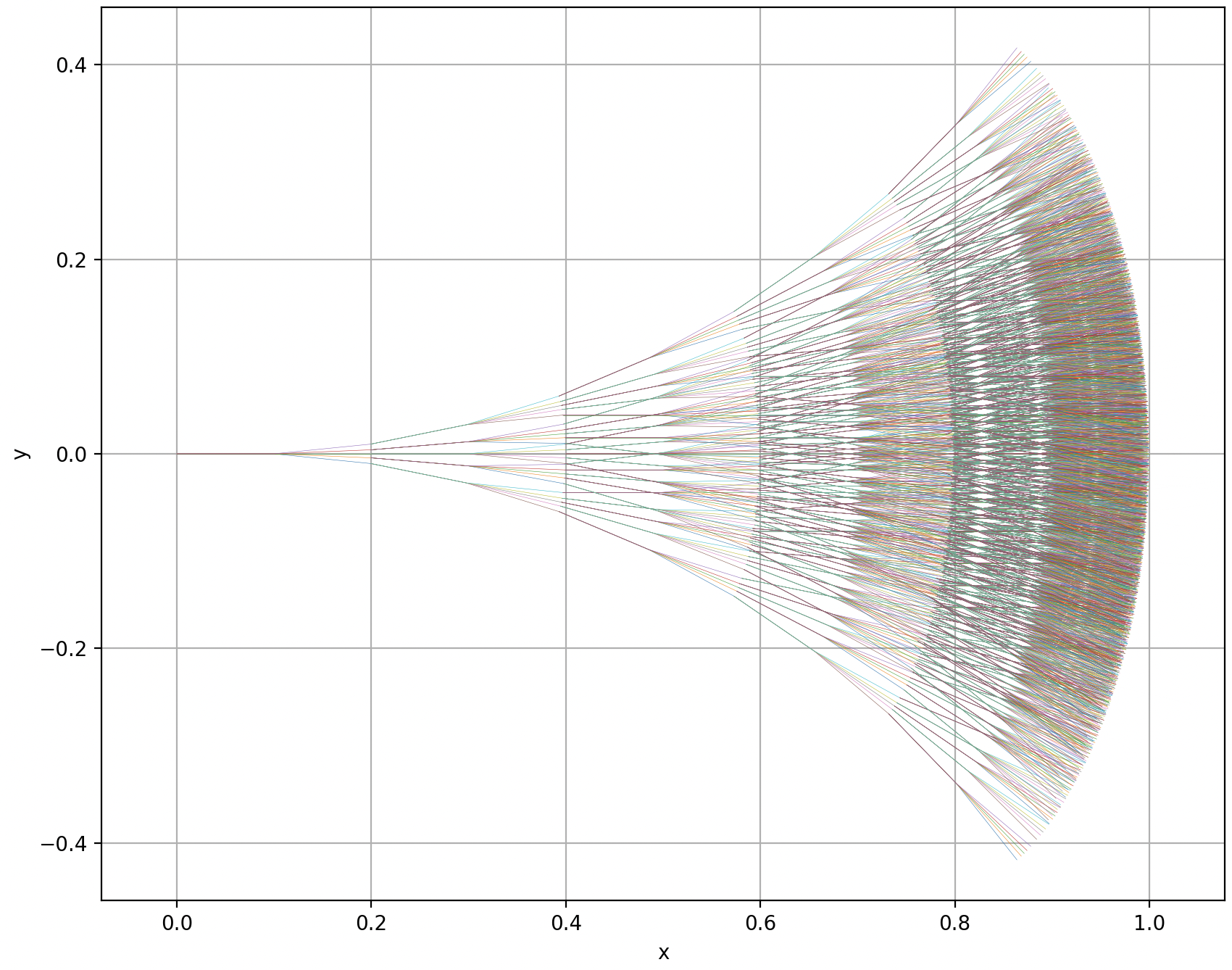}
    \caption{}
    \label{fig:motion_primitives_c}
  \end{subfigure}
  \caption{Visualization of motion primitives: (a) 5 discrete motion primitives with the length of 1 meter, (b) Sequence of 5 steps of motion primitives with the length of 1 meter, and (c) Sequence of 5 steps of motion primitives with the length of 0.2 meter.}
  \label{fig:motion primitives}
\end{figure}

Figure~\ref{fig:motion_primitives_a} visualises a set of five discrete primitives with $\delta_{\min}=-30^{\circ}$ and $\delta_{\max}=30^{\circ}$ at constant speed. The concatenation of these primitives over multiple time steps (Figure~\ref{fig:motion_primitives_b}) demonstrates how even a minimal set efficiently approximates the continuous reachability space. Figure~\ref{fig:motion_primitives_c} illustrates the effect of a reduced primitive length (0.2\,m), which enables higher spatial resolution but expands the search space. The chosen primitive length therefore represents a trade-off between spatial fidelity and computational complexity that can be calibrated to application-specific requirements.
\subsubsection{Graph Construction and Search}\label{sec:graph_construction_and_search}
Following the generation of motion primitives, we formulate the trajectory-planning problem as a combinatorial search over a directed graph. We construct a directed graph where vertices represent vehicle configurations and edges correspond to kinematically feasible transitions via motion primitives. Each configuration $\mathbf{q}\in\mathcal{Q}$ is the tuple $(x,y,\theta)$. Our approach incorporates several key algorithmic innovations that enhance computational efficiency and provide configurability through a novel adaptation of the A* search algorithm (Appendix A). Firstly, unlike standard graph-search algorithms that operate on pre-computed graphs, we employ a lazy graph representation with on-demand constrained expansion. Vertices and edges are generated progressively during the search process, while the search explores only collision-free configurations that comply with traffic rules. To achieve this, we define a neighbour function that identifies valid adjacent nodes on demand during graph construction. The neighbour function $\mathcal{N}:\mathcal{Q}\rightarrow 2^{\mathcal{Q}\times\mathbb{R}^+}$ is
\begin{equation}
    \mathcal{N}(\mathbf{q}) = \{(\mathbf{q}', c) \mid \mathbf{q}' = T(\mathbf{q}, m_k),\; c = g(\mathbf{q}, m_k),\; m_k \in \mathcal{M},\; \text{IsValid}(\mathbf{q}, m_k)\},
\end{equation}
where $\mathbf{q}' \in Q$ is the successor configuration obtained by applying motion primitive $m_k$ to $\mathbf{q}$, $c \in \mathbb{R}^+$ is the associated transition cost, $T: Q \times \mathcal{M} \rightarrow Q$ is a transition map that computes the resulting configuration after applying motion primitive $m_k$ to configuration $\mathbf{q}$, $g: Q \times \mathcal{M} \rightarrow \mathbb{R}^+$ is the cost function that computes the cost of applying motion primitive $m_k$ from configuration $\mathbf{q}$, $\mathcal{M} = \{m_1, \ldots, m_K\}$ is the set of motion primitives, and $\text{IsValid}(\mathbf{q}, m_k)$ is the constraint validation function. Constraint validation ensures collision-free transitions that also comply with traffic rules:
\begin{equation}
    \text{IsValid}(\mathbf{q}, m_k) \Leftrightarrow \{R_{\theta_q} \mathbf{p} + (x_q, y_q)^{\top} \mid \mathbf{p} \in \mathcal{P}_{m_k}\} \cap \Bigl(\bigcup_{i=1}^{N_{ob}} O_i \cup \bigcup_{j=1}^{N_{tr}} P_j\Bigr) = \emptyset,
\end{equation}
where $R_{\theta_q}$ is the planar rotation matrix using the heading angle $\theta_q$ from configuration $\mathbf{q}$, $\mathbf{p}$ represents individual collision-checking points, $(x_q, y_q)^{\top}$ are the position coordinates from configuration $\mathbf{q}$, $\mathcal{P}_{m_k}$ are the local collision-checking points of motion primitive $m_k$, $N_{ob}$ is the number of static obstacles, $O_i$ is the $i$-th static obstacle, $N_{tr}$ is the number of regions forbidden by traffic rules, and $P_j$ is the $j$-th region forbidden by traffic rules. This formulation reduces collision detection and traffic-rule compliance to efficient matrix operations while maintaining kinematic feasibility.

Secondly, we design a multi-criteria heuristic function, which extends beyond conventional Euclidean-distance metrics. The heuristic function $h$ is composed of a weighted formulation, simultaneously considering spatial proximity, orientation alignment and control-effort estimation:
\begin{equation}
h(\mathbf{q}) = w_{h,d} h_d(\mathbf{q}, \mathcal{G}) + w_{h,\theta} h_\theta(\mathbf{q}, \mathcal{G}) + w_{h,\varphi} h_\varphi(\mathbf{q})
\end{equation}
where $w_{h,d}$, $w_{h,\theta}$, and $w_{h,\varphi}$ are the weights for spatial proximity, orientation alignment, and control-effort estimation, respectively. Here, $h_d(\mathbf{q}, \mathcal{G})$ estimates the spatial distance from configuration $\mathbf{q}$ to the goal region $\mathcal{G}$, $h_\theta(\mathbf{q}, \mathcal{G})$ estimates the orientation alignment cost from configuration $\mathbf{q}$ toward the goal region $\mathcal{G}$ to ensure the search focuses on areas directed toward the goal, and $h_\varphi(\mathbf{q})$ estimates the control effort required for the configuration $\mathbf{q}$. This multi-criteria approach accelerates convergence by exploring targeted aras and avoiding branches that require harsh or unrealistic maneuvers to improve kinematic efficiency. 

Similar to the heuristic function, the cost function in the proposed algorithm balances multiple objectives through weighted component costs:
\begin{equation}
g(\mathbf{q},m_k) = w_{c,d} g_d(\mathbf{q},m_k) + w_{c,\varphi} g_\varphi(\mathbf{q},m_k) + w_{c,c} g_c(\mathbf{q},m_k)
\end{equation}
where $w_{c,d}$, $w_{c,\varphi}$, and $w_{c,c}$ are the weights for path length, steering effort, and obstacle clearance, respectively. Here, $g_d(\mathbf{q},m_k)$ quantifies the path length of motion primitive $m_k$ applied from configuration $\mathbf{q}$, $g_\varphi(\mathbf{q},m_k)$ quantifies the steering effort required for the transition from configuration $\mathbf{q}$ using motion primitive $m_k$, and $g_c(\mathbf{q},m_k)$ quantifies the obstacle clearance cost associated with the transition, penalizing trajectories that pass close to obstacles. The weights can be calibrated to prioritise shorter paths versus smoother trajectories, adapting the framework to different driving styles or vehicle characteristics.

The algorithm maintains a priority queue sorted by $f(n)=g(n)+h(n)$, exploring nodes with the lowest cost first. Search terminates when the goal function identifies a satisfactory solution:
\begin{equation}
\text{is\_goal}(q)=
   \begin{cases}
       \text{true} & \text{if } d\bigl({\mathbf{q}}_{xy},\mathcal{G}\bigr)\le\epsilon \;\text{ and }\; |\theta-\theta_{\mathcal{G}}|\le\delta_{\theta}\\[4pt]
       \text{false} & \text{otherwise}
   \end{cases}
\end{equation}
where $d(\mathbf{q}_{xy},\mathcal{G})$ is the distance from the planar position of $\mathbf{q}$ to the goal area $\mathcal{G}$, $\epsilon$ is the positional tolerance and $\delta_{\theta}$ is the allowable heading deviation. Optimality is guaranteed when the heuristic is admissible ($h(\mathbf{q}) \le h^*(\mathbf{q})$) and consistent ($h(\mathbf{q}) \le g(\mathbf{q},m_k) + h(\mathbf{q}')$ for any valid transition from $\mathbf{q}$ to $\mathbf{q}'$ via motion primitive $m_k$). Feasibility is ensured through motion primitives that respect kinematic constraints and explicit collision checking, ensuring both kinematic feasibility and collision-free paths. The algorithm is resolution-complete within the discretised configuration space.

Upon reaching the goal, the algorithm reconstructs the path by tracing predecessor relationships backward, yielding discrete configurations $[\mathbf{q}_0,\mathbf{q}_1,\ldots,\mathbf{q}_n]$. The complete trajectory is generated by concatenating the transformed motion primitives for each transition:
\begin{equation}
\mathcal{T} = \bigoplus_{i=0}^{n-1} \{R_{\theta_i}\mathbf{p} + (x_i,y_i)^{\top} \mid \mathbf{p} \in m_{k_i}\}
\end{equation}
where $\bigoplus$ denotes ordered concatenation, $m_{k_i}$ is the motion primitive used for the transition from configuration $\mathbf{q}_i$ to $\mathbf{q}_{i+1}$, $R_{\theta_i}$ is the planar rotation matrix using heading angle $\theta_i$, and $(x_i,y_i)$ are the position coordinates of configuration $\mathbf{q}_i$. This transformation places each point $\mathbf{p}$ of motion primitive $m_{k_i}$ from its local coordinate frame into the global coordinate frame. This process efficiently reconstructs a continuous trajectory that respects kinematic constraints while avoiding obstacles. The complete pseudocode implementation is presented in Algorithm~\ref{alg:A*} of the Appendix A, with the different modules detailed in the following subsections.

\subsubsection{Geometric Feasibility Checking}\label{sec:vehicle representation planner}
For efficient graph construction, vehicle configurations that intersect with road boundaries or prohibited areas must be eliminated. Due to the non-convex nature of road boundaries, this necessitates an efficient implementation to check if vehicle positions at graph nodes collide with (overlap) the infrastructure. Although vehicle geometry could be represented using concave polygons for exact geometry or convex polygons for simplified geometry, both approaches incur significant computational overhead from complex angle calculations or intensive vertex operations. To optimize computational efficiency while maintaining acceptable accuracy, we propose a dual-circle approximation of vehicle geometry. This representation encapsulates the vehicle's extremities using two circles $C_1$ and $C_2$, as illustrated in Figure~\ref{fig:Vehicle Representation}. This representation reduces the collision detection to efficient vector operations as follows:
\begin{figure}[t]
    \centering
    \includegraphics[width=0.2\textwidth]{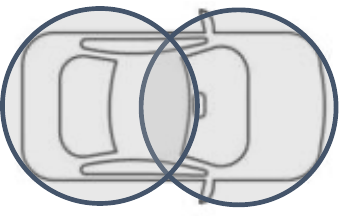}
    \caption{Representation of vehicle using two circles for collision checking.}
    \label{fig:Vehicle Representation}
\end{figure}
Let each circle $C_i$ be defined by its center $(x_{ci}, y_{ci})$ and radius $r_i$. For infrastructure elements represented as line segments between points $L_1(x_{L_1}, y_{L_1})$ and $L_2(x_{L_2}, y_{L_2})$, we formulate collision detection through vector projection and clamping operations (Figure~\ref{fig:collision_vehicle_infra}).
\begin{figure}[t]
    \centering        \includegraphics[width=0.3\textwidth]{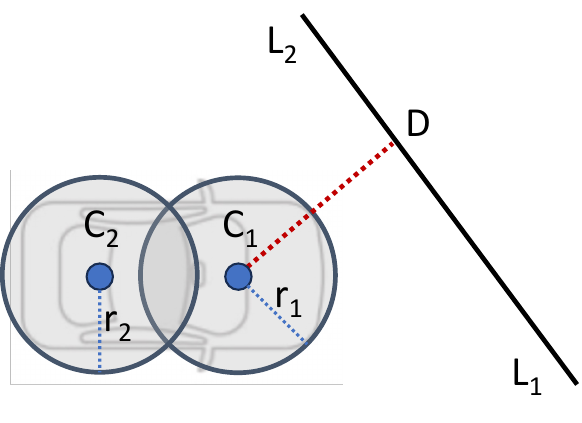}
    \caption{Collision-checking between a vehicle and the  infrastructure}
    \label{fig:collision_vehicle_infra}
\end{figure}
We first define the directional and relative position vectors:
\begin{equation}
    \vec{L_1L_2}=L_2-L_1=
        \begin{pmatrix}
            x_{L_2}-x_{L_1}\\
            y_{L_2}-y_{L_1}
        \end{pmatrix},\qquad
    \vec{L_1C_i}=
        \begin{pmatrix}
            x_{ci}-x_{L_1}\\
            y_{ci}-y_{L_1}
        \end{pmatrix}
\end{equation}
The normalized projection scalar $t$ is computed as
\begin{equation}
    t=\frac{\vec{L_1C_i}\cdot\vec{L_1L_2}}{\|\vec{L_1L_2}\|^{2}}.
\end{equation}
To ensure the projected point lies on the finite line segment, we apply the clamping $t=\max(0,\min(1,t))$, and then find the closest point $D$ on the line segment as $D=L_1+t\cdot\vec{L_1L_2}$. A collision is detected if either circle intersects with the line segment, formalized as
\begin{equation}
\text{Collision}\;\Leftrightarrow\;\exists\,i\in\{1,2\}:\;
\|\mathbf{C}_i-\mathbf{D}\|=
\sqrt{(x_{di}-x_{ci})^{2}+(y_{di}-y_{ci})^{2}}\le r_{ci},
\end{equation}
where $(x_{di},y_{di})$ are the coordinates of $D$, $(x_{ci},y_{ci})$ are the coordinates of the centre of $C_i$ and $r_{ci}$ is its radius. These efficient matrix operations significantly decrease computational overhead in the global path-planning phase.

\subsection{Motion Controller}\label{sec:motion controller}
For the motion controller, we implement a Model Predictive Control approach using a linearized kinematic vehicle model. The proposed controller architecture decouples trajectory tracking from collision prediction and avoidance to deliver computational efficiency, enhance modularity, and enable independent verification of planning and prediction algorithms. This section presents the kinematic model representation, its linearization for computational efficiency, the MPC formulation, and the collision avoidance strategy.

\subsubsection{Vehicle Kinematics}\label{sec:vehicle dynamics}
\begin{figure}[t] 
    \centering 
    \includegraphics[width=0.3\textwidth]{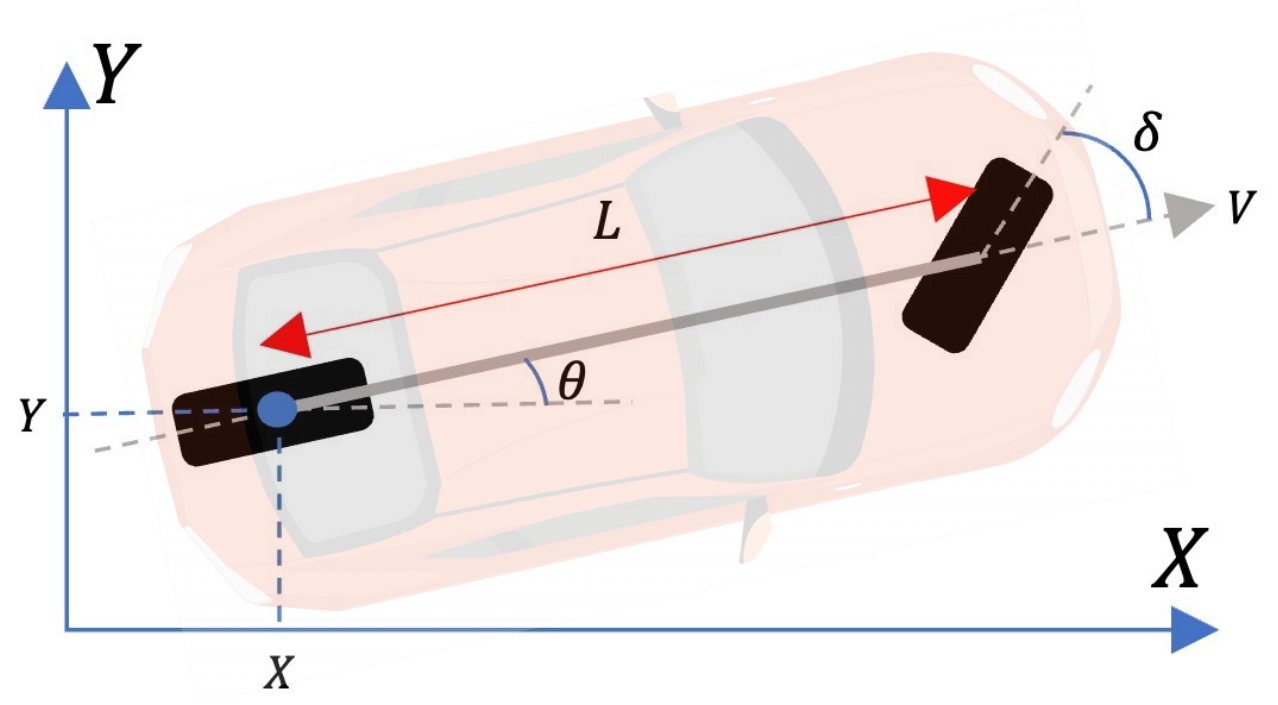} \caption{Kinematic bicycle model representation of a vehicle.} \label{fig:bicycle_model} 
\end{figure}

Figure \ref{fig:bicycle_model} depicts the kinematic vehicle model used in this study, also known as the kinematic bicycle model (\textcite{polack2017kinematic}). This representation captures essential two-dimensional vehicle motion through the state vector $\mathbf{x} = [x, y, v, \theta]^{\top}$, where $x$ and $y$ denote position coordinates (meters), $v$ velocity (m/s), and $\theta$ heading angle (radians). Accordingly, the equations of motion are achieved as follows:
\begin{gather}\label{eq:nonlinear_bicycle_model}
 \dot{\mathbf{x}}
 =
 f(\mathbf{x},\mathbf{u})
 =
 \begin{bmatrix}
         \dot{x} \\ 
        \dot{y} \\ 
        \dot{v}\\ 
        \dot{\theta}
 \end{bmatrix}
 =
  \begin{bmatrix}
        vcos(\theta) \\ 
        vsin(\theta) \\ 
        a\\ 
        \frac{v}{L}tan(\delta)
   \end{bmatrix}
\end{gather}
where control inputs $\mathbf{u} = [a, \delta]^{\top}$ represent acceleration and steering angle respectively, and $L$ denotes the vehicle wheelbase.
While this nonlinear formulation offers high fidelity in vehicle dynamics representation, its computational complexity is not optimal for microsimulation applications. For scenarios where vehicles operate below kinematic limits, as typically observed in urban intersections (speeds under 30 to 50 km/h), a linearized approximation provides sufficient accuracy with significantly reduced computational overhead (\cite{polack2017kinematic}). We perform this linearization using first-order Taylor series expansion around a nominal operating point $(\bar{\mathbf{x}}, \bar{\mathbf{u}})$ as given by Equation \ref{eq:taylor_series}:
\begin{equation}\label{eq:taylor_series}
        \dot{f}(\mathbf{x},\mathbf{u}) \approx f(\bar{\mathbf{x}}, \bar{\mathbf{u}}) + \frac{\partial f(\mathbf{x},\mathbf{u})}{\partial \mathbf{x}} |_{\bar{\mathbf{x}},\bar{\mathbf{u}}} (\mathbf{x} - \bar{\mathbf{x}})
    + \frac{\partial f(\mathbf{x},\mathbf{u})}{\partial \mathbf{u}} |_{\bar{\mathbf{x}},\bar{\mathbf{u}}} (\mathbf{u}-\bar{\mathbf{u}})     
\end{equation}
where $\bar{\mathbf{x}}$ and $\bar{\mathbf{u}}$ represent the state and input elements of the operating point. Letting \( A = \frac{\partial f}{\partial \mathbf{x}}\Big|_{\bar{\mathbf{x}}, \bar{\mathbf{u}}} \) and \( B = \frac{\partial f}{\partial \mathbf{u}}\Big|_{\bar{\mathbf{x}}, \bar{\mathbf{u}}} \) and re‑arranging, the familiar linear state‑space form emerges:
\begin{equation}\label{eq:linear_state_space}
\dot{\mathbf{x}} = A\,\mathbf{x} + B\,\mathbf{u} + \mathbf{d},
\end{equation}
with
\[
A =
\begin{bmatrix}
0 & 0 & \cos\bar\theta & -\bar v\sin\bar\theta\\
0 & 0 & \sin\bar\theta & \;\;\bar v\cos\bar\theta\\
0 & 0 & 0 & 0\\
0 & 0 & \tfrac{\tan\bar\delta}{L} & 0
\end{bmatrix},
\qquad
B =
\begin{bmatrix}
0 & 0\\
0 & 0\\
1 & 0\\
0 & \tfrac{\bar v}{L\cos^2\bar\delta}
\end{bmatrix},
\qquad
\mathbf{d} =
\begin{bmatrix}
\bar v\sin\bar\theta\,\bar\theta\\
-\bar v\cos\bar\theta\,\bar\theta\\
0\\
-\tfrac{\bar v\bar\delta}{L\cos^2\bar\delta}
\end{bmatrix}.
\]
where, the operating point of $\theta$, $v$, and $\delta$ are $\bar{\theta}$, $\bar{v}$, and $\bar{\delta}$, respectively. Derivations of $A$, $B$, and $d$ are provided in Appendix B. For digital implementation we sample the system with a zero‑order hold
at a fixed interval of \(T_s = 0.1\,\text{s}\).
Using the forward‑Euler (first‑order) ZOH mapping,
\begin{equation}
    \mathbf{x}(k+1) = A_{d}\mathbf{x}(k) + B_{d}\mathbf{u}(k) + \mathbf{d}_{d}
    \label{eq:discretized_state_space}
\end{equation}
with
\[
A_d = I + A\,T_s,\qquad
B_d = B\,T_s,\qquad
\mathbf{d}_d = \mathbf{d}\,T_s.
\]
where \(k\) is the discrete‑time index (\(t = kT_s\)).

\subsubsection{MPC Formulation}
\begin{figure}[h]
    \centering
    \includegraphics[width=0.6\textwidth]{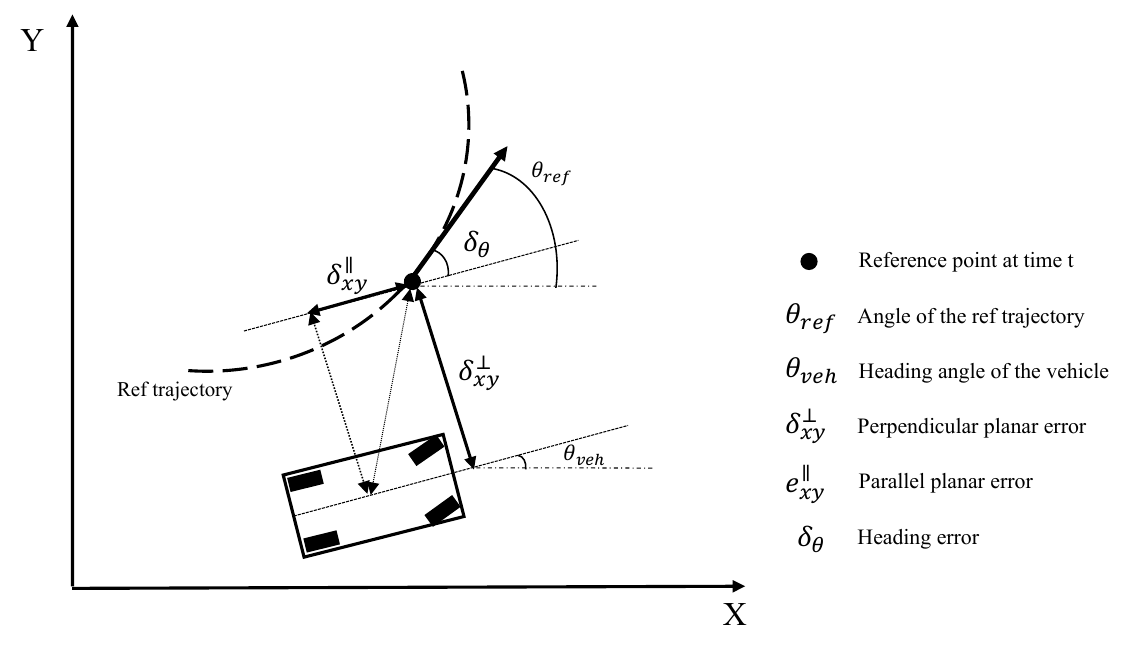}
    \caption{Demonstration of perpendicular (\( \delta_{k}^\perp \)) and parallel (\( \delta_{k}^{\parallel} \)) deviations from the reference trajectory defined in the MPC formulation.}
    \label{fig:mpc_deviation}
\end{figure}

The proposed controller is designed to address three primary objectives: trajectory tracking through minimization of positional, orientational, and velocity deviations; maneuver smoothness via penalization of control input magnitudes and their temporal derivatives; and terminal state convergence. Mathematically, the cost function of the controller is formulated as $\mathbf{J}$ over the control horizon $N_c$:
\begin{equation}
    J = \sum_{k=1}^{N_c-1} \left( \|e_{xy,k}\|^2_{Q_{xy}(k)} + \|\delta_{\theta v,k}\|^2_{Q_{\theta v}} \right) 
        + \sum_{k=0}^{N_c-1} \|\mathbf{u}(k)\|^2_{R} 
        + \sum_{k=0}^{N_c-2} \|\mathbf{u}(k+1) - \mathbf{u}(k)\|^2_{R_d} 
        + \|\mathbf{x}_{N_c} - \mathbf{x}_{\text{ref},N_c}\|^2_{Q_f}
\end{equation}
where: 
\begin{itemize}
    \item $e_{xy,k}$ is the planar position error $(\mathbf{x}_{xy}(k)-\mathbf{x}^{\text{ref}}_{xy}(k))$ at prediction step~$k$, where $\mathbf{x}_{xy}(k)=[x(k),y(k)]^{\top}$ is the actual position, and $\mathbf{x}^{\text{ref}}_{xy}(k)$ is its reference.
    \item $\delta_{\theta v,k}$ is the heading and speed error $(\mathbf{x}_{\theta v}(k)-\mathbf{x}^{\text{ref}}_{\theta v}(k))$, where $x_{\theta v}(k)=[\theta(k),v(k)]^{\top}$ is the actual heading and speed, and $\mathbf{x}^{\text{ref}}_{\theta v}(k)$ is the reference.  Minimising this error aligns the vehicle’s orientation with the path tangent and maintains the desired velocity profile to ensure the vehicle follows the planned while maintaining the desired speed and heading. 
    \item $\mathbf{u}(k)=[\delta(k),\,a(k)]^{\top}$ is the control input, where $\delta$ is the steering angle and $a$ is the longitudinal acceleration. Penalizing this term prevents overly aggressive control inputs, and its differential form $\mathbf{u}(k+1)-\mathbf{u}(k)$ suppress abrupt changes in the control input for passenger comfort.
    \item $N_c$ is the control horizon; $k=0$ denotes the current sample and $k=N_c$ the terminal stage.
    \item $\mathbf{x}_{N_c}$ and $\mathbf{x}_{\mathrm{ref},N_c}$ are the actual and reference state vectors at stage~$N_c$. Their weighted quadratic difference, scaled by $Q_{f}$, drives the optimiser towards the desired final pose and speed, ensuring closed-loop convergence.
    \item $Q_{xy}$, $Q_{\theta v}$, $R$, $R_{d}$ and $Q_{f}$ are constant positive-semidefinite weighting matrices. By adjusting these matrices the designer balances path-tracking accuracy, control effort, input smoothness and terminal precision, thereby shaping the overall driving behaviour.
\end{itemize}
The position-error weight $Q_{xy}$ blends perpendicular (cross-track) and parallel (along-track) components as shown in Figure \ref{fig:mpc_deviation} and elaborated below:
\begin{equation}
        Q_{xy}(k)=w_{\perp}\,Q_{\perp}(k)+w_{\parallel}\,Q_{\parallel}(k),
    \qquad
    w_{\perp}\ge w_{\parallel}>0,
\end{equation}
where
\begin{equation}
        Q_{\perp}(k)=
    \begin{bmatrix}
      \sin^{2}\theta_{\text{ref},k} & -\sin\theta_{\text{ref},k}\cos\theta_{\text{ref},k}\\
      -\sin\theta_{\text{ref},k}\cos\theta_{\text{ref},k} & \cos^{2}\theta_{\text{ref},k}
    \end{bmatrix},
    \quad 
    Q_{\parallel}(k)=
    \begin{bmatrix}
      \cos^{2}\theta_{\text{ref},k} &  \cos\theta_{\text{ref},k}\sin\theta_{\text{ref},k}\\
      \cos\theta_{\text{ref},k}\sin\theta_{\text{ref},k} &  \sin^{2}\theta_{\text{ref},k}
    \end{bmatrix}.
\end{equation}
Both matrices are rank-one projectors:  
$Q_{\parallel}=n_{\parallel}n_{\parallel}^{\!\top}$ with
$n_{\parallel}=[\cos\theta_{\text{ref},k},\,\sin\theta_{\text{ref},k}]^{\!\top}$ (tangential direction),  
and $Q_{\perp}=n_{\perp}n_{\perp}^{\!\top}$ with
$n_{\perp}=[-\sin\theta_{\text{ref},k},\,\cos\theta_{\text{ref},k}]^{\!\top}$ (normal direction).  
Consequently,
$e_{xy,k}^{\!\top}Q_{xy}(k)e_{xy,k}=w_{\perp}e_{\perp,k}^{2}+w_{\parallel}e_{\parallel,k}^{2}$,
so lateral accuracy can be emphasised by selecting $w_{\perp}>w_{\parallel}$. Because all weighting matrices are positive-semidefinite, the cost is convex, and the resulting MPC problem can be solved
efficiently with standard quadratic programming methods.

The final optimization problem is achieved by adding the constraints as follows:
\begin{equation}\label{eq:optimization_problem}
\begin{aligned}
\min \quad & J\\
\text{s.t.}\quad & \mathbf{x}(0)=\mathbf{x}_{\text{init}},\\
                 & \mathbf{x}(k+1)=A_{d}\mathbf{x}(k)+B_{d}\mathbf{u}(k)+\mathbf{d}_{d},\qquad k=0,\dots,N_c-1,\\
                 & |\dot{\delta}(k)|\le\dot{\delta}_{\max},\qquad |\delta(k)|\le\delta_{\max},\\
                 & v_{\min}\le v(k)\le v_{\max},\qquad a_{\min}\le a(k)\le a_{\max}.
\end{aligned}
\end{equation}
Here $\mathbf{x}_{\text{init}}$ is the initial state; $\dot{\delta}_{\max}$ and $\delta_{\max}$ bound the steering rate and steering angle; $v_{\min}$ and $v_{\max}$ bound the speed; and $a_{\min}$ and $a_{\max}$ bound the acceleration. These constraints ensure that the predicted motion respects the kinematic and non-holonomic limits of the linearised bicycle model used in this study. Because the constraints are linear and all weighting matrices are positive-semidefinite, the resulting finite-horizon problem is a convex quadratic programme that can be solved efficiently at every sampling instant. To solve the formulated convex quadratic optimization problem at each sampling instant, we used the Python-based optimization library CVXPY with the ECOS solver, which provides efficient performance for real-time applications. Collision avoidance constraints are excluded from the above MPC formulation because moving obstacle handling is delegated to a dedicated module, described in the next section. 

\subsubsection{Collision Avoidance}
The decentralized collision avoidance module operates in real-time and in parallel with the MPC controller. This module relies on finite-horizon state predictions for both the ego vehicle and proximate vehicles to ensure safe navigation. Two key parameters are defined: (1) a detection range \( R_{\text{detect}} \), specifying the perceptual region within which the ego vehicle monitors state vectors \(\mathbf{x}_j = [x_j, y_j, v_j, \theta_j]^\top\) of nearby vehicles, and (2) a prediction horizon \( T_{\text{pred}} \), defining the future time window for state predictions. This horizon is discretized into \( N_{\text{pred}} \) timesteps of duration \(\Delta t\), starting from the current timestep \(t_k\), such that \(t_{k+\ell} = t_k + \ell\Delta t\) for \(\ell = 0, 1, \dots, N_{\text{pred}}\). The prediction horizon exceeds the MPC control horizon \(N_c\), providing a temporal buffer to anticipate and mitigate collision risks beyond the optimization window \(t_k + N_c\). For this study, a straightforward collision avoidance strategy is implemented. Surrounding vehicles are modeled with constant velocity and steering angle assumptions, while the ego vehicle predicts its states assuming acceleration toward a desired speed \( v_{\text{desired}} \) along the reference trajectory. Recursive updates at each timestep \(\Delta t\) ensure effective collision mitigation despite simplifications.

\begin{figure}[b]
    \centering
    \includegraphics[width=0.3\textwidth]{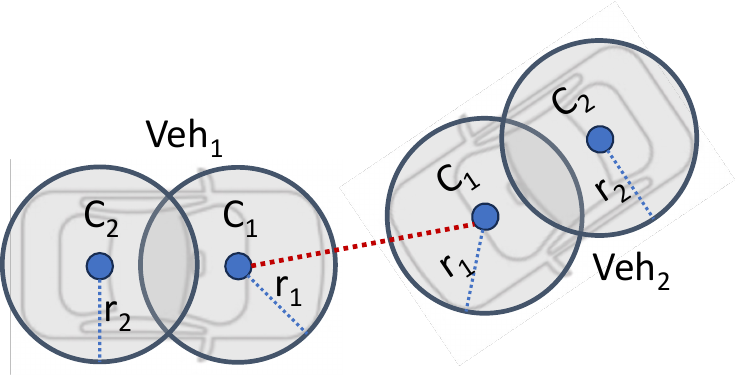}
    \caption{Dual-circle representation used for efficient collision detection between vehicles.}
    \label{fig:collision_vehicle_vehicle}
\end{figure}

\begin{figure}[t]
\centering
\includegraphics[width=0.5\textwidth]{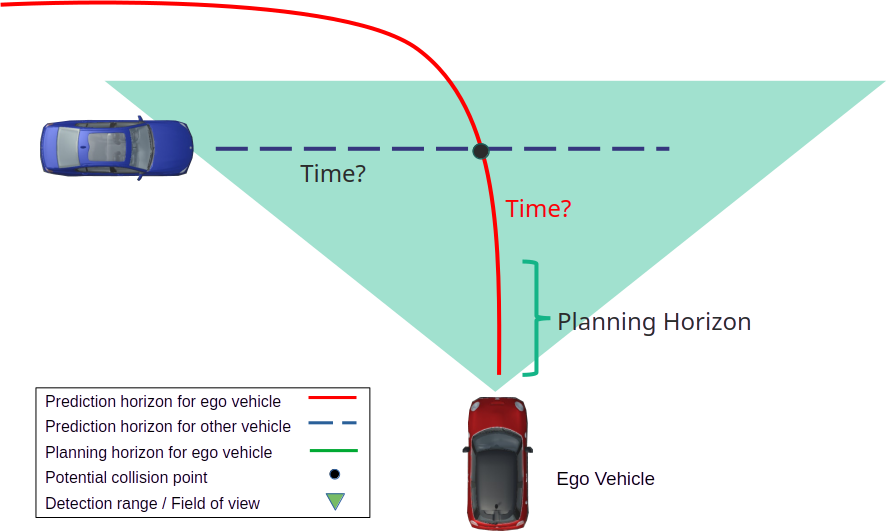}
\caption{Different horizons in the collision detection process.}
\label{fig:collision_detection}
\end{figure}

Future states are computed at each timestep \(t_k\) using a discrete-time approximation of the linearized bicycle model:
\begin{equation}
  \mathbf{x}(t_k + \ell\Delta t) = f(\mathbf{x}(t_k + (\ell - 1)\Delta t), \mathbf{u}(t_k + (\ell - 1)\Delta t)), \quad \ell = 1, \dots, N_{\text{pred}}  
\end{equation}
Predicted positions form reachable sets, represented using dual-circle vehicle approximations (see Figure~\ref{fig:collision_vehicle_vehicle}):
\begin{equation}
    Veh_i(t) =
\begin{bmatrix}
x_{if}(t) & y_{if}(t) & r_i \\
x_{ir}(t) & y_{ir}(t) & r_i
\end{bmatrix},
\end{equation}
where front and rear circle centers of vehicle \(i\) at time \(t\) are \((x_{if}(t), y_{if}(t))\) and \((x_{ir}(t), y_{ir}(t))\), each with radius \(r_i\).
Collision detection between the ego vehicle \(ego\) and any other vehicle \(j\) is performed efficiently via vectorized operations. Specifically, the Euclidean distances between all pairs of circles from the ego and other vehicles are computed as:
\begin{equation}
    D_{mn}(t) = \| Veh_{ego,m}(t) - Veh_{j,n}(t) \|_2,
\end{equation}
for \(m,n \in \{f,r\}\), where \(Veh_{ego,m}(t)\) denotes the position coordinates of the \(m\)-th circle (front or rear) of the ego vehicle at time \(t\). A collision at time \(t\) occurs if:
\begin{equation}
    D_{mn}(t) \leq r_{ego,m} + r_{j,n}, \quad \text{for any combination of } m,n \in \{f,r\}.
\end{equation}
Optionally, an analytic test acts as an early-exit filter:
\begin{equation}
    \| \mathbf{r}_{ego,j}(0) + \tau \dot{\mathbf{r}}_{ego,j} \|_2^2 = (r_{ego} + r_j + \epsilon)^2,
\end{equation}
solved for \(\tau\) to assess collision imminence, where \(\mathbf{r}_{ego-j}(0)\) and \(\dot{\mathbf{r}}_{ego-j}\) are the initial relative position and velocity vectors, and \(\epsilon\) is a safety margin. Upon detecting a potential collision, the ego vehicle's trajectory is modified to stop before the conflict point by solving:
\begin{equation}
    \min_{a(t)} J = \frac{1}{2}\int_{t_k}^{t_k + \tau^*} a(t)^2\,dt, \quad \text{subject to } v(t_k + \tau^*) = 0,
\end{equation}
yielding a constant deceleration \(a^* = -v_0/\tau^*\), where \(v_0\) is the initial speed and \(\tau^*\) is the time to stop, determined from the predicted collision time. The reference speed is updated to:
\begin{equation}
    v_{\text{ref}}(t) = v_0 + a^*(t - t_k),
\end{equation}
and the geometric path is truncated at the stopping distance \(s_{\text{stop}} = \frac{1}{2}v_0\tau^*\), with zero speed assigned beyond this point, forming a collision-free trajectory \(\hat{\tau}_{ego}\) (see Figure~\ref{fig:collision_detection}). This approach ensures robust, real-time collision avoidance without complicating the MPC formulation, balancing theoretical rigor and practical efficiency.

\section{Experiments and Results}\label{sec:experiments}
This section evaluates the proposed framework for modeling two-dimensional vehicular maneuvers at unsignalized intersections through numerical and simulation experiments. The experiments are organized into two primary categories: assessments of the global planner’s ability to efficiently generate kinematically feasible trajectories, and the evaluations of the motion controller and collision avoidance module.

\subsection{Evaluation of Global Planner}
Experiments in this section assess four key aspects: trajectory plausibility, adaptability to varied road configurations, computational efficiency, and sensitivity to different parameters and the customizability of the global planner. The set of parameters used for conducting the first three experiments is summarized in Table \ref{tab:global_planner_parameters}.

\begin{table}[h]
\centering
\caption{Parameters used in the global planner for plausibility and generalizability analyses}
\begin{tabular}{ll}
\hline
\textbf{Parameter Name} & \textbf{Value} \\
\hline
Vehicle dimension & 4 meters \\
Number of MPs & 9 \\
Safety margin & 0.5 meters \\
\(w_{h,d}\) & 1.0 \\
\(w_{h,\theta}\) & 2.7 \\
\(w_{h,\phi}\) & 15.0 \\
\(w_{c,d}\) & 1.0 \\
\(w_{c,\phi}\) & 5.0 \\
\hline
\end{tabular}
\label{tab:global_planner_parameters}
\end{table}

\subsubsection{Plausibility Analysis}
This analysis evaluates the global planner's adherence to non-holonomic kinematic constraints during reference trajectory generation. Plausibility verification is fundamental from multiple theoretical and practical perspectives. Theoretically, it ensures that generated trajectories belong to the vehicle's reachable set and satisfy differential constraints imposed by the kinematic model, preventing paths that violate curvature continuity or exceed steering rate limitations. Computationally, constraint adherence enables the graph search algorithm to operate within the feasible configuration space, systematically pruning kinematically inadmissible regions and reducing search complexity. From a control-theoretic standpoint, kinematically consistent reference trajectories ensure superior tracking performance by operating within natural system dynamics, minimizing control effort while reducing actuator saturation likelihood and maintaining trajectory fidelity.

\begin{figure}[t]
    \centering
    \begin{subfigure}[b]{0.43\textwidth}
        \centering
        \includegraphics[width=\textwidth]{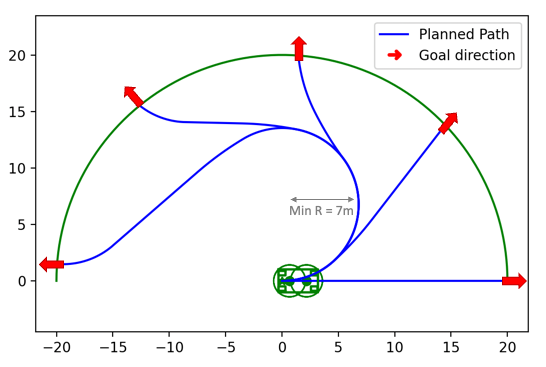}
        \caption{Turning maneuvers}
        \label{fig:plausibility_a}
    \end{subfigure}
    \hfill % Optional: add some space between the two figures
    \begin{subfigure}[b]{0.4\textwidth}
        \centering
        \includegraphics[width=\textwidth]{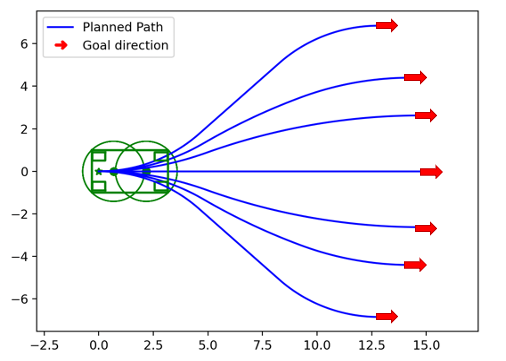}
        \caption{Lane changing and lateral maneuvers}
        \label{fig:plausibility_b}
    \end{subfigure}
    \caption{Plausibility analysis of the global planner in two common scenarios}
    \label{fig:plausibility}
\end{figure}

To validate kinematic feasibility across diverse maneuvering scenarios, two complementary experimental configurations were designed (Figure~\ref{fig:plausibility}). The first configuration, depicted in Figure~\ref{fig:plausibility_a}, includes a radial test. Here the vehicle is positioned at the center of a circular domain while the terminal states are distributed along the perimeter at varied orientations (Figure 11a). This arrangement evaluates the planner's capability to generate diverse maneuvers that are expected at intersections and roundabouts, including U-turns. The generated trajectories in Figure~\ref{fig:plausibility_a} (blue lines) extend from the initial vehicle state (represented by the dual-circle geometric approximation) to goal locations (indicated by red directional vectors). The results demonstrate that the planner successfully generates smooth and kinematically feasible trajectories for all tested orientations. Each path respects the vehicle's turning radius constraints and exhibiting continuous curvature profiles. 

The second configuration adopts a forward-projection approach, evaluating lateral maneuvering, which can be observed when vehicles change their lane or their lateral position. This experiment is designed by setting the goal states positioned frontward at various lateral offsets (Figure~\ref{fig:plausibility_b}). The results demonstrate that the planner effectively generates smooth lateral transitions while maintaining proper curvature constraints and final heading directions. These findings confirm the planner's ability to produce kinematically feasible paths for lateral displacement scenarios, which are commonly encountered at intersection approaches and departure zones. Collectively, these experimental validations confirm that the adopted motion primitives and the designed graph search algorithm effectively discretize and search through the continuous vehicle reachable space while preserving kinematic feasibility. This confirms the planner's theoretical foundation and practical applicability for intersection modeling.

\subsubsection{Generalizability Across Different Intersection Designs}
This experiment evaluates the planner's adaptability across diverse junction configurations to verify its generalizability beyond specific intersection typologies. Four scenarios are constructed, comprising a T-intersection, dual four-leg intersections (single-lane and multi-lane configurations), and a roundabout. Various maneuvers were examined within each topological class.

\begin{figure}[t]
    \centering
    \includegraphics[width=1\textwidth]{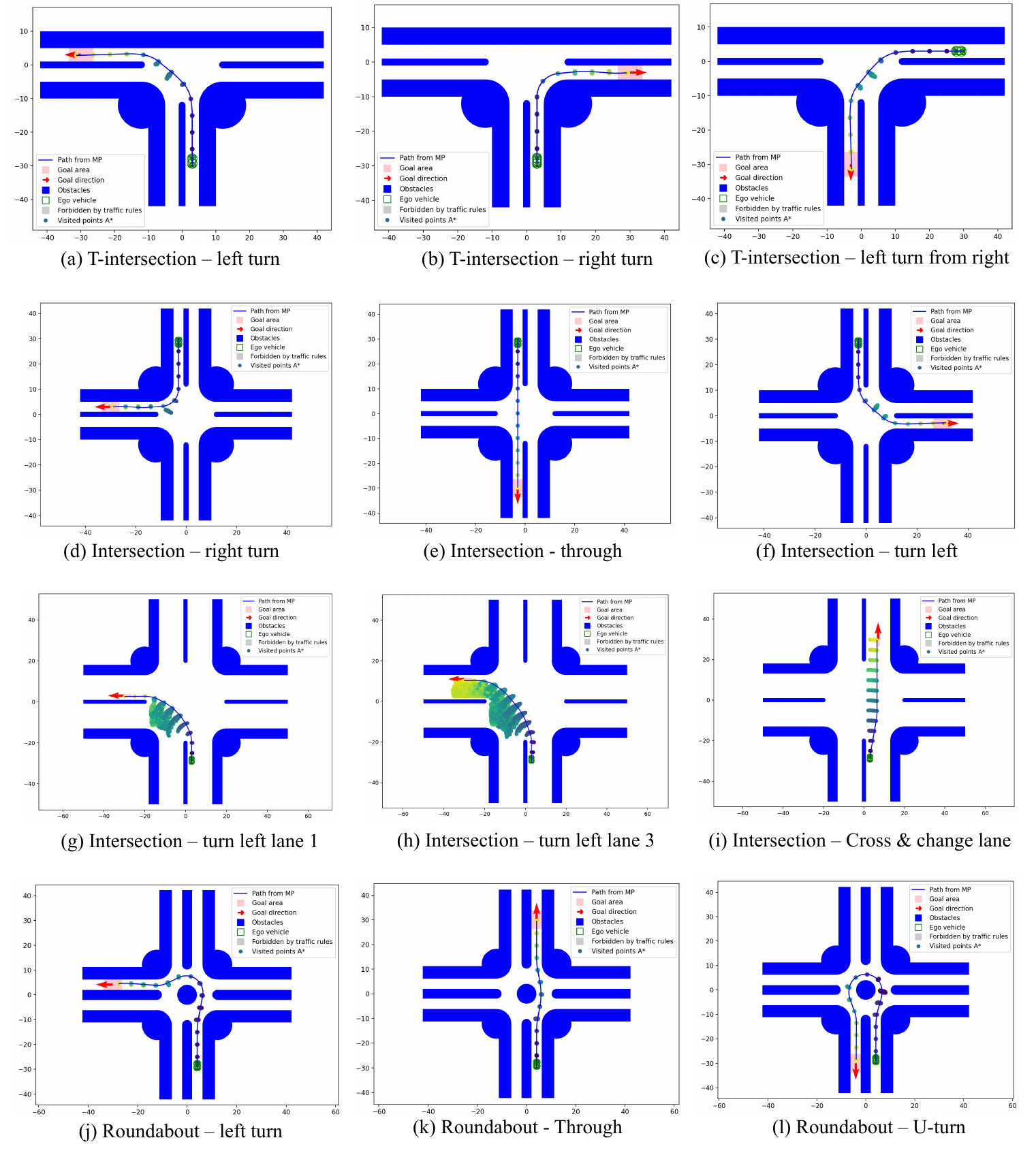}
    \caption{Generalizability of the proposed planner in different scenarios}
    \label{fig:generalizability}
\end{figure}

Figure \ref{fig:generalizability} presents the results for this experiment. In each subfigure, blue geometric shapes (rectangles and circles) represent environmental obstacles that represent the intersection geometry. The solid blue curves indicate the generated reference trajectory from the location of the ego vehicle (double-green circles) to the goal area indicated by red rectangular and arrow. Scattered colored dots represent nodes explored by the search algorithm, with their colors corresponding to heuristic function values. 
Figures \ref{fig:generalizability}.a-c demonstrate various maneuvers at a T-intersection. The generated reference trajectories and explored search nodes, represented by colored dots, indicate that the global planner successfully identifies smooth and traversable paths from initial positions to goal regions while maintaining a constrained search space. This demonstrates the effectiveness of the designed heuristic function. Figures \ref{fig:generalizability}.d-f present trajectory synthesis results for a four-leg intersection. These results show that the planner generates both smooth curved trajectories for turning maneuvers and straight paths when required, such as for through movements. Similar to the previous scenario, the consistently bounded search space indicates the computational efficiency of the heuristic-guided search methodology.

The results of path planning in more challenging scenarios are presented in Figures \ref{fig:generalizability}.g-l. These scenarios include a multi-lane intersection requiring concurrent lane changes while crossing the intersection, as well as a roundabout where the vehicle must avoid the central island, respect traffic regulations, and execute complex maneuvers such as U-turns. The multi-lane intersection results in Figures \ref{fig:generalizability}.g-h demonstrate that the planner successfully synthesizes smooth and lane-specific trajectories from the origin lane to the destination lane. Figure~\ref{fig:generalizability}.i, specifically, reveals the planner's capability to execute smooth lane transitions during intersection crossing, combining longitudinal and lateral motion in a single maneuver. The distribution of the colored nodes reveals that despite the larger feasible and drivable space, the search algorithm maintains computational efficiency through strategic exploration towards the goal area. The results for the roundabout scenarios also indicate that the generated trajectories fully comply with junction design and traffic regulations although the central island creates a non-convex feasible region for the search algorithm. Notably, Figure~\ref{fig:generalizability}.k shows the planner synthesizing a curved trajectory to circumnavigate the central island in order to maintain appropriate safety margins.

Collectively, these findings demonstrate the planner's consistent generation of kinematically feasible and regulation-compliant trajectories across varied topological configurations. The generated paths exhibit realistic vehicle dynamics without sharp angles and close proximity to the obstacles. This experimental validation confirms the global planner's capability to generate realistic and consistent trajectories across diverse intersection typologies without the need for parametric recalibration. However, it's worth noting that parametric optimization depending on the problem-specific design may improve the quality of computational burden of the planner.

\subsubsection{Sensitivity Analysis}\label{sec:planner_sensitivity}
The global planner is formulated with parametric sensitivity to user-defined optimization criteria through a configurable parameter vector. This property distinguishes the proposed framework from previous studies that employ deterministic, pre-planned trajectories (\cite{zhan2017spatially, liu2017speed}).

\begin{figure}[ht!]
\centering
\begin{subfigure}[b]{0.35\textwidth}
    \centering
    \includegraphics[width=\linewidth]{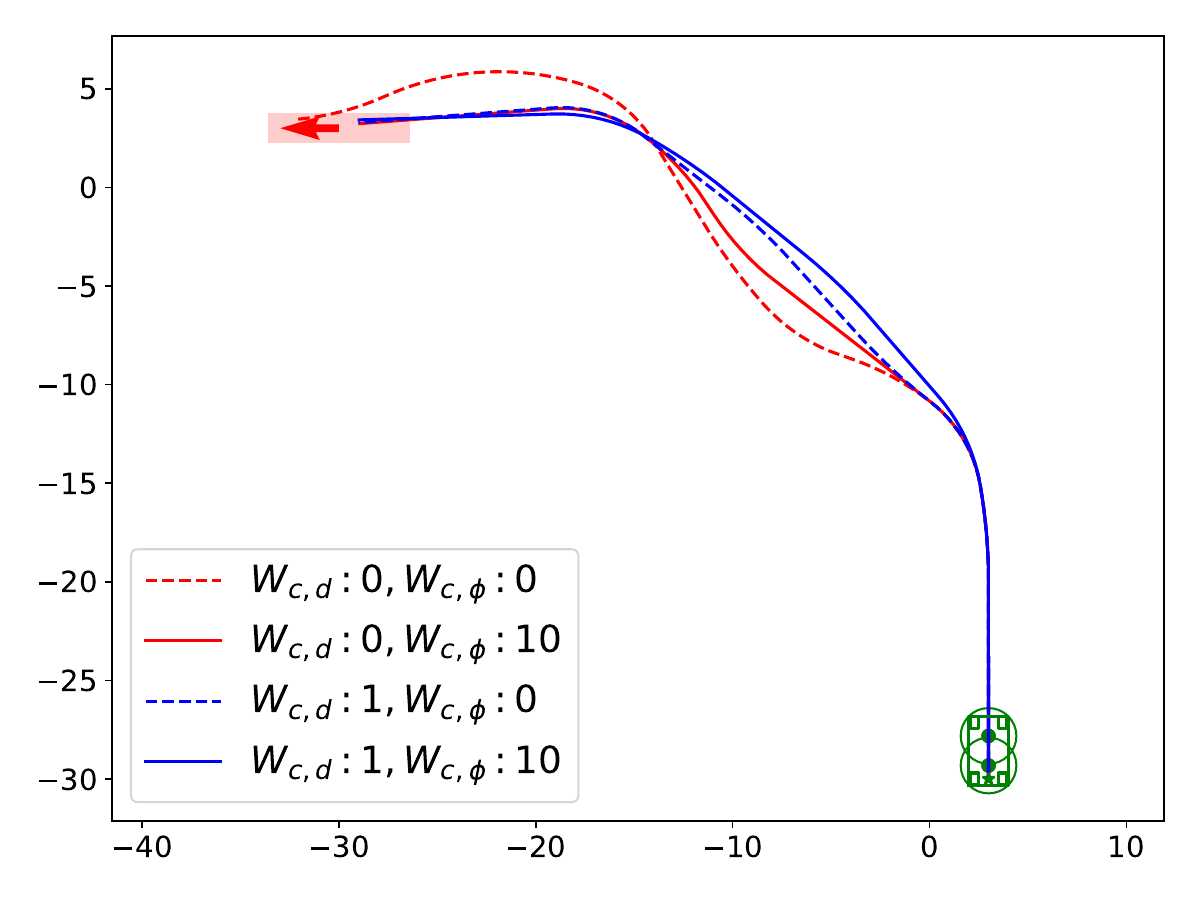}
    \caption{Real cost parameters}
    \label{fig:sensitivity_true_cost}
\end{subfigure}
\hspace{0.05\textwidth} % Adjust the horizontal space between the plots as needed
\begin{subfigure}[b]{0.35\textwidth}
    \centering
    \includegraphics[width=\linewidth]{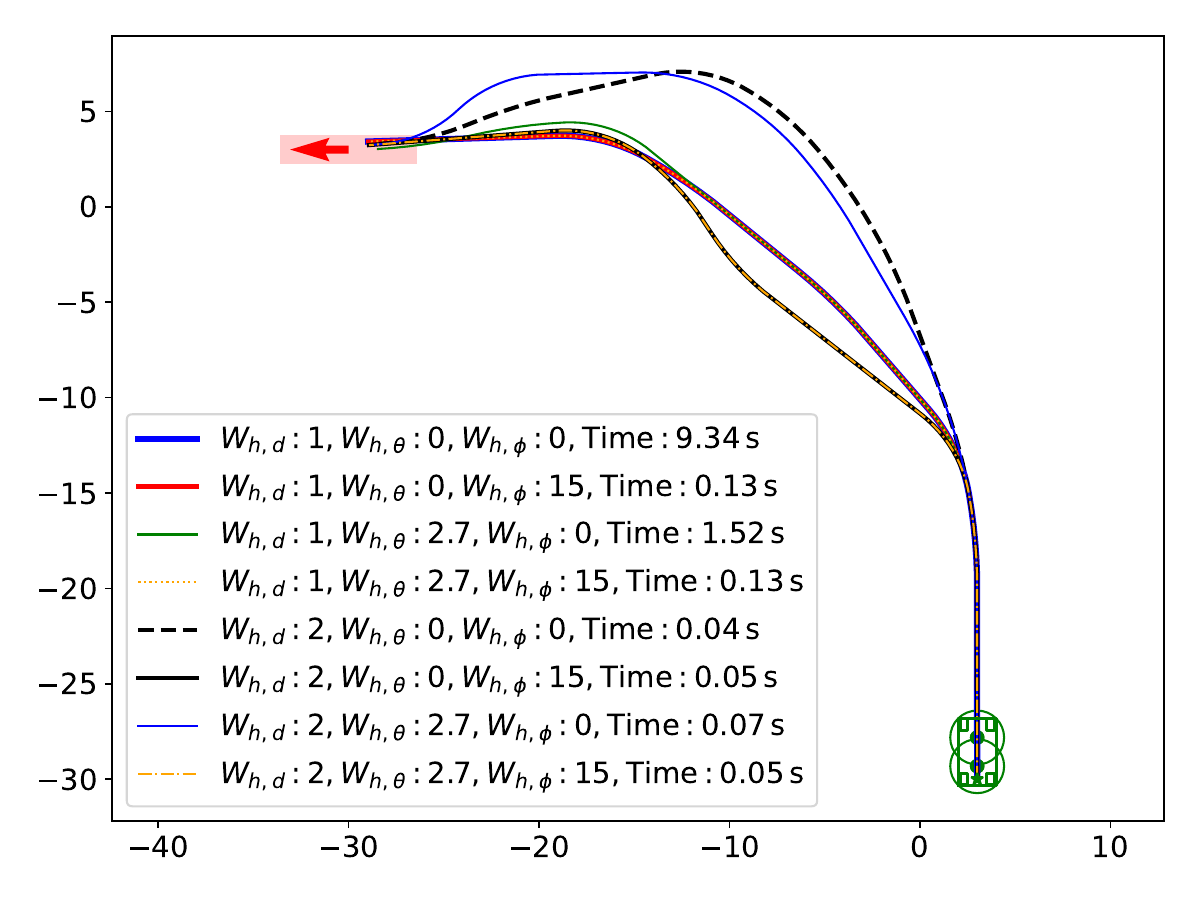}
    \caption{Heuristic cost parameters}
    \label{fig:sensitivity_heuristic}
\end{subfigure}

\caption{Sensitivity Analysis on Global Planner's Parameters}
\label{fig:sensitivity_planner}
\end{figure}

To quantify the influence of different parameters of the global planner on solution characteristics and algorithmic performance, we implemented a sensitivity analysis examining both cost function parameters and heuristic function weights on the synthesized trajectory (Figure \ref{fig:sensitivity_planner}). For experimental consistency  and analytical tractability, all experiments were conducted under fixed initial configuration space coordinates and terminal region constraints to isolate parameter effects. These boundary conditions are represented by the dual-circle vehicle representation and the transparent red rectangular goal region, respectively, in Figure~\ref{fig:sensitivity_planner}.

Figure~\ref{fig:sensitivity_true_cost} demonstrates the influence of cost function parameters $w_{c,d}$ and $w_{c,\phi}$ on optimal trajectory topology. The absence of the steering effort penalty $w_{c,\phi}$ results in oscillatory trajectory characteristics, while the elimination of the path length weight $w_{c,d}$ produces suboptimal solutions with respect to the defined cost functional. Figure~\ref{fig:sensitivity_heuristic} presents the analysis of heuristic parameter effects, establishing that appropriate weight calibration substantially reduces computational complexity while preserving solution optimality. Specifically, the variation of the orientation alignment weight $w_{h,\theta}$ and control effort estimation weight $w_{h,\phi}$ achieves computational efficiency improvements of two orders of magnitude (reducing execution time from 9.34s to 0.05s) while generating geometrically equivalent optimal trajectories. This validates the theoretical foundation that informed heuristic design enables efficient exploration of high-dimensional configuration spaces while maintaining optimality guarantees (provided the heuristic function remains admissible).

\subsubsection{Run-time and Computational Efficiency}
This section evaluates the computational efficiency of the global planner. We present a comparative analysis of our optimized A* algorithm against standard implementations, quantifying the performance gains achieved through our custom multi-criteria heuristic function. Our algorithm is benchmarked against Dijkstra's algorithm and several A* variants using single-criterion heuristics. Performance metrics are summarized in Table \ref{tab:search_algorithm_comparison}. The proposed algorithm demonstrates significant computational efficiency, exploring only 22 nodes with a runtime of 0.02 seconds, approximately 300 times faster than the naive A* implementation  that utilizes only Euclidean distance to the goal area as heuristic.
\begin{table}[ht]
\centering
\caption{Performance comparison of different pathfinding algorithms}
\label{tab:search_algorithm_comparison}
\begin{tabular}{@{}lccc@{}}
\toprule
Algorithm & No. Nodes Visited & Runtime  \\ \midrule
Dijkstra & 36420 & 19.823  \\
Naive A* (Direct distance) & 12430 & 6.245  \\
A* (Direction-guided) & 2101 & 1.78 \\
A* (Steering-guided) & 1540 & 1.24  \\
Proposed Optimized A* & 22 & 0.02  \\
\bottomrule
\end{tabular}
\end{table}

\begin{figure}[t]
    \centering
    \includegraphics[width=1\textwidth]{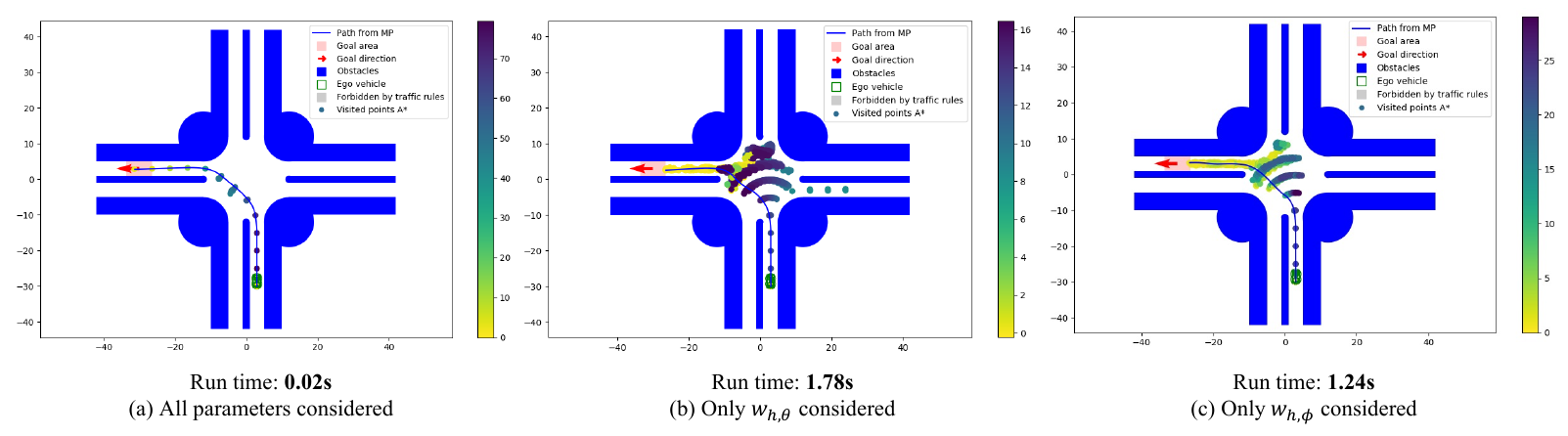}
    \caption{Run-time analysis of the proposed planner with different heuristic functions}
    \label{fig:planner_efficiency}
\end{figure}

This experiment also proves that the combined effect of considering multiple criteria in the A* algorithm is effective, whereas each criterion alone is not sufficient to achieve an efficient search algorithm. Beyond the quantitative performance metrics, Figure~\ref{fig:planner_efficiency} presents a comparative visualization of search behavior to provide deeper insights into the underlying mechanisms driving the algorithmic efficiency improvements. The colored dots represent nodes visited during path finding process to visualize the search space exploration patterns of different heuristic functions. Our proposed algorithm (Fig. \ref{fig:planner_efficiency}a) exhibits highly directed search behavior with minimum number of nodes visited, while single-criterion implementations (Fig. \ref{fig:planner_efficiency}b,c) explore substantially larger portions of the configuration space. This visual evidence, combined with the quantitative results in Table \ref{tab:search_algorithm_comparison}, confirms that the integration of multiple criteria in the heuristic function has been effective in achieving computational efficiency in kinodynamic planning tasks.

The results presented in this section confirm that the global planner achieves its design objectives of generating optimal trajectories with real-time computational performance across varied intersection topologies. While these findings validate the upper level of the bi-level architecture, successful framework implementation requires accurate trajectory execution by the motion controller in dynamic environments. The following section evaluates the motion controller's performance within the proposed framework.

\subsection{Evaluation of Motion Controller}
This section evaluates the motion controller's trajectory tracking precision and collision avoidance capabilities through a systematic parameter calibration and multi-scenario validation. The evaluation addresses three critical aspects: trajectory tracking fidelity under kinematic constraints, collision avoidance robustness in multi-vehicle scenarios, and computational efficiency for real-time implementation.

\subsubsection{Parametric Analysis and Configuration}
The motion controller's performance is directly influenced by the calibration of cost function weights that balance different objectives. The parametric analysis examined controller sensitivity across multiple weight parameters governing trajectory tracking, control effort, and dynamic response characteristics. The investigation focused on perpendicular deviation weights ($W_\perp$), parallel deviation weights ($W_\parallel$), control input penalties ($R_{acc}$, $R_{steer}$), and rate-of-change penalties ($R_{dacc}$, $R_{dsteer}$) through systematic parameter variation while maintaining a fixed experimental setup. In this analysis, the reference trajectory was recorded from an accurate and nonlinear model of TU Delft Toyota Prius (\cite{Spahn_gym_envs_urdf}) performing a left-turn intersection maneuver, ensuring realistic scenario. Each parameter variation maintained fixed experimental conditions including prediction horizon (13 steps), maximum steering angle (30°), and acceleration/deceleration limits ($-10$ m/s$^2$, $2.0$ m/s$^2$,) to isolate parametric effects. It is worth noting that these settings are selected as proof of concept and may not be used directly for real-world applications. The complete sensitivity analysis, including quantitative trajectory deviation metrics and kinematic profiles, is presented in Appendix C. For conciseness, this section presents the principal findings and parameter thresholds.

The analysis revealed several important parametric relationships for controller calibration. Perpendicular deviation weights demonstrated that sufficiently high values ensure accurate trajectory following. Values of $W_\perp = 20$ and above provided consistent tracking performance, while increases beyond this range did not lead to significant differences in trajectory following and speed-acceleration dynamics. Parallel deviation weights exhibited a trade-off between trajectory adherence and longitudinal control flexibility. Low values of $W_\parallel$ restricted speed adherence via sufficient and timely acceleration and resulted in sluggish vehicle response. However, while higher values enabled smoother acceleration and deceleration transitions, excessive $W_\parallel$ values compromised trajectory tracking as $W_\perp$ importance diminished. Based on the analysis, values between 1 and 5 proved appropriate for $W_\parallel$ in the tested scenarios.

Control input weights ($R_{acc}$, $R_{steer}$) demonstrate the fundamental trade-off between maneuver smoothness and controller responsiveness. Extreme $R_{acc}$ values prevent adequate acceleration/deceleration, creating safety and efficiency concerns. Values between 0.1 and 1 provide appropriate $R_{acc}$ performance. Control rate penalties $R_{d_{acc}}$ and $R_{d_{steer}}$ showed influence on both trajectory deviation and dynamic response and comfort characteristics. Low values of $R_{d_{acc}}$ produced abrupt acceleration changes, which could lead to passenger discomfort, while excessive values prevented the vehicle from utilizing its kinematic capabilities. Values between 5 and 10 demonstrated appropriate balance for $R_{d_{acc}}$. Different values of $R_{d_{steer}}$ did not show notable impact on kinematic profiles but its higher values resulted in more accurate lateral trajectory tracking. Values higher than 1 showed to be appropriate for this parameter.

This parametric sensitivity analysis confirms the importance of balanced parameter calibration in the proposed MPC controller and validates the framework's flexibility in accommodating different performance objectives. While detailed field calibration is outside the scope of this study, the results provide a solid basis for choosing controller parameters in the validation experiments that follow. The parameter values used in all subsequent experiments are listed in Table \ref{tab:mpc_parameters}.
\begin{table}[t]
 \caption{Selected parameters for the motion controller evaluation experiments}
  \centering
  \begin{tabular}{lll}
    \toprule
    Parameter     & Description     & Value (unit) \\
    \midrule
    \(H\) & Horizon length & 13 (steps) \\
    \(W_{\perp}\) & Weight for perpendicular deviation & 20.0 \\
    \(W_{\parallel}\) & Weight for parallel deviation & 1.0 \\
    \(R_{acc, steer}\) & Control input weights & [0.1, 0.01] \\
    \(Rd_{acc, steer}\) & Weights for change in control input & [10, 1.0] \\
    \(Q_{v, \theta}\) & Weights for speed and yaw & [0.0, 0.5] \\
    \(Q_{f}\) & Final state weights & [1.0, 1.0, 0.0, 0.5] \\
    \(Max_{Steer}\) & Maximum steering angle & 30.0 (degrees) \\
    \(Max_{acc}\) & Maximum acceleration & 2.0 (m/s\(^2\)) \\
    \(Max_{dec}\) & Maximum deceleration & -10 (m/s\(^2\)) \\
    \bottomrule
  \end{tabular}
  \label{tab:mpc_parameters}
\end{table}

\subsubsection{Isolated Scenario - Evaluation of Trajectory Follower}
The isolated scenario experiments are designed to evaluate the fundamental capabilities of the controller in tracking the reference trajectory without the influence of external factors such as other vehicles or dynamic obstacles. The experiment setup includes different scenarios aiming to cover possible maneuvers at prevailing urban junctions, including left turn, crossing, and right turn at intersections, as well as left turn, moving through, and U-turn at roundabouts. The trajectory used as the reference is generated by the global planner using the motion primitives achieved by modeling an accurate and nonlinear model of the TU Delft Toyota Prius vehicle (\cite{Spahn_gym_envs_urdf}).  

\begin{figure}[htbp]
    \centering
    \begin{subfigure}[b]{0.45\textwidth}
        \centering
        \includegraphics[width=\textwidth]{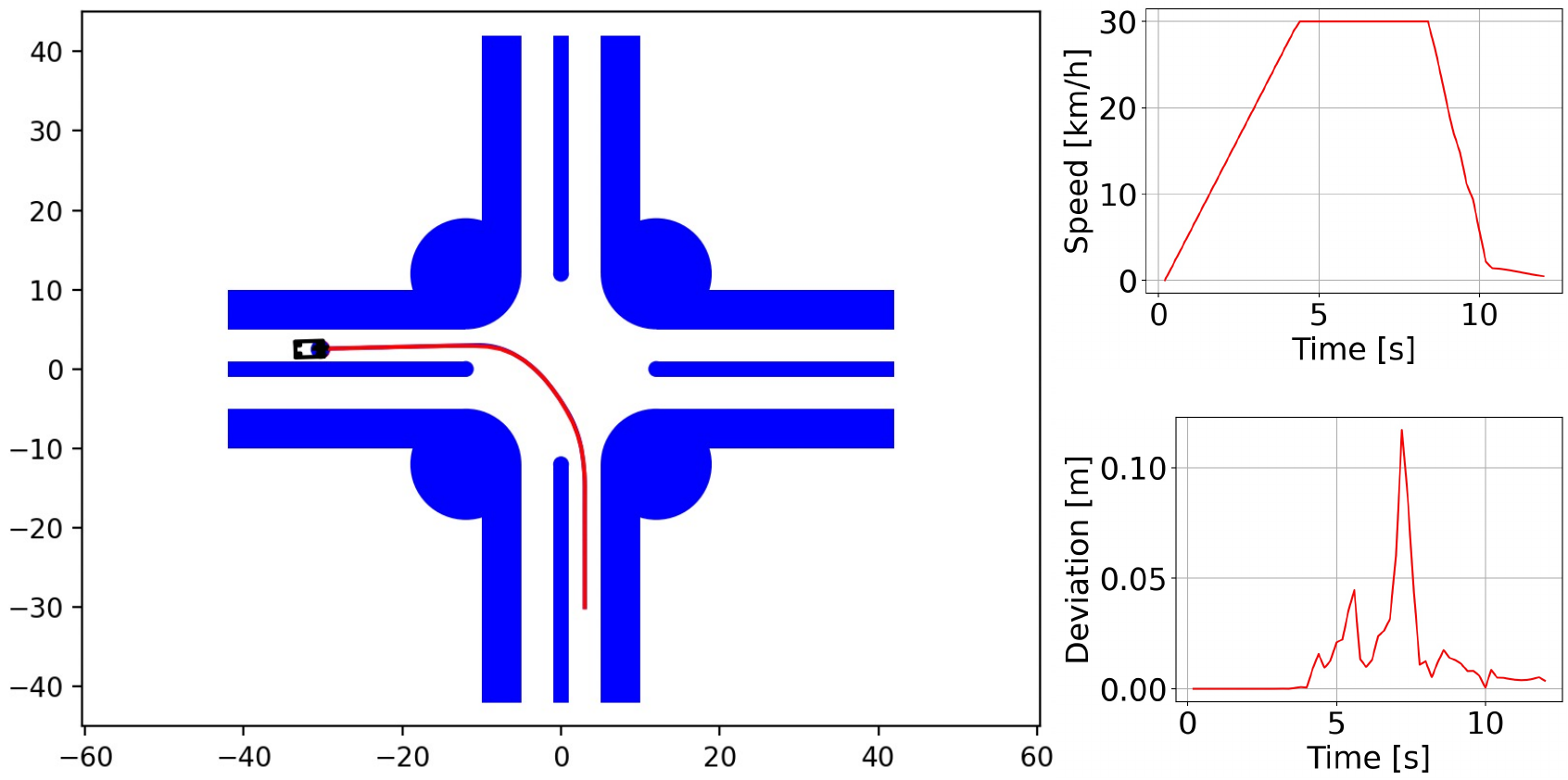}
        \caption{left-turn - Intersection}
        \label{fig:sub1}
    \end{subfigure}
    \hfill
    \begin{subfigure}[b]{0.45\textwidth}
        \centering
        \includegraphics[width=\textwidth]{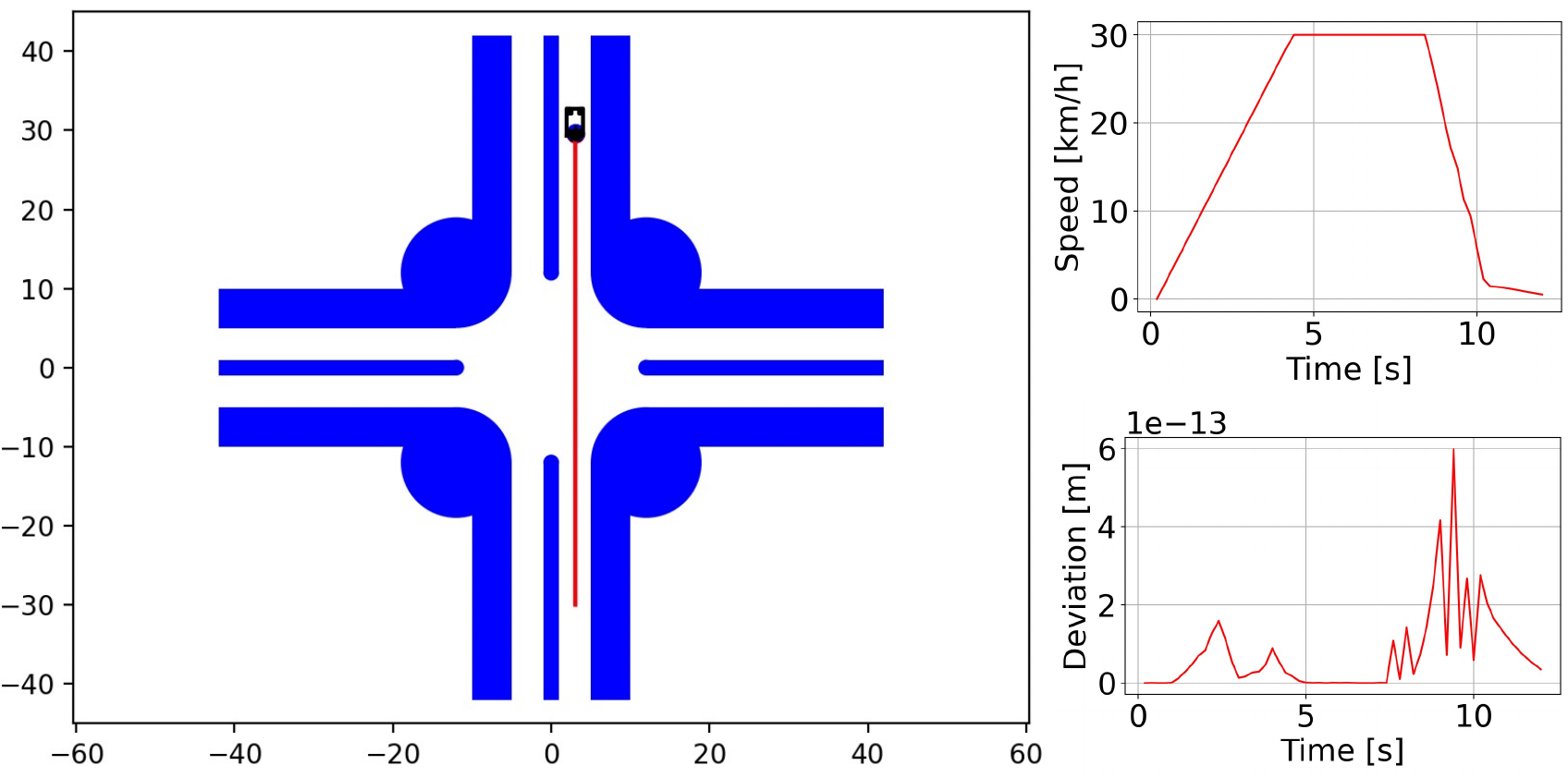}
        \caption{Crossing - Intersection}
        \label{fig:sub2}
    \end{subfigure}

    \vspace{0.3cm}

    \begin{subfigure}[b]{0.45\textwidth}
        \centering
        \includegraphics[width=\textwidth]{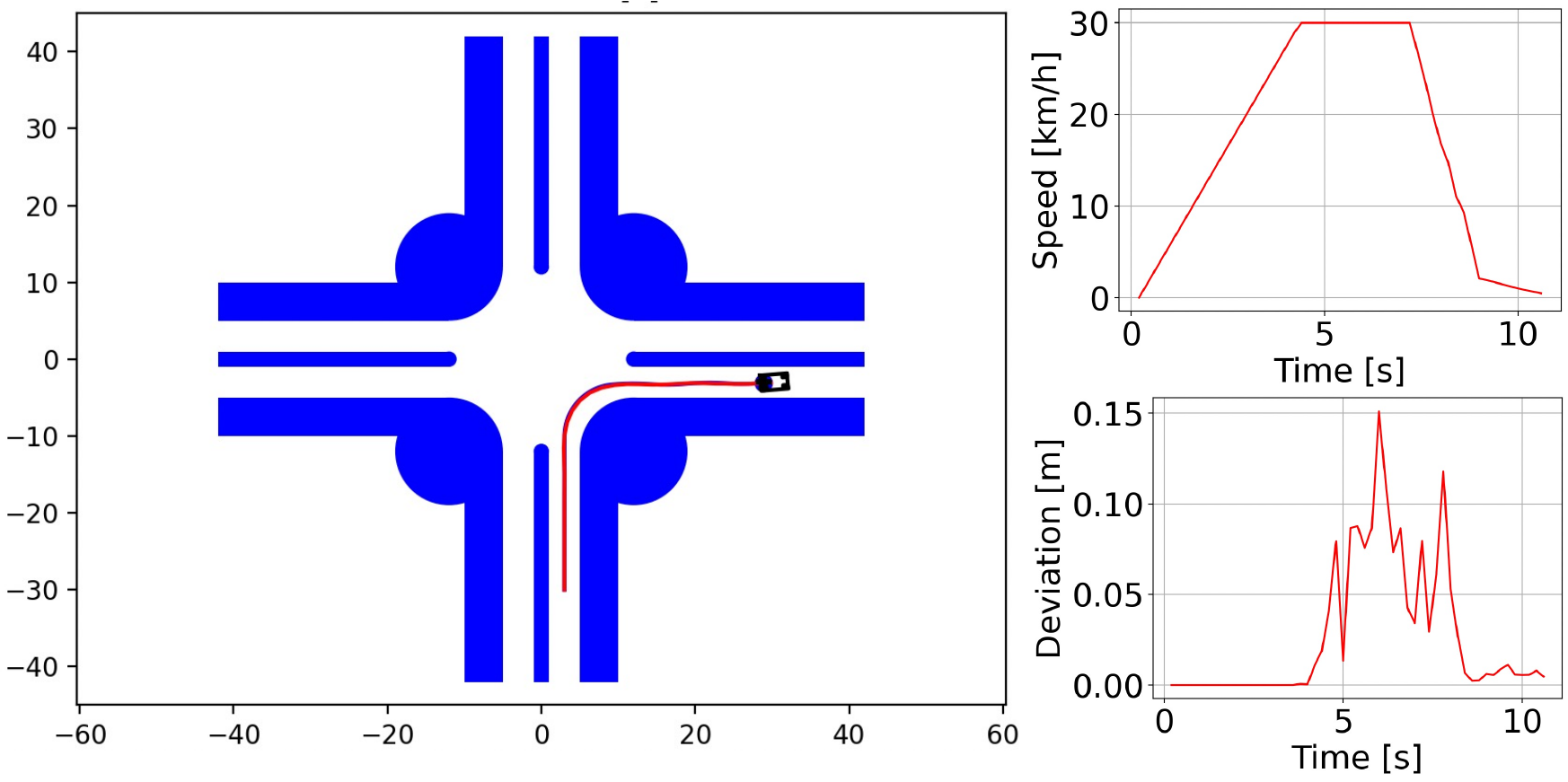}
        \caption{right-turn - Intersection}
        \label{fig:sub3}
    \end{subfigure}
    \hfill
    \begin{subfigure}[b]{0.45\textwidth}
        \centering
        \includegraphics[width=\textwidth]{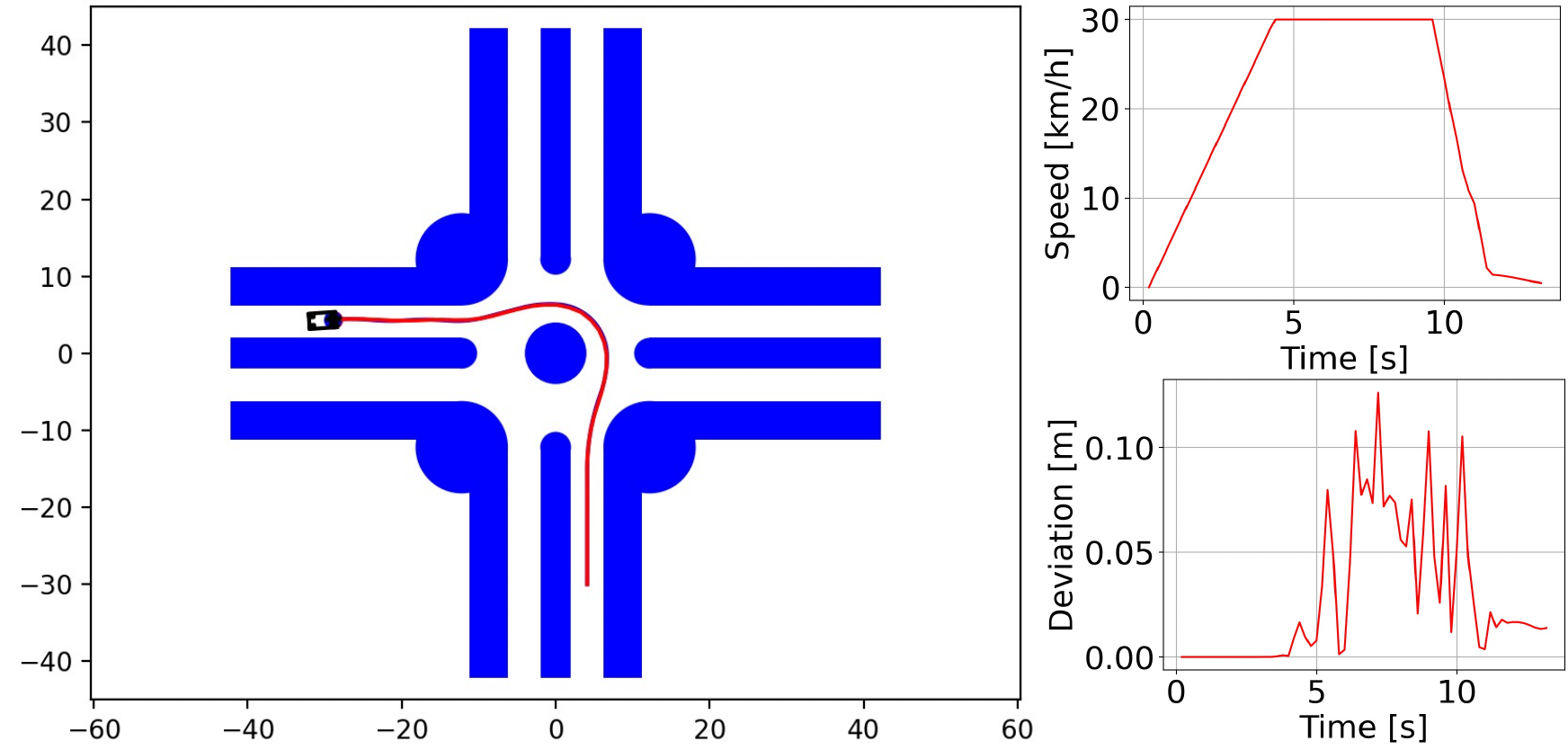}
        \caption{left-turn - Roundabout}
        \label{fig:sub4}
    \end{subfigure}

    \vspace{0.3cm}

    \begin{subfigure}[b]{0.45\textwidth}
        \centering
        \includegraphics[width=\textwidth]{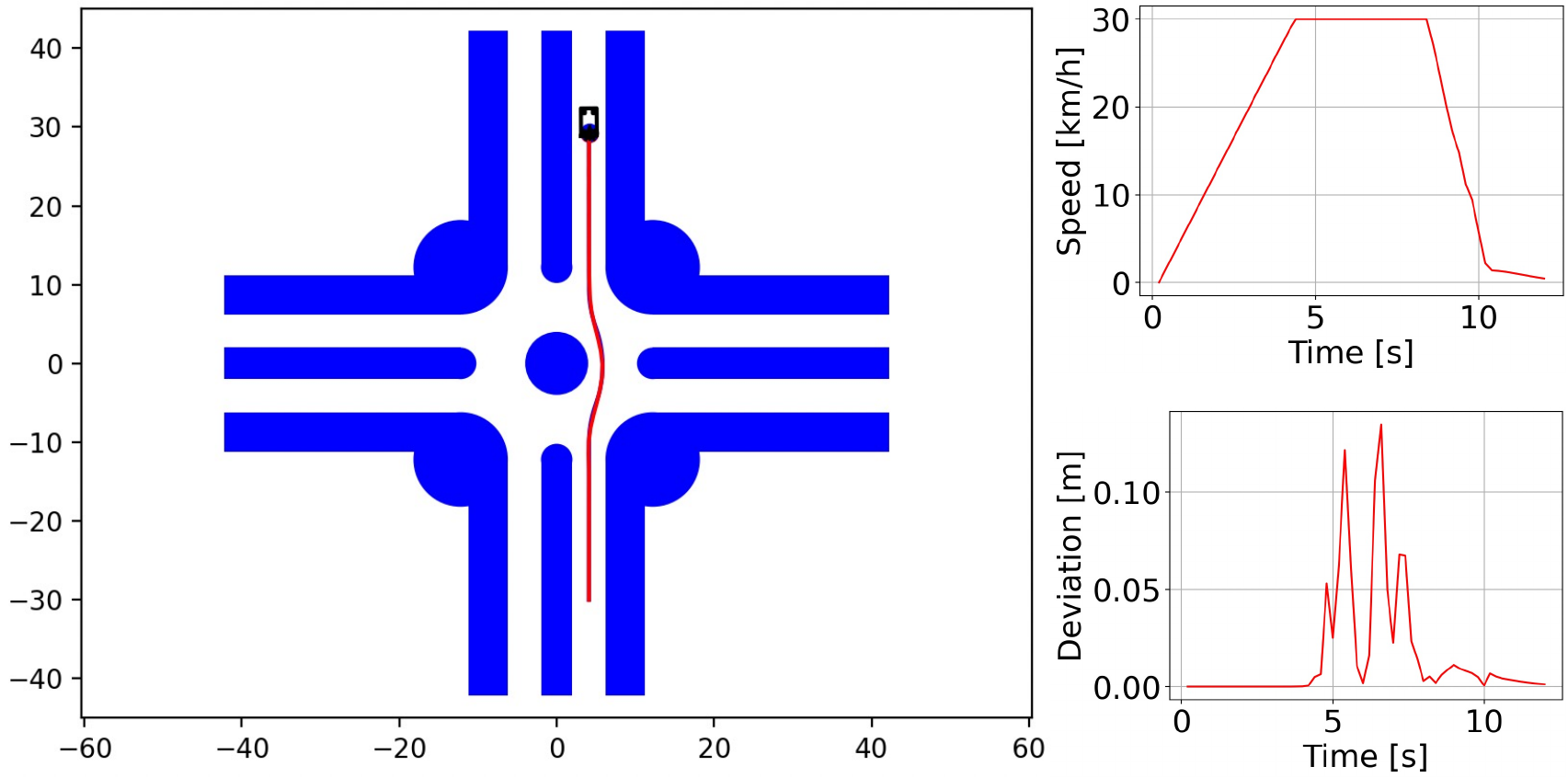}
        \caption{Move through - Roundabout}
        \label{fig:sub5}
    \end{subfigure}
    \hfill
    \begin{subfigure}[b]{0.45\textwidth}
        \centering
        \includegraphics[width=\textwidth]{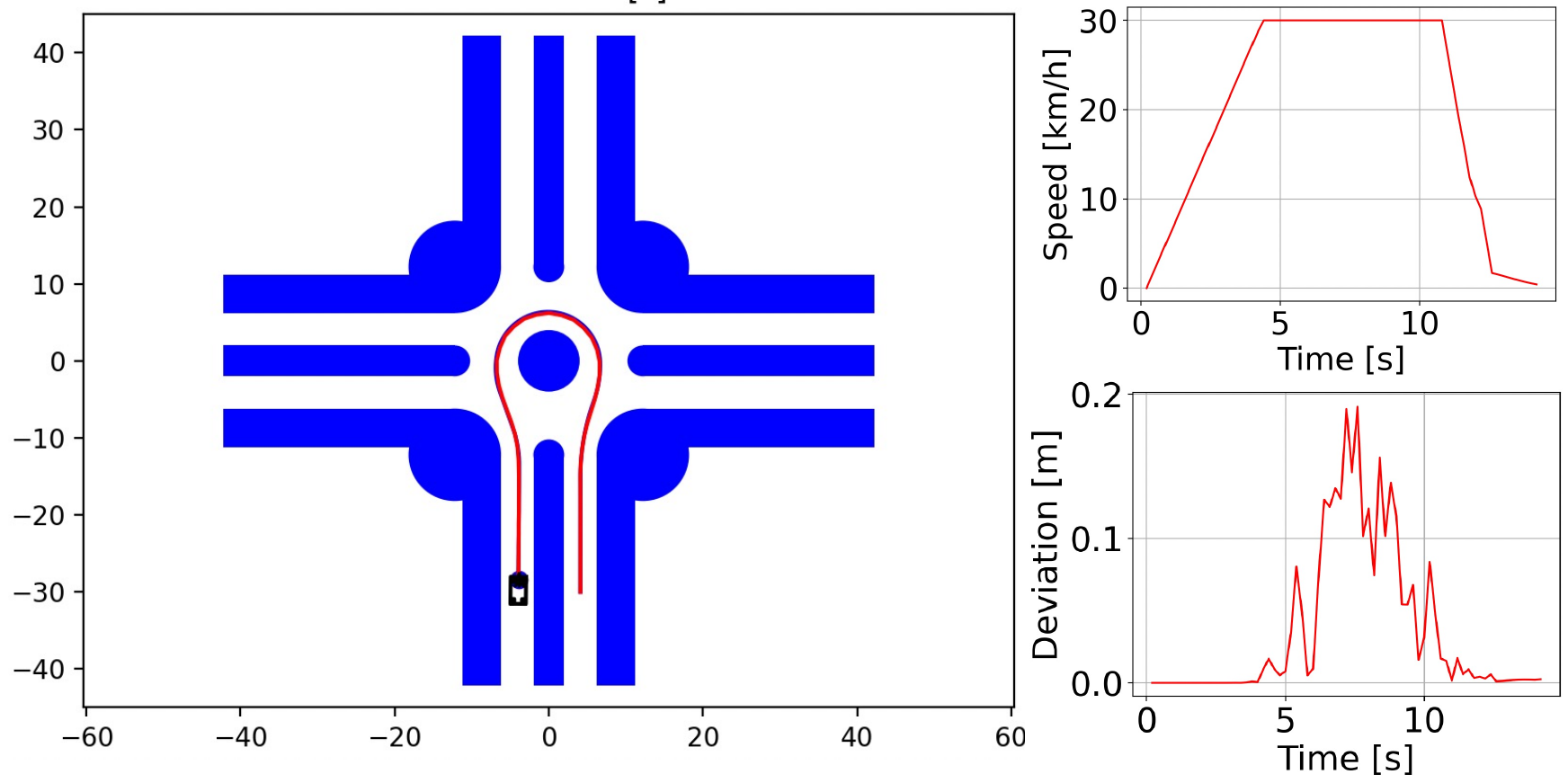}
        \caption{U-turn - Roundabout}
        \label{fig:sub6}
    \end{subfigure}

    \caption{Performance evaluation of the controller in isolated scenarios}
    \label{fig:scenarios_single_agent}
\end{figure}

Figure \ref{fig:scenarios_single_agent} illustrates the experimental validation results. Each sub-figure comprises three distinct plots. The left plot provides a scenario snapshot in which the vehicle, depicted as a black rectangle, is positioned in the goal area after following the reference trajectory. Within these plots, the reference trajectory is shown as a blue curve, while the executed trajectory is marked in red, enabling a clear visual comparison of intended and actual paths. The top-right plot displays the speed profile of the vehicle from the initial to the goal position during the simulation, which serves to assess the controller’s capability to regulate and adhere to the desired speed while tracking the reference trajectory. The bottom-right plot presents the deviation from the reference trajectory at each simulation timestep, quantified in meters. Together, these plots offer both visual and quantitative evaluations of the controller’s effectiveness in accurately tracking the reference trajectory.

Velocity profiles in Figure \ref{fig:scenarios_single_agent} demonstrate appropriate kinematic transitions, with smooth acceleration, steady cruising at the desired speed of 30 km/h, and controlled deceleration as the vehicle approaches the goal states. Quantitative analysis of the deviation from the reference trajectory reveals maximum discrepancies of 0.1 to 0.2 meters, which is negligible relative to the 25-meter intersection diameter (0.4–0.8\%). that the reference (global) trajectory is generated using a high-fidelity nonlinear vehicle model, whereas the controller operates based on a linearized bicycle model. The observed results demonstrate that the linearized bicycle model is sufficiently accurate for capturing the essential dynamics of the nonlinear system at low operational speeds. This finding underscores the effectiveness of the linearized approach for real-time control applications, while maintaining computational tractability and ensuring reliable trajectory tracking performance. Further improvements in performance may be attainable through further parametric tuning of the optimal control formulation.

\subsubsection{Multi-agent Scenarios}
Multi-agent scenarios extend the isolated scenarios by introducing two additional vehicles that interact with the ego vehicle at an urban junction. The experiments assess controller performance under interactive conditions and validate decentralized collision avoidance algorithms through emergent behaviors from independent decision-making processes.

\begin{figure}[t]
    \centering
    \begin{subfigure}[b]{0.45\textwidth}
        \centering
        \includegraphics[width=\textwidth]{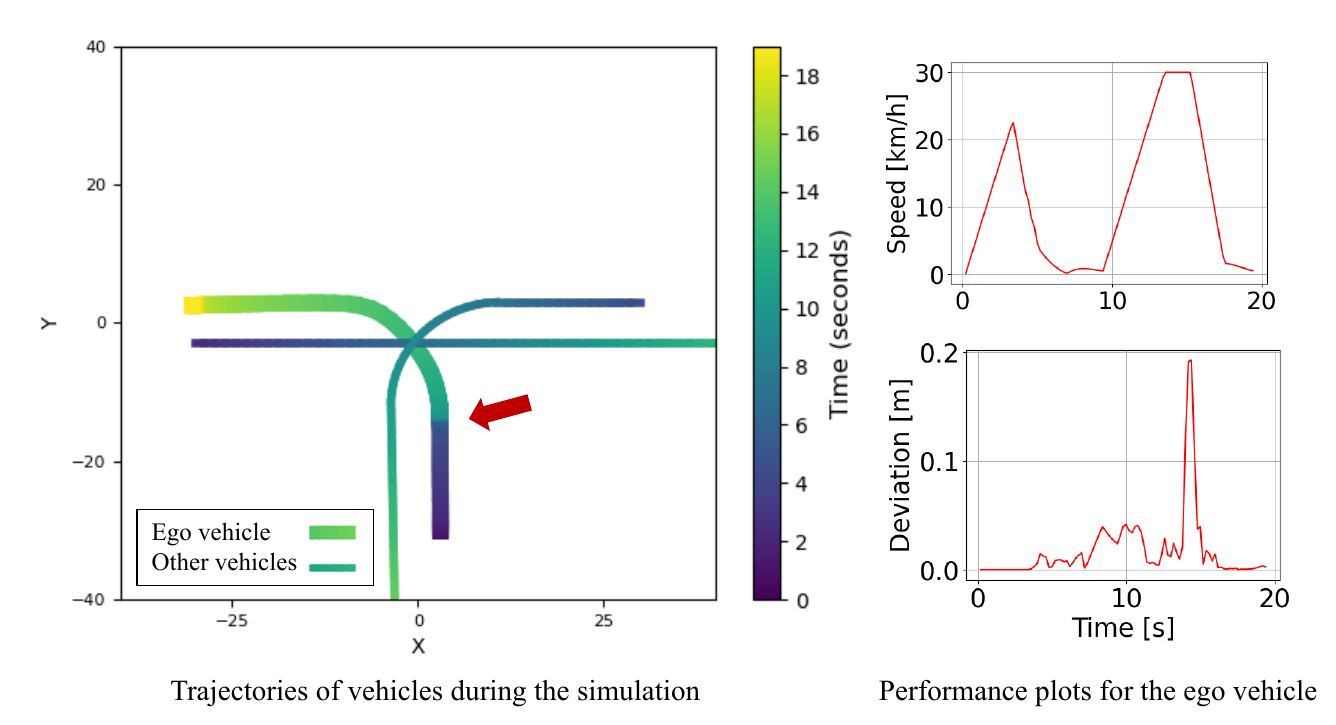}
        \caption{left-turn - Intersection}
        \label{fig:sub1}
    \end{subfigure}
    \hfill
    \begin{subfigure}[b]{0.45\textwidth}
        \centering
        \includegraphics[width=\textwidth]{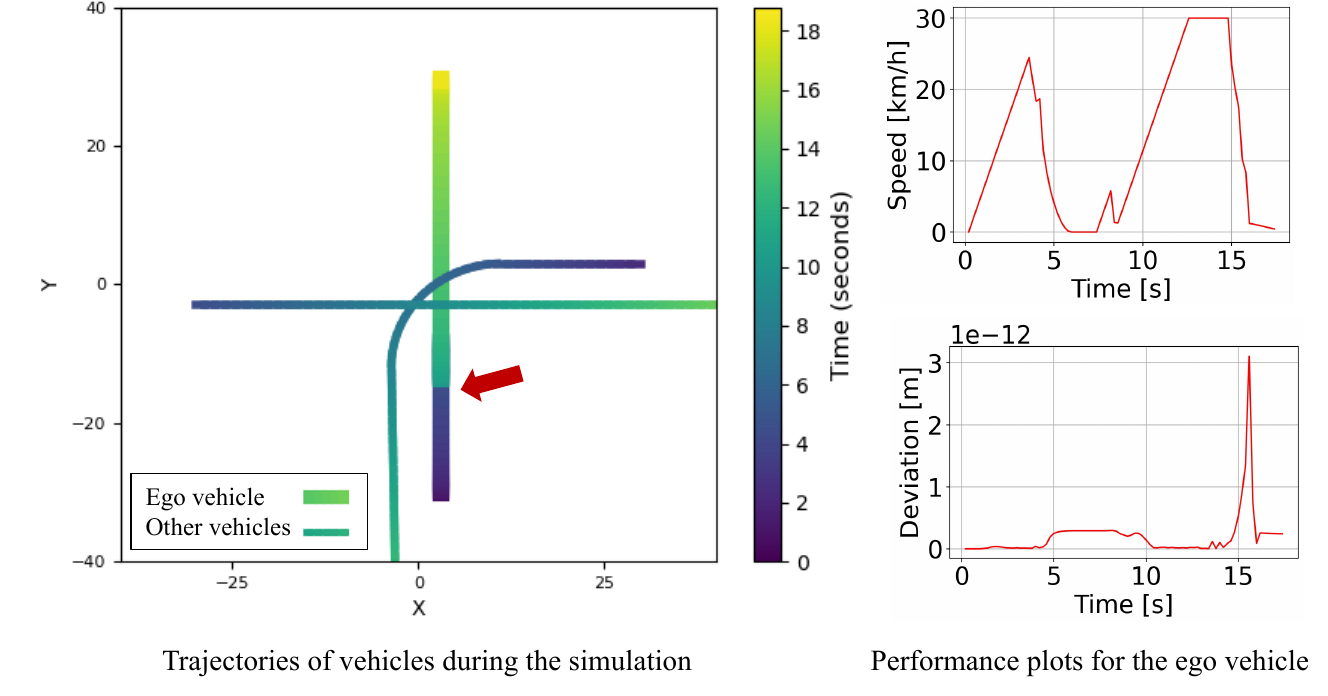}
        \caption{crossing - Intersection}
        \label{fig:sub2}
    \end{subfigure}

    \vspace{0.3cm}

    \begin{subfigure}[b]{0.45\textwidth}
        \centering
        \includegraphics[width=\textwidth]{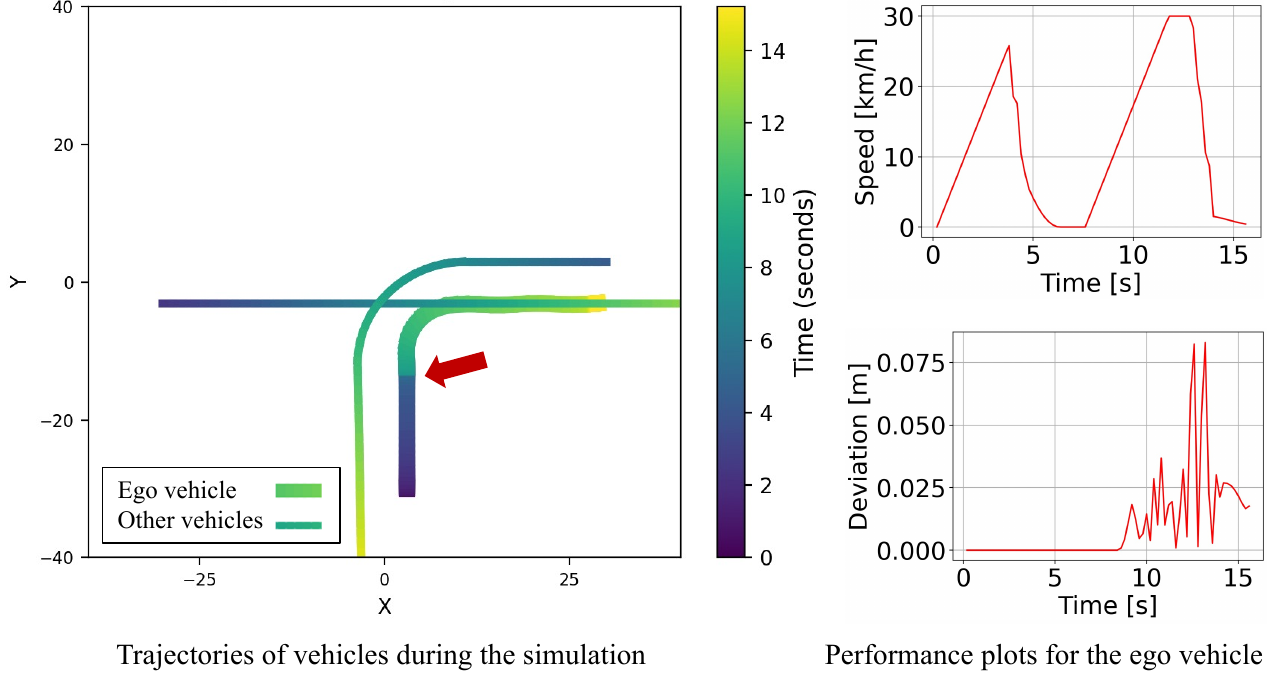}
        \caption{right-turn - Intersection}
        \label{fig:sub3}
    \end{subfigure}
    \hfill
    \begin{subfigure}[b]{0.45\textwidth}
        \centering
        \includegraphics[width=\textwidth]{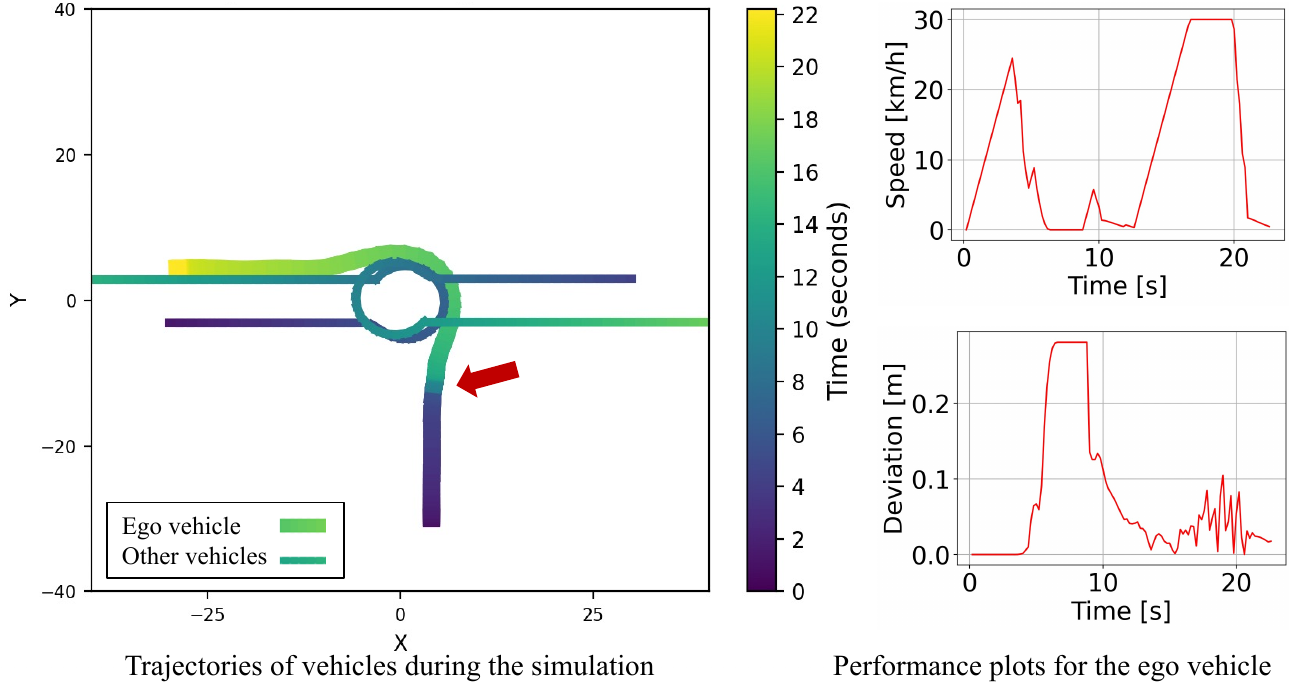}
        \caption{Left-turn - Roundabout}
        \label{fig:sub4}
    \end{subfigure}

    \vspace{0.3cm}

    \begin{subfigure}[b]{0.45\textwidth}
        \centering
        \includegraphics[width=\textwidth]{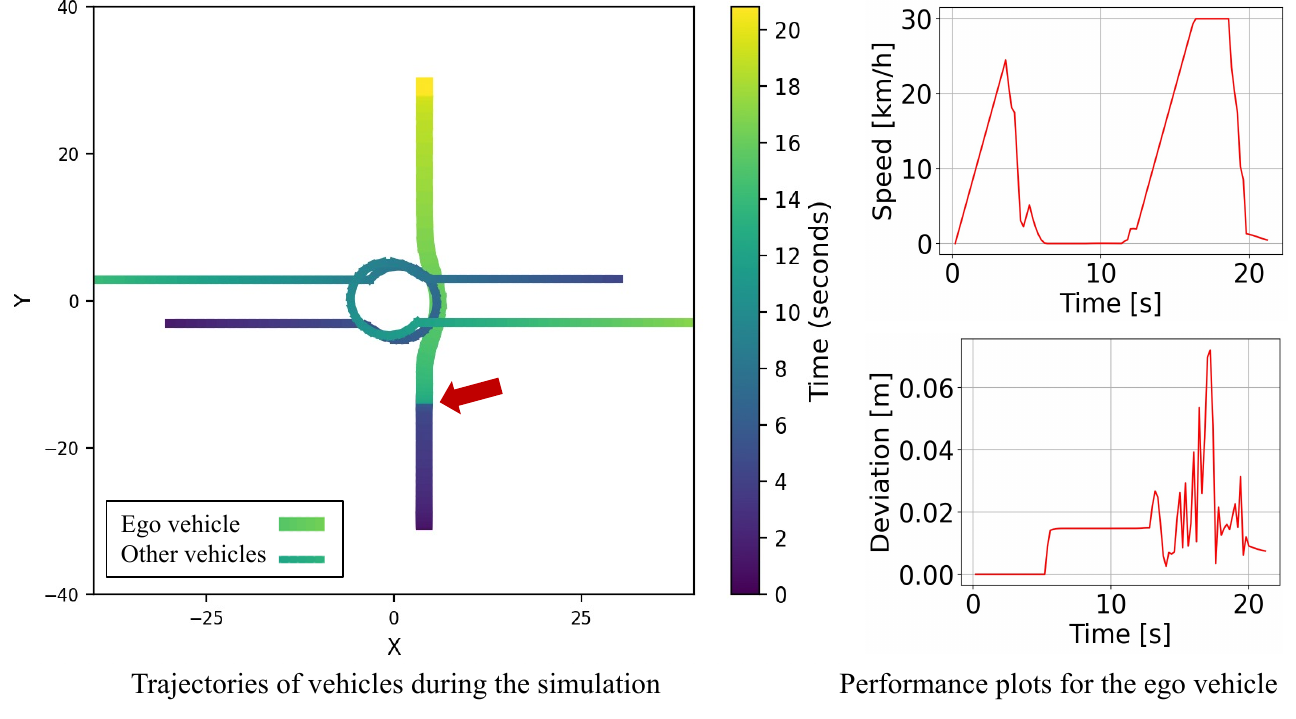}
        \caption{Move through - Roundabout}
        \label{fig:sub5}
    \end{subfigure}
    \hfill
    \begin{subfigure}[b]{0.45\textwidth}
        \centering
        \includegraphics[width=\textwidth]{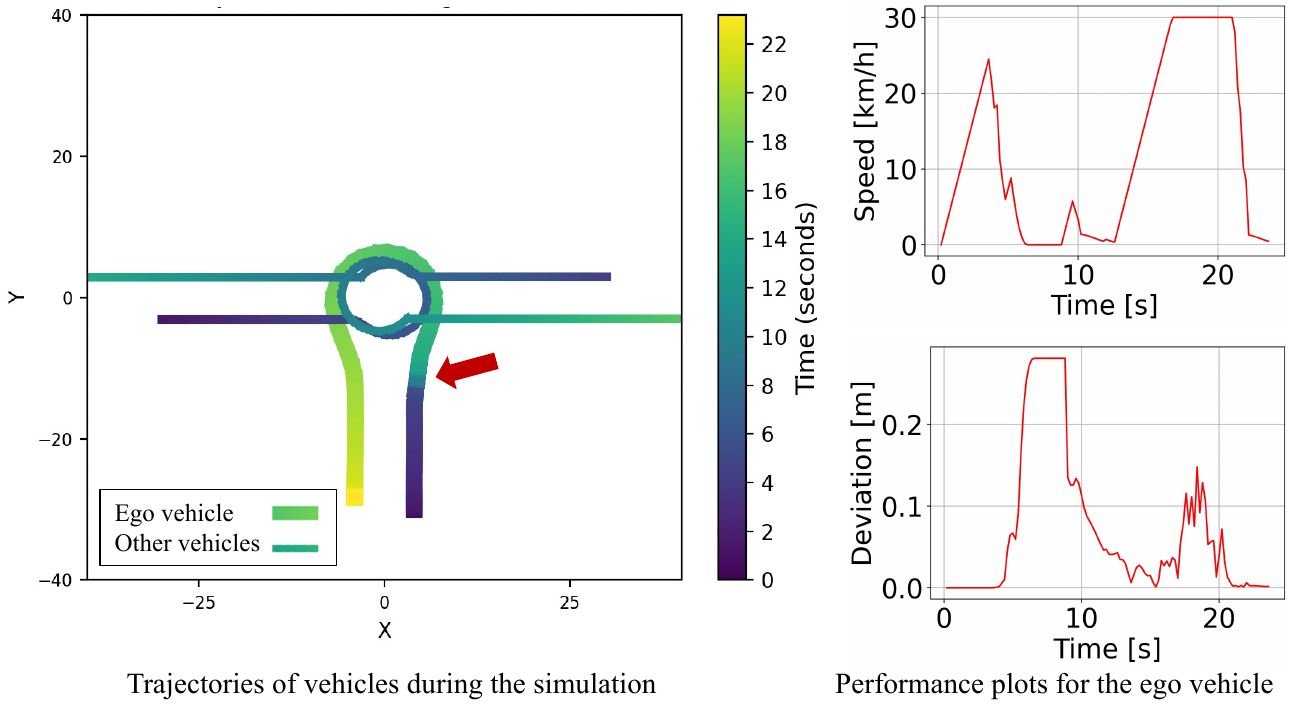}
        \caption{U-turn - Roundabout}
        \label{fig:sub6}
    \end{subfigure}

    \caption{Performance evaluation of the controller in multi-agent scenarios}
    \label{fig:scenarios_multi_agent}
\end{figure}

Figure~\ref{fig:scenarios_multi_agent} presents the evaluation results through temporally color-coded vehicle trajectories and speed and deviation profiles. In each sub-figure, the plot on the left represents the progress of agents' trajectories throughout the simulation, where the ego vehicle trajectory is distinguished by increased line weight to facilitate identification. The color gradient from blue to yellow represents temporal progression throughout each simulation to enable the analysis of spatial-temporal interactions. The intersection of trajectories with different temporal gradients demonstrates successful collision avoidance by indicating that vehicles traverse identical spatial locations at distinct temporal instances. Red arrows indicate collision avoidance intervention points where the ego vehicle has stopped or reduced its speed to prevent potential collisions, demonstrating the controller's conflict detection and response capabilities. The speed-time profiles are depicted on the top right part of each sub-figure to demonstrate the changes of the speed and adherence to desired speed during the interactions. The deviation plots on the bottom right part of the sub-figures demonstrate the divergence from the reference trajectories throughout the interactions to evaluate the performance of the motion controller in tracking the reference trajectory while performing collision-avoidance maneuvers. 

Visual examination of the coloured trajectories confirms the effective collision avoidance performance. The different color gradients at the intersection points of trajectories indicates that vehicles have avoided simultaneous occupation of the same locations. The red arrows identify specific locations where the ego vehicle has detected potential collisions and triggered avoidance responses through deceleration or complete stops. The positioning of these points demonstrates accurate conflict identification considering a safety margin, which can be configured based on specific safety criteria and operational constraints. 

The speed profiles demonstrate controlled velocity management across multi-agent scenarios. During normal operation in the beginning of the scenarios, the ego vehicle maintains smooth accelerations to achieve the reference velocities according to the planned trajectory. Upon detecting potential collisions, the controller executes relatively smooth deceleration to mitigate collision risk by reducing speed. Following successful conflict avoidance, the vehicle smoothly accelerates to resume planned velocities. Finally, upon arrival at the goal area, the vehicle decelerates to stop at the final state. These velocity transitions occur without excessive acceleration rates, confirming adherence to vehicle dynamic constraints and passenger comfort considerations. The deviation measurements provide quantitative validation of trajectory tracking performance under dynamic conditions. Maximum deviations from reference paths remain consistently below 0.2 meters across all experimental scenarios, demonstrating precise path-following capability despite simultaneous collision avoidance maneuvers. This performance validates the bi-level architecture's effectiveness in maintaining trajectory accuracy while executing real-time safety maneuvers. Supplementary video, available on \url{https://github.com/SaeedRahmani/MPC_for_AV_at_Intersection} provides additional validation of real-time controller performance in multi-agent settings.

These experimental findings demonstrate the effectiveness of the proposed bi-level framework for vehicular maneuvering at unsignalized intersections operating under decentralized conditions. The global planner consistently generates kinematically feasible reference trajectories across diverse intersection configurations, while the motion controller maintains trajectory accuracy during real-time multi-agent interactions. The collision avoidance implementation, although employing simplified assumptions, successfully prevents conflicts through independent vehicle decision-making without requiring inter-vehicle communication or centralized coordination. 

\subsection{Computational Efficiency}
The computational efficiency of the proposed bi-level framework is assessed through a comparative runtime analysis against standard monolithic MPC controllers with varying prediction horizon lengths, as well as the optimal control method described by (\cite{zhao2023microscopic}). Here, monolithic MPC refers to a single optimization formulation that jointly handles trajectory tracking and collision avoidance constraints within an extended receding horizon. Monolithic MPC implementations require substantially extended horizons ($H = 25$ to $100$ steps) to ensure reliable collision avoidance in urban scenarios (\cite{yu2021model}), which leads to a quadratic increase in decision-space dimensionality and computational cost. By decoupling collision avoidance from the receding horizon controller, the proposed architecture enables the use of much shorter horizons without compromising safety. This comparison quantifies the computational benefits of the bi-level design relative to these baseline approaches. Table~\ref{tab:_framework_runtime} summarizes the results of this comparative analysis. 

The proposed framework achieves real-time performance, requiring only 0.096 seconds per simulation step and completing an 18-second scenario in 7.26 seconds. In contrast, pure MPC implementations with extended horizons exhibit substantially higher computational demands, with runtimes increasing from 0.176 seconds per step for $H=25$ to 0.84 seconds per step for $H=100$. The optimal control method by \cite{zhao2023microscopic} is significantly more computationally intensive, with a total runtime of 150 seconds for the same scenario. It is worth noting our experiments indicate that a pure MPC controller with $H=25$ not only doubles the runtime relative to the proposed framework, but also fails to consistently guarantee collision-free operation, highlighting the limitations of monolithic approaches with shorter horizons. Moreover, the approach proposed by \cite{zhao2023microscopic} relies on the sharing of utility functions among all vehicles to enable centralized coordination and thus avoiding collisions through joint optimization. However, such centralized controllability comes at the cost of significantly increased computational burden as the complexity of the joint optimization scales rapidly with the number of interacting agents.

These results demonstrate that the proposed bi-level framework achieves notable computational savings, making it well-suited for applications in multi-vehicle scenarios. The next section presents a brief discussion on the framework’s capability to simulate controlled collision scenarios, highlighting its potential for safety analysis studies.
\begin{table}[h!]
    \centering
    \small
    \caption{Comparison of absolute run times and percentage \emph{increase} relative to the proposed controller}
    \label{tab:framework_runtime}
    \resizebox{\textwidth}{!}{%
        \begin{tabular}{p{\dimexpr\textwidth-10cm} *{5}{p{2cm}}} \hline
            Run time (sec) 
                & H=14 \newline Proposed 
                & H=25 \newline Pure MPC 
                & H=50 \newline Pure MPC 
                & H=100 \newline Pure MPC 
                & \textcite{zhao2023microscopic} \\ \hline
                
            % ---- total run‑time ----
            Total run‑time 
                & 7.26 
                & 16.95 
                & 32.61 
                & 81.49 
                & $\sim$150 \\ 
            \% increase vs proposed\footnotesize\hspace{0.1em}
                & -- 
                & 133\,\% 
                & 349\,\% 
                & 1\,023\,\% 
                & 1\,966\,\% \\ \hline
                
            % ---- per‑step run‑time ----
            Run‑time per simulation step 
                & 0.096 
                & 0.176 
                & 0.340 
                & 0.840 
                & -- \\ 
            \% increase vs proposed 
                & -- 
                & 83\,\% 
                & 254\,\% 
                & 775\,\% 
                & -- \\ \hline
                
            % ---- remaining metrics ----
            Simulation duration (sec) 
                & 18 & 18 & 18 & 18 & 150 \\ \hline
            Real time 
                & Yes 
                & Yes* 
                & No 
                & No 
                & No \\ \hline
                
            \multicolumn{6}{p{\textwidth}}{%
                *Although the run‑time needed for solving the MPC can be handled in real time, in cases where a new reference trajectory is required the total computation time may exceed the simulation time step. 
            } \\ 
            \multicolumn{6}{p{\textwidth}}{%
            \textbf{Note:} run times were recorded on a MacBook Pro (M1, 16GB RAM) using CPU only with the power adapter connected. 
            Percentage increase is calculated as $(\text{other}-\text{proposed})/\text{proposed}\times100$.%
} \\ \hline
        \end{tabular}
    }
    \label{tab:_framework_runtime}
\end{table}

\subsection{Accidents and Safety Analysis}
The proposed framework supports the analysis of safety-critical interactions by explicitly modeling perception constraints, kinematic limits, and finite response times by enjoying a decentralized control structure. As a result, unlike methods that assume perfect information sharing and enforce collision-free operation, this approach permits emergent behaviors, including collisions, that arise naturally from realistic sensing and actuation limitations.

\begin{figure}
  \centering
  % First Image
  \begin{subfigure}[b]{0.48\textwidth}
    \includegraphics[width=\textwidth]{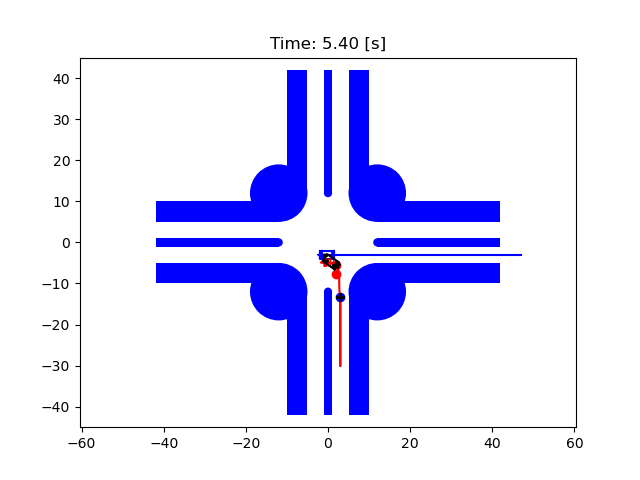}
    \caption{}
    \label{fig:crash_late_detection}
  \end{subfigure}
  \hfill
  % Second Image
  \begin{subfigure}[b]{0.48\textwidth}
    \includegraphics[width=\textwidth]{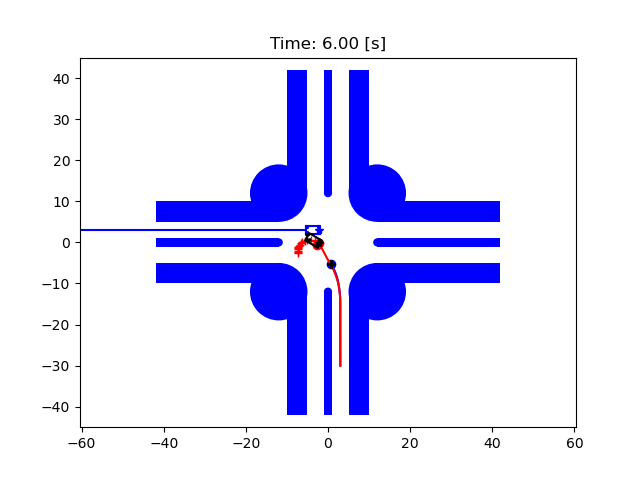}
    \caption{}
    \label{fig:crash_high_speed}
  \end{subfigure}
  
  \caption{Two manually-crafted scenarios to generate accidents (a) late detection of the other vehicle; (b) Extreme speed of the interacting vehicle (200km/h)}
\label{fig:accident_scenarios}
\end{figure}

To demonstrate this capability, we construct two experimental scenarios that explore distinct mechanisms leading to unsafe outcomes. In the first configuration (Figure~\ref{fig:crash_late_detection}), the ego vehicle accelerates northbound with a cruising speed of \(50~\mathrm{km\,h^{-1}}\) (\(13.9~\mathrm{m\,s^{-1}}\)) while a second vehicle approaches from the west. Here, the ego vehicle’s perception is artificially limited to a 10-meter detection range, representing degraded sensing capability or environmental occlusion. A reaction latency of \(0.5~\mathrm{s}\) is simulated by withholding sensor updates from the MPC during that interval. During this delay, the vehicle covers \(d_{\text{react}} = v_{\text{des}}\Delta t \approx 13.9 \times 0.5 = 7.0~\mathrm{m}\). Subsequent braking at the maximum deceleration of \(-10~\mathrm{m\,s^{-2}}\) requires \(d_{\text{brake}} = v_{\text{des}}^{2} / \bigl(2|a_{\text{max}}|\bigr) \approx 9.6~\mathrm{m}\). The total stopping distance, \(
d_{\text{stop}} = d_{\text{react}} + d_{\text{brake}} \approx 16.6~\mathrm{m},
\) exceeds the available distance required to avoid collision. Consequently, the ego vehicle enters the intersection, leading to a collision due to the kinematic constraints. This is also illustrated in Figure~\ref{fig:crash_late_detection}, where the two vehicles has crashed into each other. The red curve in this figure shows the traversed trajectory by the ego vehicle, while the blue rectangle shows the other crossing vehicle. 

In the second scenario (Figure~\ref{fig:crash_high_speed}), the ego vehicle has a $50\mathrm{m}$ detection range; however, the conflicting vehicle approaches at the speed of $200\mathrm{km/h}$ ($55.6,\mathrm{m,s^{-1}}$). Upon detection, the approaching vehicle reaches the conflict point in roughly $t_{\text{close}} = 50 / 55.6 \approx 0.9,\mathrm{s}$. The ego cruises at $50,\mathrm{km,h^{-1}}$ ($13.9,\mathrm{m,s^{-1}}$), experiencing a $0.5,\mathrm{s}$ perception–reaction delay, needs $d_{\text{react}} = 13.9 \times 0.5 \approx 7.0,\mathrm{m}$ before braking begins. With the remaining $0.4,\mathrm{s}$ and a deceleration limit of $a_{\max} = -10,\mathrm{m,s^{-2}}$, the ego can shed only $d_{\text{brake}} \approx 4.8,\mathrm{m}$, giving a total clearance of $d_{\text{total}} \approx 11.8,\mathrm{m}$. However, the baseline stopping distance at this speed and latency is $d_{\text{stop}} = 7.0 + 9.6 = 16.6,\mathrm{m}$, so the available distance falls short by nearly $5,\mathrm{m}$. Consequently, even with perfect detection and an optimal response, the spatial–temporal margin is insufficient and the feasible set for collision avoidance is empty.

These controlled experiments validate the framework's ability to generate collision scenarios through realistic inter-vehicle dynamics and sensor limitations. The results confirm that safety-critical events can be systematically studied within the proposed mathematical structure, providing a foundation for quantitative risk assessment in intersection environments. Although these demonstrations serve primarily as proof-of-concept illustrations, they establish the framework's potential for safety analysis applications where accident modeling is essential for understanding traffic system vulnerabilities and evaluating countermeasure effectiveness.
\section{Conclusions}\label{sec:conclusion}
This study presents a bi-level mathematical framework for modeling vehicular maneuvers at urban junctions, aiming to address critical limitations in existing approaches for mixed-autonomy traffic analysis. The key theoretical contribution lies in decomposing the complex multi-vehicle trajectory optimization problem into computationally tractable subproblems while preserving kinematic feasibility and enabling decentralized operation.  Unlike existing methods that assume cooperation or central control, our approach operates under realistic decentralized conditions, enabling emergence of safety-critical interactions essential for AV system validation.

Experimental validation demonstrates the framework generates kinematically feasible trajectories across diverse intersection topologies with real-time computational performance. The approach achieves linear scalability with vehicle numbers and maintains trajectory tracking accuracy below 0.2 meters while executing collision avoidance maneuvers. The theoretical advances position this framework as a significant methodological contribution for intersection modeling, particularly valuable for AV development where detailed analysis of vehicular interactions under realistic operational constraints is essential for safe deployment in mixed-traffic environments. While the proposed framework demonstrates significant theoretical and computational advantages, there are several areas where it can be improved. Empirical calibration using real-world trajectory data and scalability assessment under high-density traffic conditions represent essential next steps. The current collision avoidance formulation exhibits conservative characteristics that may limit throughput efficiency. Advanced prediction and negotiation mechanisms, potentially incorporating machine learning methodologies, offer promising directions for addressing this limitation while maintaining safety guarantees. Integration of stochastic uncertainty modeling could enhance the framework's applicability to real-world deployment scenarios with sensor noise and environmental variability.

\section*{Authors Contributions}
The authors confirm their contribution to the paper as follows: study conception and design: S. Rahmani, S.C. Calvert, B. van Arem; analysis and interpretation of results: S. Rahmani, S.C. Calvert, B. van Arem; draft manuscript preparation: S. Rahmani, S.C. Calvert; manuscript draft revision:  S.C. Calvert, B. van Arem. All authors reviewed the results and approved the final version of the manuscript.

\printbibliography

\newpage
\appendix
\section{Appendix A: Global Planner Algorithm}
The pseudo-code for the developed graph construction and optimized A* search algorithm is depicted in Algorithm \ref{alg:A*}.

\begin{algorithm}[h]
\caption{Proposed Graph Construction and Search Algorithm}
\begin{algorithmic}[1]
\label{alg:A*}
\REQUIRE start node $start$, goal node $goal$, neighbor function $neighbors$, heuristic function $h$, and actual cost function $g$
\ENSURE path from $start$ to $goal$

\STATE openList $\leftarrow$ empty priority queue
\STATE closedList $\leftarrow$ empty set
\STATE $g[start] \leftarrow 0$
\STATE $f[start] \leftarrow h(start)$
\STATE add $start$ to openList with priority $f[start]$

\WHILE{openList is not empty}
    \STATE current $\leftarrow$ node in openList with the lowest $f$ value
    \IF{current is $goal$}
        \STATE \textbf{return} reconstructPath(cameFrom, current)
    \ENDIF
    \STATE remove current from openList
    \STATE add current to closedList
    
    \FOR{each neighbor in neighbors(current)}
        \STATE tentative\_gScore $\leftarrow g[current] + c(current, neighbor)$
        \IF{neighbor in closedList \AND tentative\_gScore $\geq g[neighbor]$}
            \STATE \textbf{continue}
        \ENDIF
        \IF{neighbor not in openList \OR tentative\_gScore $<$ g[neighbor]}
            \STATE cameFrom[neighbor] $\leftarrow$ current
            \STATE g[neighbor] $\leftarrow$ tentative\_gScore
            \STATE f[neighbor] $\leftarrow$ g[neighbor] + h(neighbor)
            \IF{neighbor not in openList}
                \STATE add neighbor to openList with priority $f[neighbor]$
            \ENDIF
        \ENDIF
    \ENDFOR
\ENDWHILE

\STATE \textbf{return} failure

\STATE
\STATE \textbf{Function reconstructPath(cameFrom, current):}
    \STATE totalPath $\leftarrow$ [current]
    \WHILE{current in cameFrom}
        \STATE current $\leftarrow$ cameFrom[current]
        \STATE prepend current to totalPath
    \ENDWHILE
    \STATE \textbf{return} totalPath
\end{algorithmic}
\end{algorithm}

\newpage
\section{Appendix B: Linearization and Discretization of the Vehicle Dynamics}

The linearization of the nonlinear bicycle model around the operating point $(\bar{\mathbf{x}}, \bar{\mathbf{u}})$ is carried out using a first-order Taylor series expansion, resulting in the state-space matrices as follows:

\begin{equation}
        A = \begin{bmatrix}
    \frac{\partial \dot{x}}{\partial x} & \frac{\partial \dot{x}}{\partial y} & \frac{\partial \dot{x}}{\partial v} & \frac{\partial \dot{x}}{\partial \theta} \\[6pt]
    \frac{\partial \dot{y}}{\partial x} & \frac{\partial \dot{y}}{\partial y} & \frac{\partial \dot{y}}{\partial v} & \frac{\partial \dot{y}}{\partial \theta} \\[6pt]
    \frac{\partial \dot{v}}{\partial x} & \frac{\partial \dot{v}}{\partial y} & \frac{\partial \dot{v}}{\partial v} & \frac{\partial \dot{v}}{\partial \theta} \\[6pt]
    \frac{\partial \dot{\theta}}{\partial x} & \frac{\partial \dot{\theta}}{\partial y} & \frac{\partial \dot{\theta}}{\partial v} & \frac{\partial \dot{\theta}}{\partial \theta}
    \end{bmatrix} = \begin{bmatrix}
    0 & 0 & \cos(\bar{\theta}) & -\bar{v} \sin(\bar{\theta}) \\[6pt]
    0 & 0 & \sin(\bar{\theta}) & \bar{v} \cos(\bar{\theta}) \\[6pt]
    0 & 0 & 0 & 0 \\[6pt]
    0 & 0 & \frac{\tan(\bar{\delta})}{L} & 0
    \end{bmatrix}, \quad 
B = \begin{bmatrix}
    \frac{\partial \dot{x}}{\partial a} & \frac{\partial \dot{x}}{\partial \delta} \\[6pt]
    \frac{\partial \dot{y}}{\partial a} & \frac{\partial \dot{y}}{\partial \delta} \\[6pt]
    \frac{\partial \dot{v}}{\partial a} & \frac{\partial \dot{v}}{\partial \delta} \\[6pt]
    \frac{\partial \dot{\theta}}{\partial a} & \frac{\partial \dot{\theta}}{\partial \delta}
    \end{bmatrix} = \begin{bmatrix}
    0 & 0 \\[6pt]
    0 & 0 \\[6pt]
    1 & 0 \\[6pt]
    0 & \frac{\bar{v}}{L \cos^2(\bar{\delta})}
    \end{bmatrix}
\end{equation}

\begin{equation}
        \mathbf{d} = f(\bar{\mathbf{x}}, \bar{\mathbf{u}}) - A\bar{\mathbf{x}} - B\bar{\mathbf{u}} 
    =
    \begin{bmatrix}
    \bar{v} \sin(\bar{\theta}) \bar{\theta} \\[6pt]
    -\bar{v} \cos(\bar{\theta}) \bar{\theta} \\[6pt]
    0 \\[6pt]
    -\frac{\bar{v} \bar{\delta}}{L \cos^2(\bar{\delta})}
    \end{bmatrix}
\end{equation}

\paragraph{Discretization}
Given a specified sampling interval \(T_s\), and first-order Euler approximations, the state transition matrix \(A_d\) and discrete-time input matrix \(B_d\) are derived through the following transformations:

\begin{equation}
A_d = e^{A T_s} \approx I + A T_s, \quad
B_d = \left( \int_0^{T_s} e^{A \tau} \, d\tau \right) B \approx T_s B
\label{eq:combined_discretization}
\end{equation}

Next, by substituting the continuous-time matrices \(A\) and \(B\) derived above, the following matrices are obtained:
\begin{equation}
A_d = \begin{bmatrix}
1 & 0 & T_s \cos(\bar{\theta}) & -T_s \bar{v} \sin(\bar{\theta}) \\[6pt]
0 & 1 & T_s \sin(\bar{\theta}) & T_s \bar{v} \cos(\bar{\theta}) \\[6pt]
0 & 0 & 1 & 0 \\[6pt]
0 & 0 & \frac{T_s \tan(\bar{\delta})}{L} & 1
\end{bmatrix}, \quad
B_d = \begin{bmatrix}
0 & 0 \\[6pt]
0 & 0 \\[6pt]
T_s & 0 \\[6pt]
0 & \frac{T_s \bar{v}}{L \cos^2(\bar{\delta})}
\end{bmatrix}
\mathbf{d}_d = \begin{bmatrix}
T_s \bar{v} \sin(\bar{\theta}) \bar{\theta} \\[6pt]
- T_s \bar{v} \cos(\bar{\theta}) \bar{\theta} \\[6pt]
0 \\[6pt]
-\frac{T_s \bar{v} \bar{\delta}}{L \cos^2(\bar{\delta})}
\end{bmatrix}
\label{eq:combined_matrices}
\end{equation}

The resulting discrete-time state-space equations are formally represented as follows:
\begin{equation}
    \mathbf{x}(k+1) = A_{d}\mathbf{x}(k) + B_{d}\mathbf{u}(k) + \mathbf{d}_{d}
    \label{eq:discretized_state_space}
\end{equation}

where \(k\) denotes the discrete time step, \(A_{d}\), \(B_{d}\), and \(\mathbf{d}_{d}\) are the discrete versions of matrices \(A\), \(B\), and \(\mathbf{d}\) derived above. The vehicle's states and inputs at each time instant \(t = k T_s\), \(k = 0, 1, 2, \dots\), are represented by \(\mathbf{x}(k)\) and \(\mathbf{u}(k)\), respectively, where \(T_s\) represents the sampling interval used for discretizing the continuous-time system described by Equation~\ref{eq:linear_state_space}. The system here is assumed to have full-state feedback, implying that all states can be measured.

\section{Appendix C: MPC Parameters Analysis}
Figure \ref{fig:mpc_sensitivity_1} presents the results of parametric sensitivity analysis on different parameters of the MPC controller. The analysis systematically examines the influence of weight parameters on trajectory tracking performance, control smoothness, and dynamic response characteristics across multiple scenarios. These results provide guidance for parameter calibration and demonstrate the trade-offs between different performance objectives in the proposed control framework.

\begin{landscape}
\begin{figure}[h]
    \centering
    \begin{tabular}{>{\centering\arraybackslash}m{1cm} c c c c}
        % & Column 1 & Column 2 & Column 3 & Column 4 \\
        \(W_{\perp}\) & 
        \begin{subfigure}{0.3\textwidth}\label{mpc_sensitivity_a}
            \includegraphics[width=\linewidth]{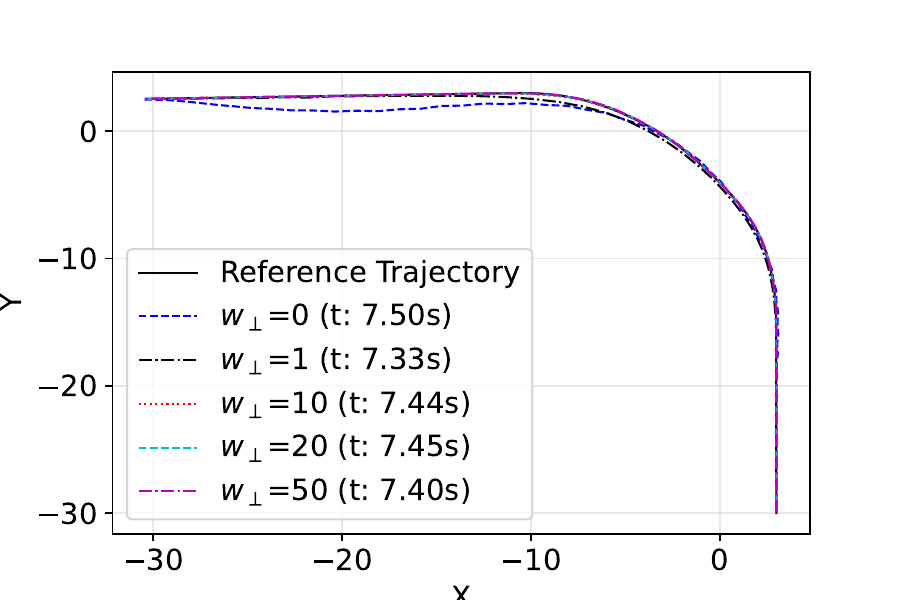}
            \caption{}
        \end{subfigure} & 
        \begin{subfigure}{0.3\textwidth}\label{mpc_sensitivity_b}
            \includegraphics[width=\linewidth]{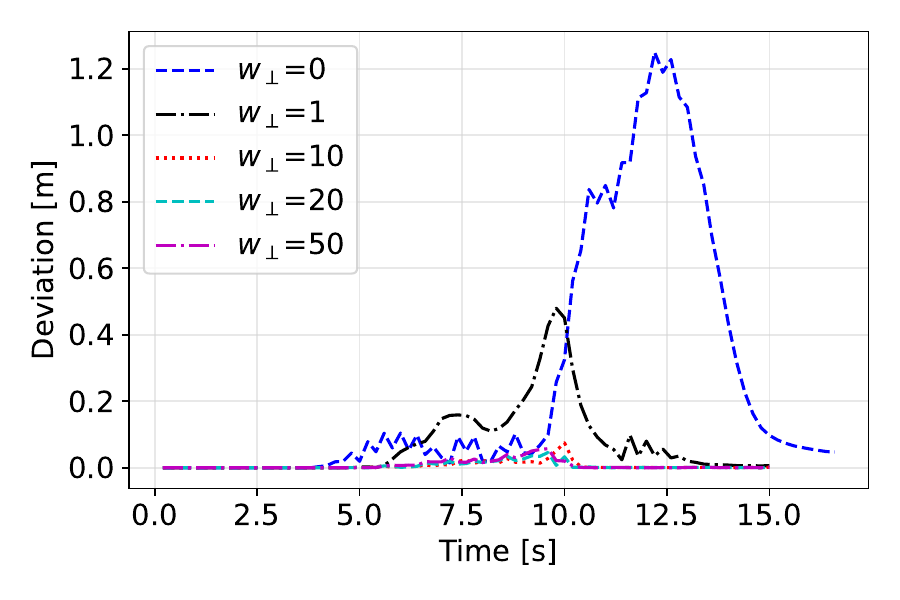}
            \caption{}
        \end{subfigure} & 
        \begin{subfigure}{0.3\textwidth}\label{mpc_sensitivity_c}
            \includegraphics[width=\linewidth]{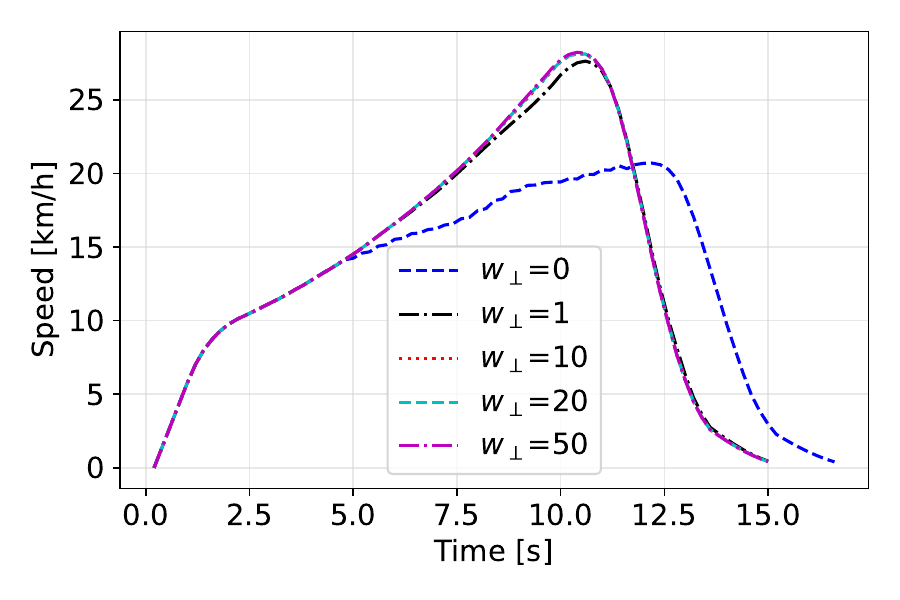}
            \caption{}
        \end{subfigure} & 
        \begin{subfigure}{0.3\textwidth}\label{mpc_sensitivity_d}
            \includegraphics[width=\linewidth]{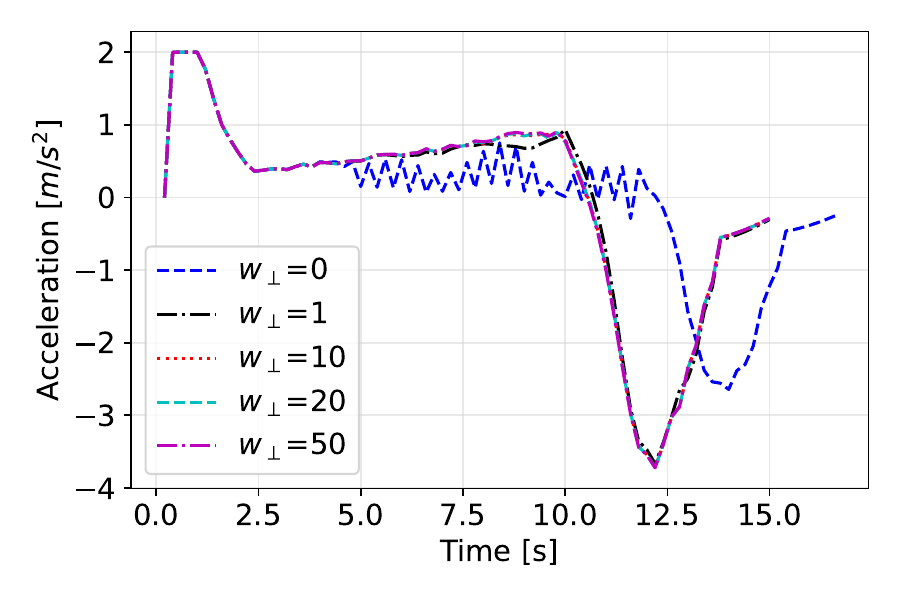}
            \caption{}
        \end{subfigure} \\
        \(W_{\parallel}\) & 
        \begin{subfigure}{0.3\textwidth}\label{mpc_sensitivity_e}
            \includegraphics[width=\linewidth]{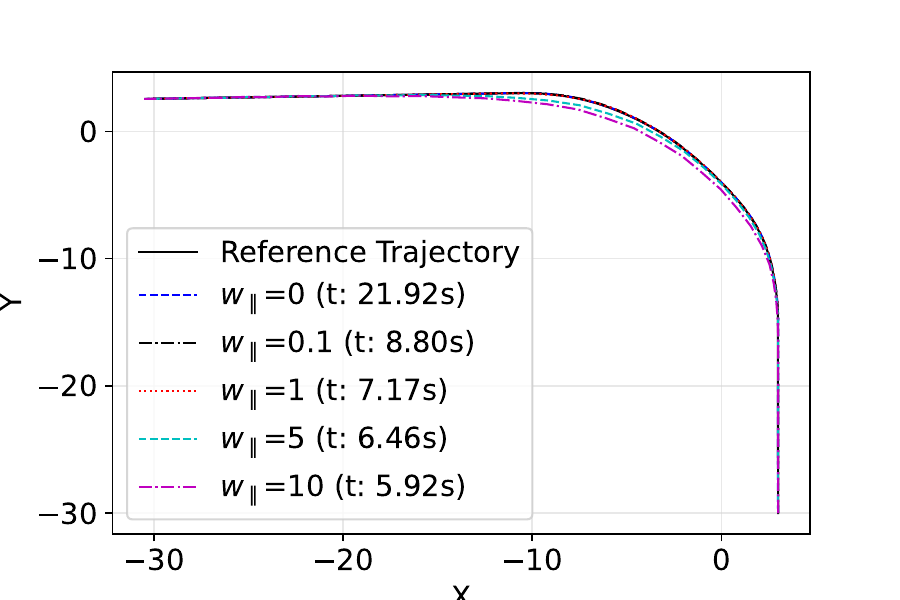}
            \caption{}
        \end{subfigure} & 
        \begin{subfigure}{0.3\textwidth}\label{mpc_sensitivity_f}
            \includegraphics[width=\linewidth]{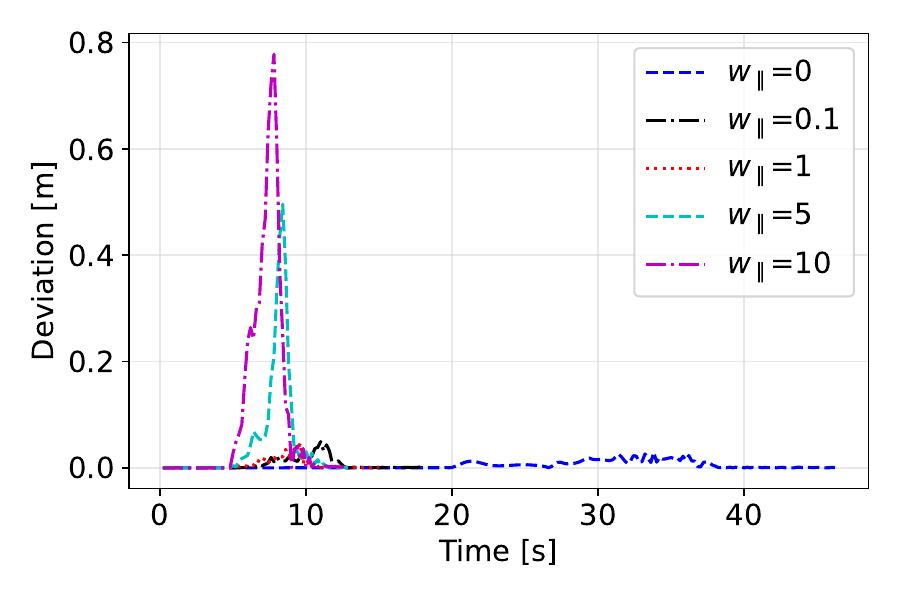}
            \caption{}
        \end{subfigure} & 
        \begin{subfigure}{0.3\textwidth}\label{mpc_sensitivity_g}
            \includegraphics[width=\linewidth]{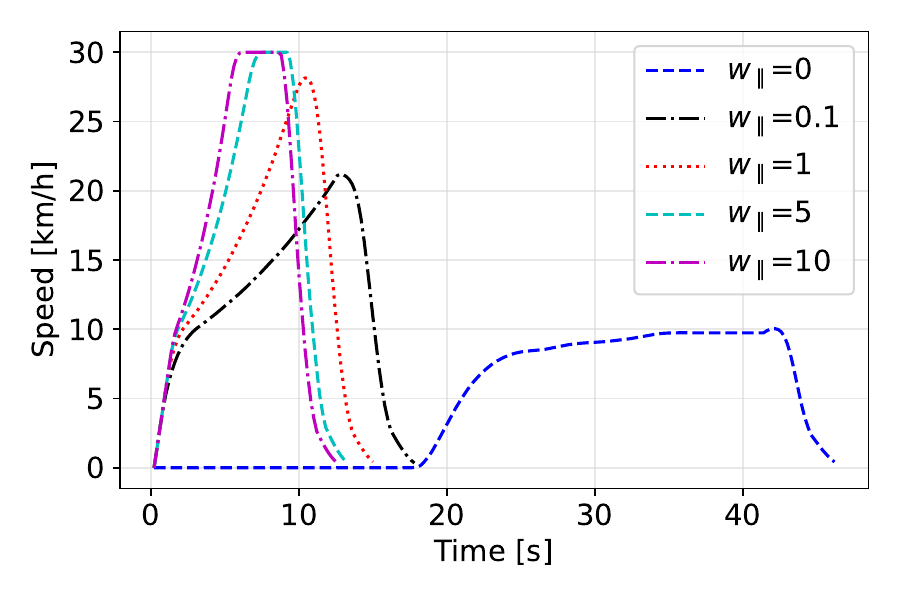}
            \caption{}
        \end{subfigure} & 
        \begin{subfigure}{0.3\textwidth}\label{mpc_sensitivity_h}
            \includegraphics[width=\linewidth]{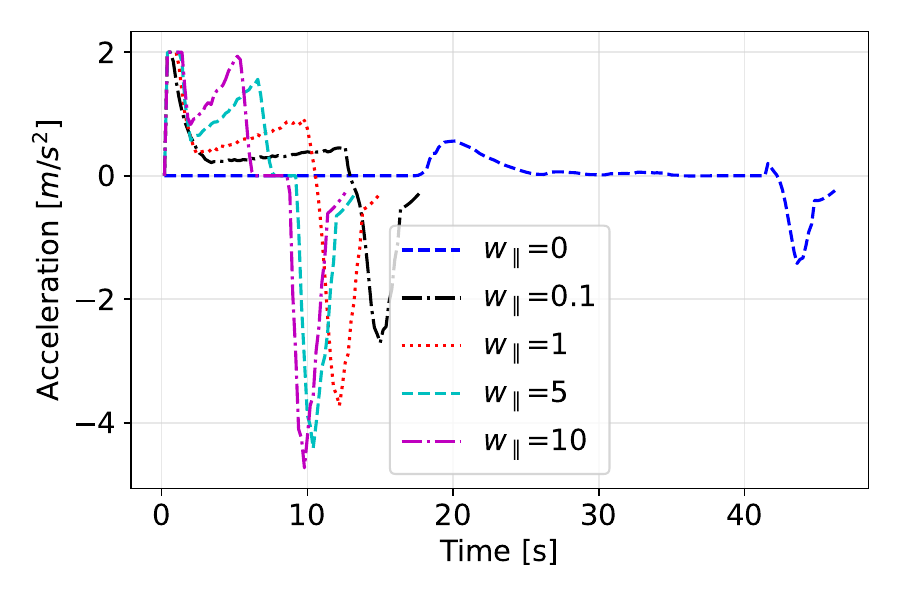}
            \caption{}
        \end{subfigure} \\
        \(R_{acc}\) & 
        \begin{subfigure}{0.3\textwidth}\label{mpc_sensitivity_i}
            \includegraphics[width=\linewidth]{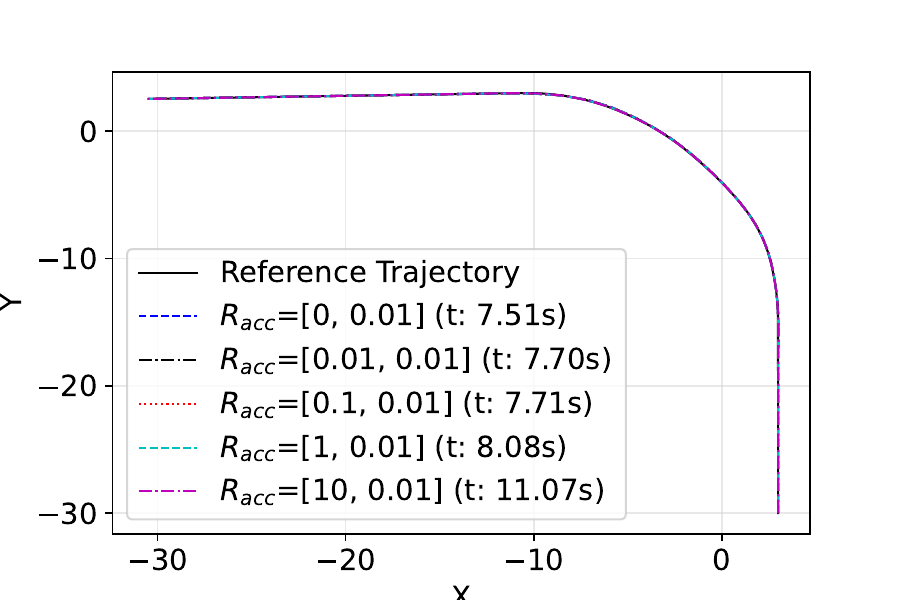}
            \caption{}
        \end{subfigure} & 
        \begin{subfigure}{0.3\textwidth}\label{mpc_sensitivity_j}
            \includegraphics[width=\linewidth]{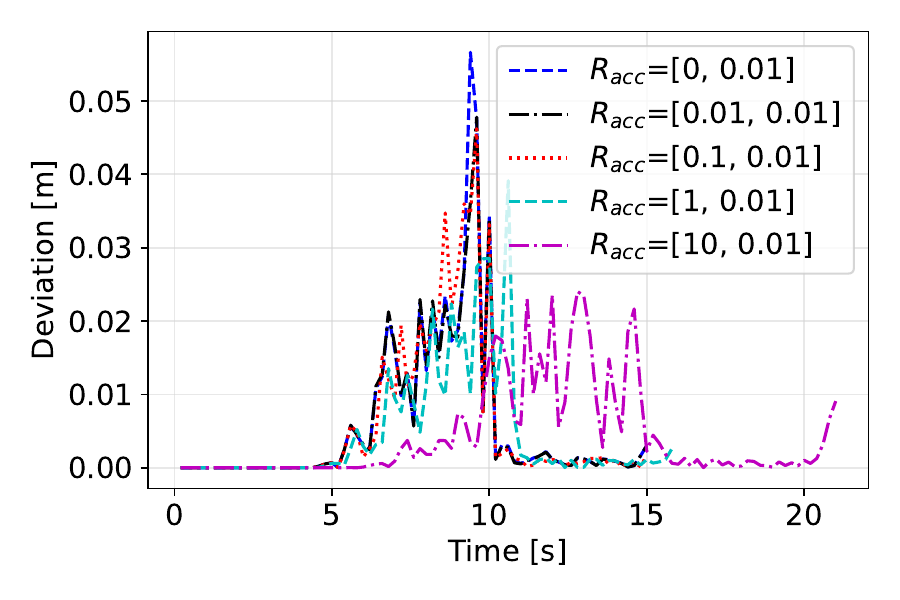}
            \caption{}
        \end{subfigure} & 
        \begin{subfigure}{0.3\textwidth}\label{mpc_sensitivity_k}
            \includegraphics[width=\linewidth]{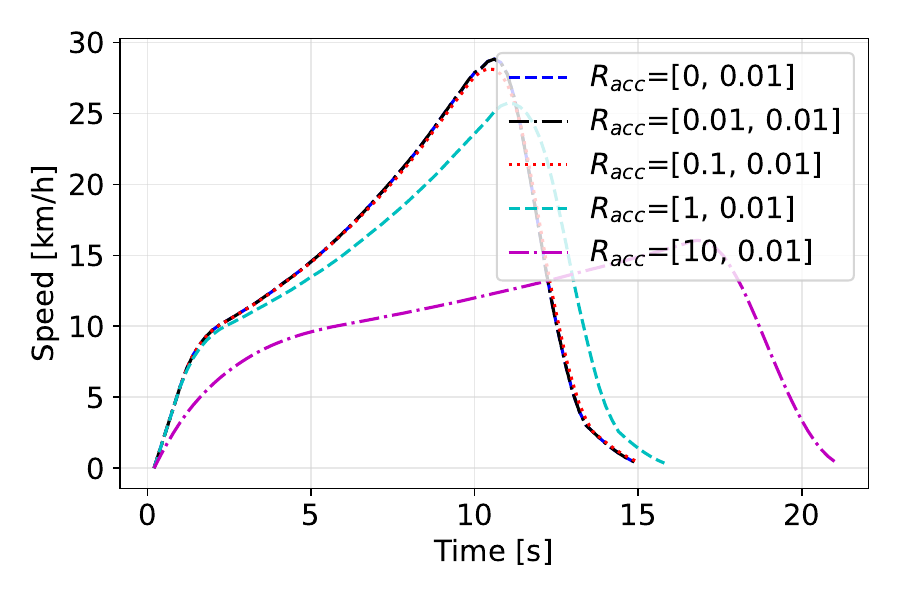}
            \caption{}
        \end{subfigure} & 
        \begin{subfigure}{0.3\textwidth}
            \includegraphics[width=\linewidth]{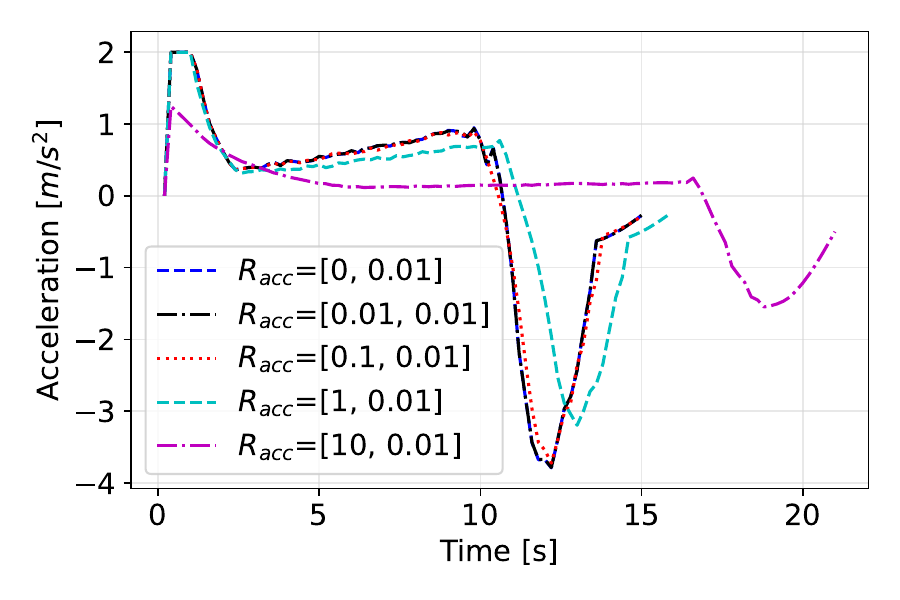}
            \caption{}
            \label{mpc_sensitivity_l}
        \end{subfigure} \\
        & Trajectory Following & Deviation from Reference & Speed Dynamics & Acceleration Dynamics \\
    \end{tabular}
    \caption{Sensitivity analysis of controller parameters}
    \label{fig:mpc_sensitivity_1}
\end{figure}

\begin{figure}[t]\ContinuedFloat
    \centering
    \begin{tabular}{>{\centering\arraybackslash}m{1cm} c c c c}
        % & Column 1 & Column 2 & Column 3 & Column 4 \\
        \(R_{steer}\) & 
        \begin{subfigure}{0.3\textwidth}\label{mpc_sensitivity_m}
            \includegraphics[width=\linewidth]{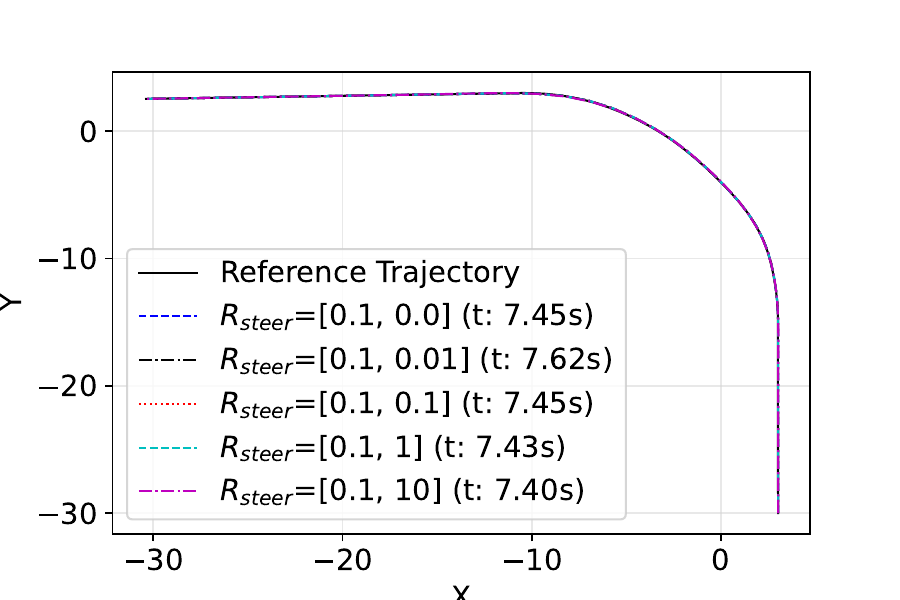}
            \caption{}
        \end{subfigure} & 
        \begin{subfigure}{0.3\textwidth}\label{mpc_sensitivity_n}
            \includegraphics[width=\linewidth]{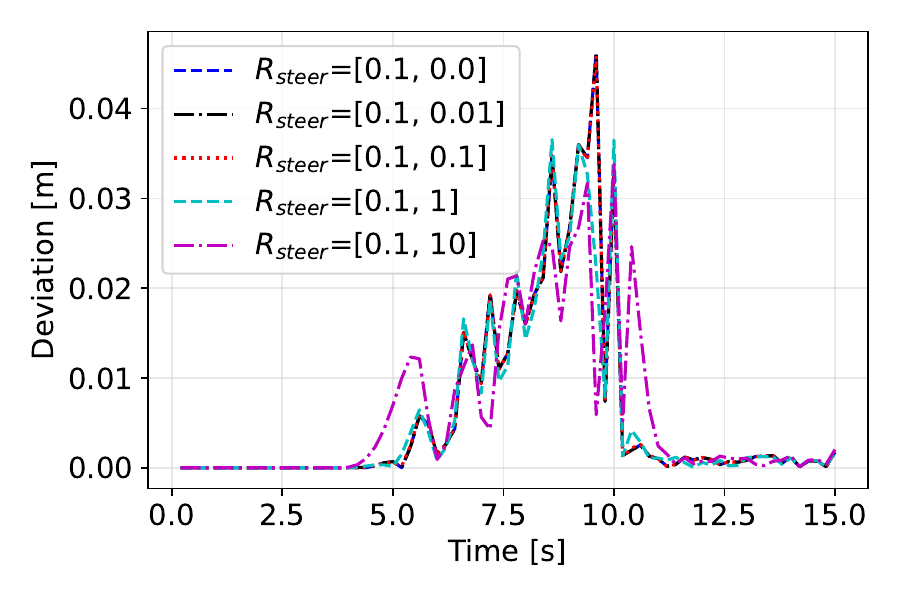}
            \caption{}
        \end{subfigure} & 
        \begin{subfigure}{0.3\textwidth}\label{mpc_sensitivity_o}
            \includegraphics[width=\linewidth]{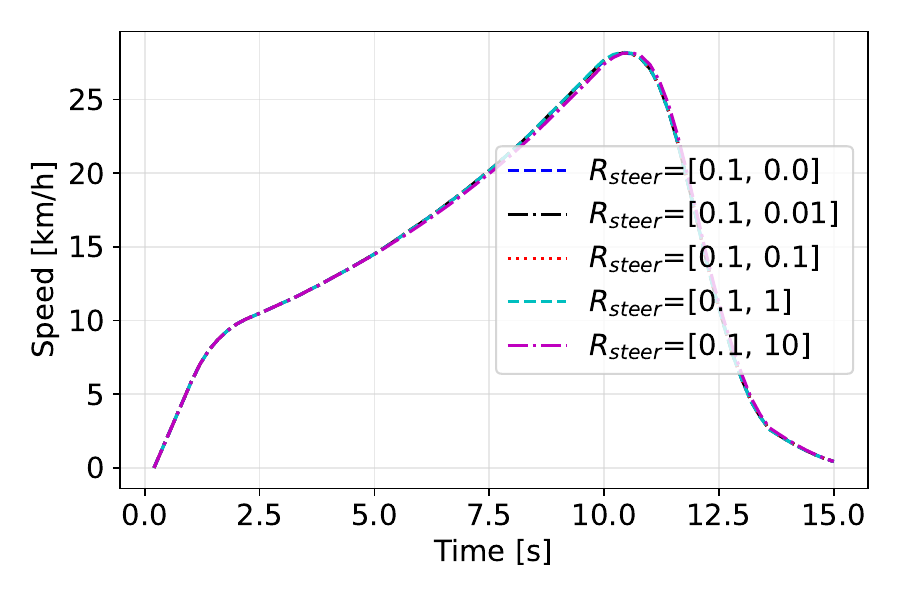}
            \caption{}
        \end{subfigure} & 
        \begin{subfigure}{0.3\textwidth}\label{mpc_sensitivity_p}
            \includegraphics[width=\linewidth]{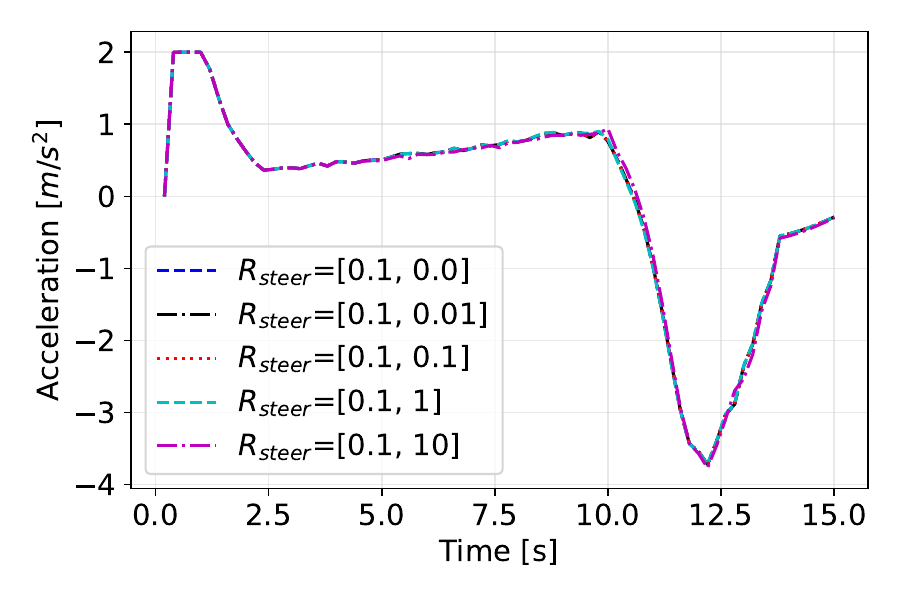}
            \caption{}
        \end{subfigure} \\
        \(Rd_{acc}\) & 
        \begin{subfigure}{0.3\textwidth}\label{mpc_sensitivity_q}
            \includegraphics[width=\linewidth]{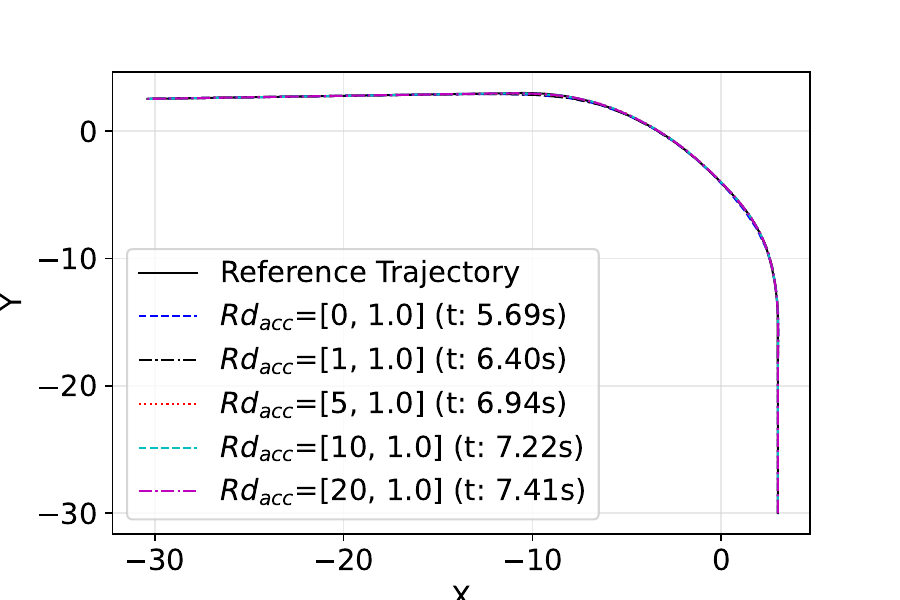}
            \caption{}
        \end{subfigure} & 
        \begin{subfigure}{0.3\textwidth}\label{mpc_sensitivity_r}
            \includegraphics[width=\linewidth]{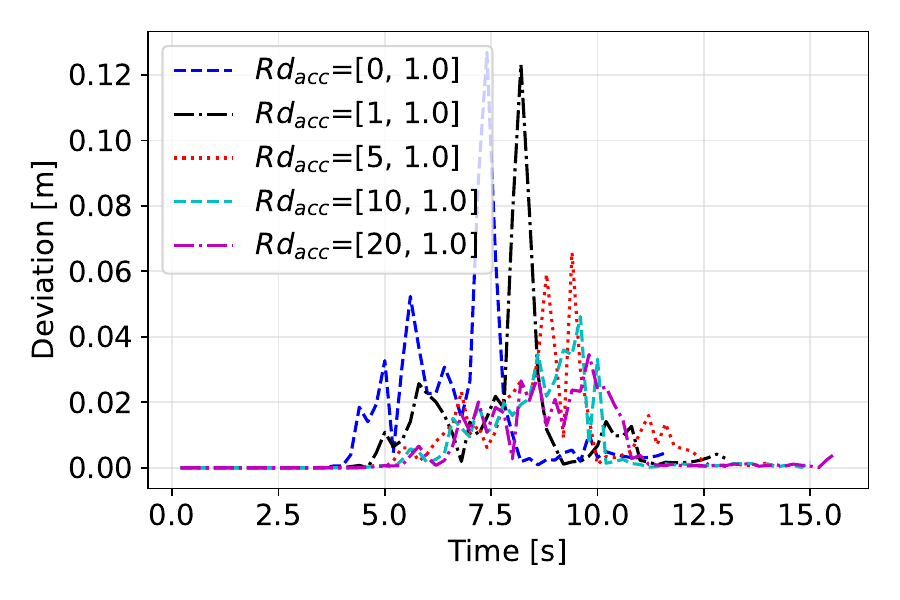}
            \caption{}
        \end{subfigure} & 
        \begin{subfigure}{0.3\textwidth}\label{mpc_sensitivity_s}
            \includegraphics[width=\linewidth]{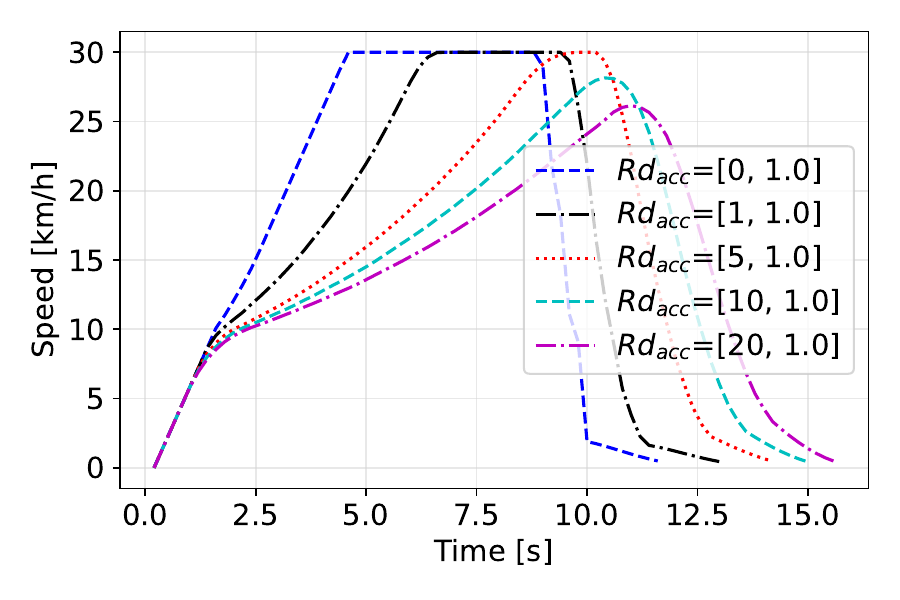}
            \caption{}
        \end{subfigure} & 
        \begin{subfigure}{0.3\textwidth}\label{mpc_sensitivity_t}
            \includegraphics[width=\linewidth]{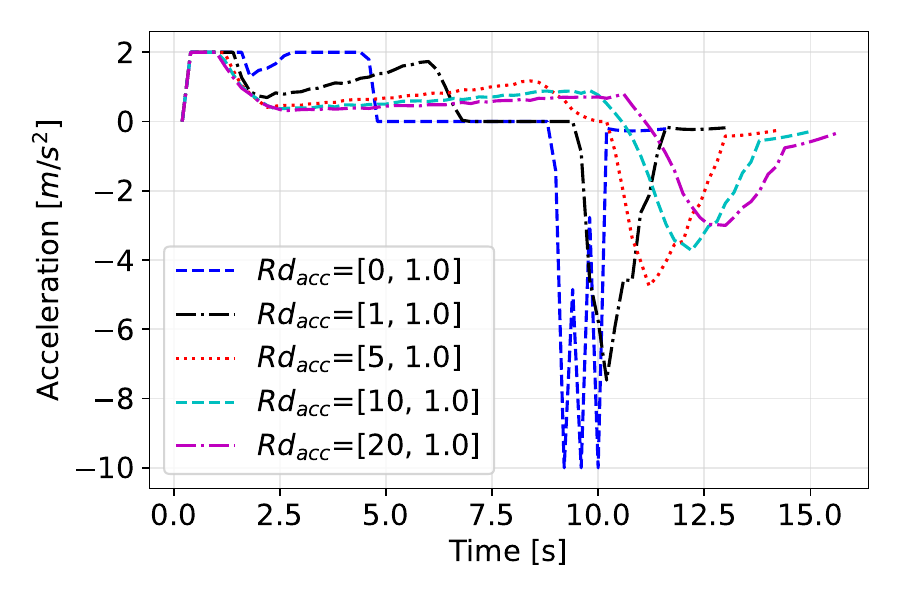}
            \caption{}
        \end{subfigure} \\
        \(Rd_{steer}\) & 
        \begin{subfigure}{0.3\textwidth}\label{mpc_sensitivity_u}
            \includegraphics[width=\linewidth]{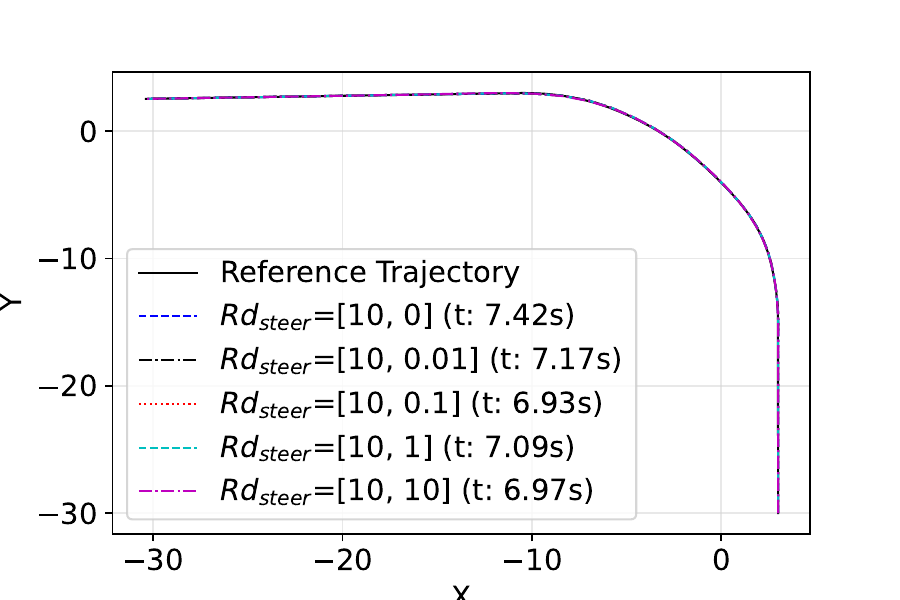}
            \caption{}
        \end{subfigure} & 
        \begin{subfigure}{0.3\textwidth}\label{mpc_sensitivity_v}
            \includegraphics[width=\linewidth]{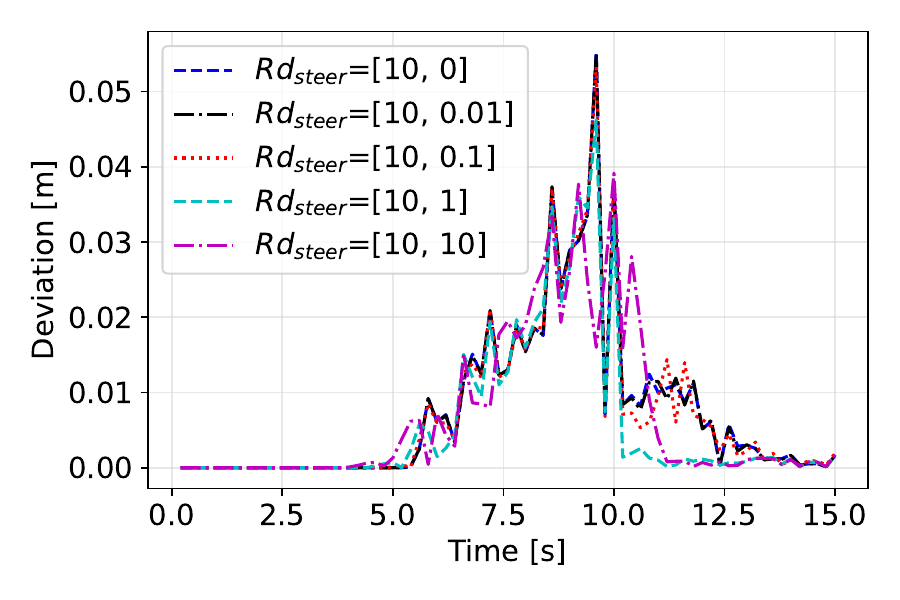}
            \caption{}
        \end{subfigure} & 
        \begin{subfigure}{0.3\textwidth}\label{mpc_sensitivity_w}
            \includegraphics[width=\linewidth]{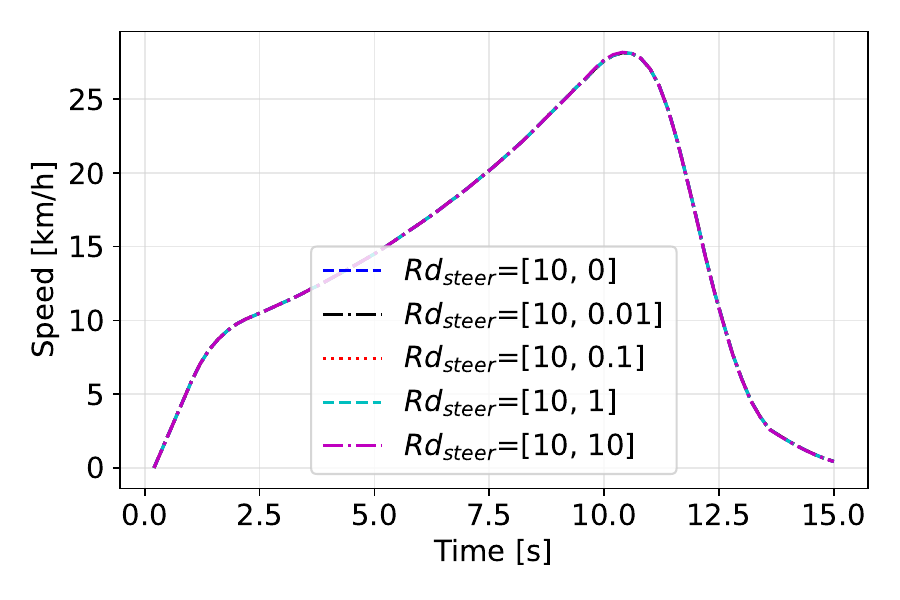}
            \caption{}
        \end{subfigure} & 
        \begin{subfigure}{0.3\textwidth}\label{mpc_sensitivity_x}
            \includegraphics[width=\linewidth]{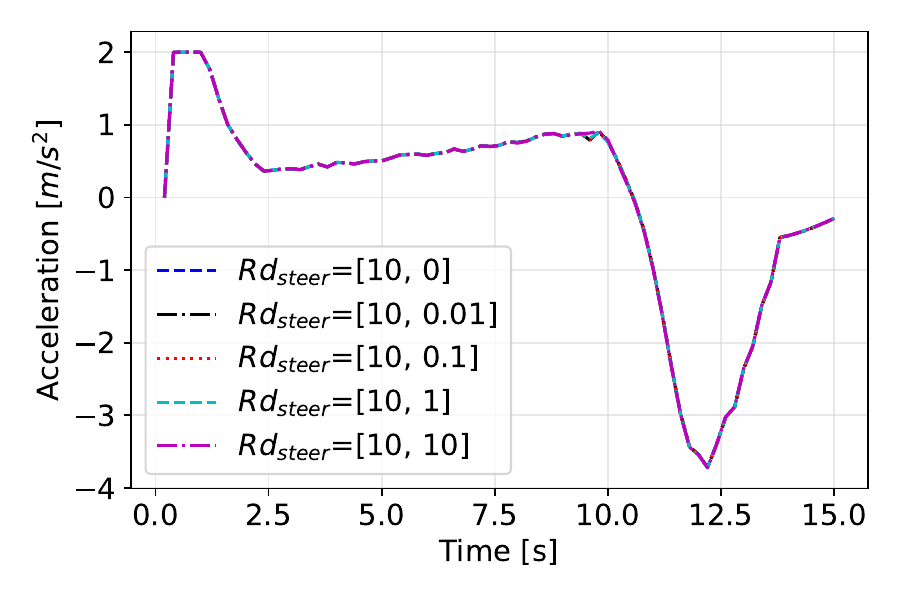}
            \caption{}
        \end{subfigure} \\
         & Trajectory Following & Deviation from Reference & Speed Dynamics & Acceleration Dynamics \\
    \end{tabular}
    \caption{Sensitivity analysis of controller parameters}
    \label{fig:mpc_sensitivity_2}
\end{figure}

\end{landscape}

\end{document}